\documentclass[a4paper,reqno]{amsart}
\usepackage[all]{xy}
\usepackage{color}
\usepackage{mathtools}
\usepackage{mathrsfs}
\usepackage{hyperref}
\usepackage{cleveref}
\usepackage{enumerate}
\usepackage{amsmath}
\usepackage{amsthm}
\usepackage{amssymb}
\usepackage{amscd}
\usepackage{graphicx}
\usepackage{epsfig}
\usepackage[english]{babel}
\usepackage[stable]{footmisc}
\usepackage[utf8]{inputenc}
\usepackage[T1]{fontenc}

\newtheorem{theorem}{Theorem}[section]
\newtheorem{lemma}{Lemma}[section]
\newtheorem{proposition}{Proposition}[section]
\newtheorem{corollary}{Corollary}[section]

\theoremstyle{definition}
\newtheorem{definition}{Definition}[section]
\newtheorem{example}{Example} [section]

\newtheorem{remark}{Remark}[section]

\begin{document}

\author{David Martínez Torres}

\address{Centro de An\'{a}lise Matem\'{a}tica, Geometria e Sistemas Din\^{a}micos,
Departamento de Matem\'atica, Instituto Superior T\'ecnico, Av. Rovisco Pais, 1049-001
Lisboa, Portugal}

\email{dfmtorres@gmail.com}

\title{The geometry of 2-calibrated manifolds}




\maketitle

%

\tableofcontents

\section{Introduction and statement of results}
In recent years there has been an enormous success in the study of symplectic
manifolds using approximately holomorphic methods. These methods -introduced by
S. Donaldson in $1996$ \cite{Do96}- amount to treating symplectic manifolds
 as generalizations of K\"ahler manifolds. To this end it is convenient to think of a symplectic manifold
  -once a compatible almost complex structure $J$ has been fixed-  as a K\"ahler manifold
   $(P,J,\Omega)$ for which the integrability condition for $J$ has been dropped.

Let $M$ be any hypersurface of the K\"ahler manifold $(P,J,\Omega)$. $M$ inherits
 on the one hand a codimension one distribution $D:=JTM\cap TM$ endowed with an integrable almost complex structure
   $J\colon D\rightarrow D$ (i.e. a CR structure of hypersurface type), and on the other hand a closed
   $2$-form $\omega:=\Omega_{\mid M}$ which is no-where degenerate when
    restricted to $D$. A 2-calibrated structure on $M$ -together with a compatible almost complex structure- is the
     structure obtained  when the integrability assumption on $J\colon D\rightarrow D$ is dropped.

Let us assume that the CR distribution of the $2n+1$-dimensional CR manifold
 (of hypersurface type)  $(M,D,J)$ is co-oriented (i.e the real line bundle $TM/D$ is trivial and a positive side has been chosen). The Levi-form is the symmetric tensor
\begin{eqnarray*}
\mathcal{L}\colon D\times D& \rightarrow & TM/D\\
(u,v)&\longrightarrow & [U,JV]/\sim\;\;,
\end{eqnarray*}
where $U, V$ are local sections of $D$ extending $u,v\in T_xM$, and we consider the class of the above Lie bracket at $x$ in the quotient real line bundle $TM/D$, where we can make sense of
positive and negative values. We can distinguish several interesting geometries according to the behavior of the Levi-form:

\begin{enumerate}
\item  If $\mathcal{L}$ is strictly positive (resp. negative) we get  a strictly
 pseudo-convex (resp. pseudo-concave) CR structure.  If we drop $J$ what remains is a co-oriented contact structure (these always carry almost complex structures along the contact distribution).
\item If $\mathcal{L}\equiv 0$ then  $D$ integrates into a codimension one foliation
 whose leaves inherit a K\"ahler structure. If $J$ is dropped what we obtain is a class of regular Poisson manifolds that include mapping tori associated
   to symplectomorphisms and more generally cosymplectic structures (defined by a closed 1-form $\alpha$
    and a closed 2-form $\omega$ such that $\alpha\wedge\omega^n$ is a volume form).
     When $n=1$ the latter are nothing but smooth taut foliations.
\item  If $n=1$ and $\mathcal{L}\geq 0$, by dropping $J$ we obtain a class of structures that include
 all taut confoliations (see section 3.5 in \cite{ET99}).
\end{enumerate}

\begin{definition}\label{def:2calstrdef} A 2-calibrated structure on $M^{2n+1}$ is a pair $(D,\omega)$, where $D$ is
 a codimension one distribution and $\omega$ a closed $2$-form no-where degenerate on $D$.

We call the triple $(M,D,\omega)$ a 2-calibrated manifold. We also say that $\omega$
 is positive on $D$. If $D$ is integrable we speak of 2-calibrated foliations.

$(M,D,\omega)$ is said to be integral if $[\omega]\in H^2(M;\mathbb{R})$ is in
the image of the integer cohomology, in which case we choose a lift $h\in H^2(M;\mathbb{Z})$ of
$[\omega]$ that we fix once and for all. The pre-quantum
 line bundle $(L,\nabla)$ is the unique  -up to isomorphism-   Hermitian line
  bundle with compatible connection with Chern class $h$ and curvature $-2\pi i\omega$.
\end{definition}

As we saw 2-calibrated structures do contain contact structures, cosymplectic structures and 3-dimensional taut confoliations.

A 2-calibrated manifold $(M,D,\omega)$ always admits compatible almost complex structures $J\colon D\rightarrow D$. The purpose of this paper is to explore how to adapt approximately holomorphic geometry to the tuple $(M,D,\omega,J)$, and see how we can apply this theory to know more about $(M,D,\omega)$.

In what follows all our manifolds will be closed and smooth, and all tensors and maps smooth unless otherwise stated.

The first application we will obtain is an analog of the existence of transverse
 cycles through any point of a 3-dimensional taut foliation.

The appropriate generalization of a transverse cycle is as follows:
\begin{definition}\label{def:submanifold}
$W$ is a 2-calibrated submanifold of $(M,D,\omega)$ if $TW\cap D$
has codimension one inside $TW$ and  $\omega$ is positive when
restricted to it. In other words, $W$ must intersect $D$
transversely and in a symplectic sub-distribution of $(D,\omega)$. \end{definition}

The existence of submanifolds -which extends the main result for contact manifolds in \cite{IMP01}- is the content of the following

 \begin{proposition}\label{pro:proposition1}
 Let $(M^{2n+1},D,\omega)$ be an integral 2-calibrated manifold and $L^{\otimes k}$ the sequence of powers of its pre-quantum line bundle (definition \ref{def:2calstrdef}). For any
  fixed point $y\in M$, any $m=1,\dots,n$, and any rank $m$ complex vector bundle
   $E\rightarrow M$, it is possible to find for all  $k\in \mathbb{N}$ large enough
   2-calibrated submanifolds $W_k$ of $M$ of codimension $2m$ through $y$ with the following properties:
 \begin{itemize}
 \item The inclusion $l\colon W_k\hookrightarrow M$ induces maps
 $l_{*}\colon \pi_j(W_k)\rightarrow \pi_j(M)$ which are isomorphisms for $j=0,\dots,n-m-1$,
  and an epimorphism for $j=n-m$. The same result holds for the homology groups.
 \item The Poincar\'e dual of  $[W_k]$ is $c_m(E\otimes L^{\otimes k}).$
 \end{itemize}
 \end{proposition}

The submanifolds in proposition \ref{pro:proposition1} are obtained  by pulling back the ${\bf 0}$ section of a vector bundle. Something similar
  can be done with the determinantal loci of a homomorphism of complex vector bundles
   (see theorem 1.6 in \cite{MPS02} and corollary 5.2 in \cite{Au01}).

\begin{proposition}\label{pro:determinantal} Let $(M,D,\omega)$ be an integral 2-calibrated
 manifold and $L^{\otimes k}$ the sequence of powers of its pre-quantum line bundle.
  Let $E$, $F$ be Hermitian vector bundles with connections and consider the sequence
  of bundles  $I_k=E^*\otimes F\otimes L^{\otimes k}=\mathrm{Hom}(E, F\otimes L^{\otimes k})$.
   Then for all $k\in \mathbb{N}$ large enough there exist  sections  $\tau_k$ of   $I_k$  for which the
    determinantal loci  $\Sigma^i(\tau_k)=\{x\in M |\; \mathrm{rank}(\tau_k(x))=i\}$
     are integral 2-calibrated submanifolds  stratifying  $M$.
\end{proposition}

The Poincar\'e Dual of the closure of $\Sigma^i(\tau_k)$ is given by the Porteous formula \cite{Po71}:

\[
 \Delta_{E,F\otimes L^{\otimes k},i} = \left| \begin{array}{cccc}
 c_{n-i} & c_{n-i+1} & \cdots & \\
 c_{n-i-1} & c_{n-i} & \cdots & \\
 & & \ddots & \\
 c_{n-m+1}& & \cdots & c_{n-i}
 \end{array}
 \right|,
\]
where $\mathrm{rank}E=m$, $\mathrm{rank}F=n$, and $c_j$ is the $j$-th Chern class
$c_j(F\otimes  L^{\otimes k}-E)$ defined by the equality
\begin{eqnarray*}
& 1+c_1(F\otimes  L^{\otimes k}-E)+c_2(F\otimes  L^{\otimes k}-E)+\cdots = &\\
& (1+c_1(F\otimes  L^{\otimes k})+c_2(F\otimes  L^{\otimes k})+\cdots)/(1+c_1(E)+c_2(E)+\cdots). &
\end{eqnarray*}
If the rank of $E$ and $F$, and $i$ are chosen so that $\Sigma^{i-1}(\tau_k)$ is empty, then
  $\Sigma^{i}(\tau_k)$ is a closed 2-calibrated submanifold.

 \begin{corollary} \label{cor:determinantal2}   Let $(M,\alpha)$, $\alpha\in \Omega^1(M)$, be an
  exact contact manifold of dimension $2n+1$. Let  $E,F$ be complex
vector bundles and let $i$ be a positive integer such that
\begin{itemize}
\item The codimension in $\mathrm{Hom}(E,F)$ of the strata of homomorphisms of rank
 $i$ is not bigger than $2n+1$.
 \item The  codimension in $\mathrm{Hom}(E,F)$ of the
  strata of homomorphisms of rank $i-1$ is  bigger than $2n+1$.
\end{itemize}
Then there exist contact submanifolds whose Poincar\'e dual is $\Delta_{E,F,i}$. In particular,
 for any even cohomology class which is a Chern class of some complex vector bundle over
  $M$, there exist a contact submanifold Poincar\'e dual to it.
\end{corollary}

\begin{remark}\label{rem:detvszero}  One is expecting that the determinantal submanifolds
 coming proposition \ref{pro:determinantal} will be more general than the zeroes of vector bundles
 coming from proposition \ref{pro:proposition1}. A more detailed discussion of this issue appears in  \ref{sec:det}.
\end{remark}

The next application is an analog for $2$-calibrated manifolds of the embedding theorem for
 symplectic manifolds of \cite{MPS02} (theorem 1.2), extending results of  \cite{MPS01} for contact manifolds.

 \begin{corollary}\label{cor:embedding} Let  $(M^{2n+1},D,\omega)$ be an integral 2-calibrated
  manifold.  Then  it is possible to find  maps  $\phi_k\colon M\rightarrow \mathbb{C}\mathbb{P}^{2n}$
   so that for all $k\in \mathbb{N}$ large enough one has:

\begin{itemize}
\item ${d\phi_k}_{\mid D}$ is injective ($\phi_k$ is an immersion along $D$).
\item $[\phi_k^*\omega_{FS}]=[k\omega]$, where $\omega_{FS}$ is the  Fubini-Study 2-form of $\mathbb{C}\mathbb{P}^{2n}$.
\end{itemize}

In particular if $(M^3,D)$ is a 3-manifold with a (smooth) taut confoliation,  it
 is possible to find immersions along $D$ in $\mathbb{C}\mathbb{P}^2$.
\end{corollary}

The previous corollary can be improved in two directions:

\begin{corollary}\label{cor:Tischler} (see \cite{MPS02}, corollary 2.6) Let  $(M^{2n+1},\mathcal{D},\omega)$
 be  a manifold with an integral  2-calibrated foliation. Then the maps of corollary
  \ref{cor:embedding} can be composed from the right with  diffeomorphisms of $M$, so that for all
    $k\in \mathbb{N}$ large enough  the equality  $[\phi_k^*\omega_{FS}]=[k\omega]$ holds also at the level of foliated
    $2$-forms, i.e. ${\phi_k^*\omega_{FS}}_{\mid\mathcal{D}}={kw}_{\mid \mathcal{D}}$.
\end{corollary}

The second improvement is that the immersion along $D$ can be perturbed to be transverse to any
finite collection of complex submanifolds of  projective space.

Another application is the existence of  Lefschetz pencil structures, introduced in  \cite{IM04}.

\begin{definition}(see section 1 in \cite{Do99}) Let $(M,D,\omega)$ be a 2-calibrated manifold and $x\in M$.
 A chart $\varphi\colon (\mathbb{C}^n\times\mathbb{R},0)\rightarrow (M,x)$ is compatible
with $(D,\omega)$ (at $x$) if at the origin it sends the foliation of $\mathbb{C}^n\times\mathbb{R}$
 by complex hyperplanes into $D$, and $\varphi^*\omega(0)$ restricted to the subspace
   $\mathbb{C}^n\times\{0\}$  is of type (1,1).
\end{definition}

\begin{definition}\label{def:pencil} (see \cite{Pr02}) A Lefschetz pencil structure for $(M,D,\omega)$
 is a triple $(f,B,\Delta)$ where $B\subset M$ is a codimension four 2-calibrated submanifold,
  and $f\colon M\backslash B\rightarrow \mathbb{C}\mathbb{P}^1$ is a smooth map such that:
\begin{enumerate}
\item $f$ is a submersion along $D$ away from $\Delta$, a 1-dimensional manifold transverse to
 $D$ where the restriction of the differential of $f$ to $D$ vanishes.
\item For any $x\in \Delta$ there exist a chart $\varphi$ compatible with
 $(D,\omega)$ at $x$ and a complex coordinate $\zeta$ of $\mathbb{C}\mathbb{P}^1$ defined about $f(x)$, such that
\[
\zeta\circ f\circ\varphi (z,s)={(z^1)}^2+\cdots +{(z^n)}^2+t(s),
\]
 where $t\in C^\infty(\mathbb{R},\mathbb{C})$.
\item  For any $x\in B$  exist a chart $\varphi$ compatible with
 $(D,\omega)$ at $x$ and a complex coordinate $\zeta$ of $\mathbb{C}\mathbb{P}^1$ defined about $f(x)$, such that
  $B\equiv z^1=z^2=0$ and  $\zeta\circ f\circ\varphi(z,s)=z^1/z^2$.
\item $f(\Delta)$ is an immersed curve with generic self intersections.
\end{enumerate}
\end{definition}

\begin{theorem}\label{thm:pencil} Let $(M,D,\omega)$ be an integral 2-calibrated manifold
  and let $h$ be an integer lift of $[\omega]$.  Then for all  $k\in \mathbb{N}$ large enough there exist Lefschetz
   pencils $(f_k,B_k,\Delta_k)$ such that:
\begin{enumerate}
\item The regular fibers are Poincar\'e dual to $kh$.
\item The inclusion $l\colon W_k\hookrightarrow M$ induces maps $l_*\colon \pi_j(W_k)\rightarrow \pi_j(M)$
(resp. \newline $l_*\colon H_j(W_k;\mathbb{Z})\rightarrow H_j(M;\mathbb{Z})$) which are isomorphisms
for $j\leq n-2$ and an epimorphism for $j=n-1$.
\end{enumerate}
\end{theorem}

All the stated results follow mostly from a general principle  of (estimated) transversality
 along $D$ (theorems \ref{thm:main1} and \ref{thm:main2}).

 In a  problem $\mathcal{P}$ of transversality along $D$ we have three ingredients:
  (i) the bundle $E\rightarrow (M,D,\omega)$, (ii) the submanifold or more generally the
   stratification $\mathcal{S}\subset E$, and (iii) the section $\tau\colon M\rightarrow E$
    to be perturbed to become transverse along $D$ to $\mathcal{S}$.

In section \ref{sec:ampleb} we will define the class of sections and bundles we will work with,
  the so called sequences of very ample bundles and approximately holomorphic  sections.

 As in the approximately holomorphic  theory  for symplectic manifolds (see \cite{Do96,Au01}),
  transversality problems will be  solved  by patching local solutions. The right strategy
  to solve the corresponding local problems for sections is to turn them into  local problems
   for approximately holomorphic functions. This will be done through the use of reference sections,
    which can be thought of as the bump functions of the theory. The necessary local analysis
     needed to construct such sections is developed in section \ref{sec:locah}.

There is a second strategy to solve  $\mathcal{P}$. It is not only true that the natural example
 of a  2-calibrated structure is a hypersurface inside a symplectic manifold, but every
  $2$-calibrated manifold ($D$ co-oriented) admits a symplectization $(M\times [-\epsilon, \epsilon],\Omega)$
   (lemma \ref{lem:sympl}). We will introduce a new transversality problem $\bar{\mathcal{P}}$
    for a stratification $\bar{\mathcal{S}}$ of a bundle
     $\bar{E}\rightarrow (M\times [-\epsilon, \epsilon],\Omega)$, so that a solution
$\bar{\tau}\colon M\times [-\epsilon,\epsilon]\rightarrow \bar{E}$
to $\bar{\mathcal{P}}$ restricts to
$\bar{\tau}_{\mid M}$ a solution to $\mathcal{P}$.  The advantage of
this point of view is that since we are in a symplectic manifold, as
long as the extension $\bar{\mathcal{P}}$ falls in the right class
of problems we can use the existing approximately holomorphic theory
for symplectic manifolds to solve it. Still, the existing
approximately holomorphic theory turns out not to be enough for our
purposes, so we need to develop further the relative
approximately holomorphic theory introduced by J.P. Mohsen
\cite{Moh01}.  We will make an  exposition of both the intrinsic and
the relative approximately holomorphic theories, and we will prove
the main transversality theorem using the latter.

In section \ref{sec:esttrans} we give an account  of the notion of estimated transversality
  of a section along a distribution. For the intrinsic theory (problem $\mathcal{P}$)
   the distribution will be $D$, whereas for the relative theory the problem $\bar{\mathcal{P}}$
    will amount to achieving transversality over $M\subset (M\times[-\epsilon,\epsilon],\Omega)$.
     We will also introduce the right class of stratifications $\mathcal{S}$ (already defined
      in the symplectic setting in \cite{Au01}), the so called approximately holomorphic
       finite Whitney stratifications, whose  strata  roughly behave as the zero section of
        a vector  bundle in the sense that locally  they will be given by approximately holomorphic
         functions and they will be transverse enough to the fibers. The fundamental technical result
         (lemma \ref{lem:localchartrans}) will be that locally estimated transversality
         along $D$ (resp. over $M$) of an approximately holomorphic  section to $\mathcal{S}$
          (resp. $\bar{\mathcal{S}}$), will be equivalent to estimated transversality
          along $D$ (resp. over $M$) to ${\bf 0}$ of a related $\mathbb{C}^l$-valued approximately holomorphic function.

Section \ref{sec:pholjets} is devoted to the study of bundles of pseudo-holomorphic jets,
 needed to obtain what we call generic  approximately holomorphic maps to projective spaces,
  constructed by projectivizing  ($m+1$)-tuples of  approximately holomorphic sections
  of powers of the pre-quantum line bundle $L^{\otimes k}$  (i.e. analogs of generic
  linear systems in complex geometry); genericity will be defined as the solution of a uniform
  strong transversality problem to a stratification $\mathcal{S}$ in these
   bundles of pseudo-holomorphic jets. Several difficulties
    have to be overcome. Firstly, since we want to obtain a strong transversality result
     the jet of the section to be perturbed has to be itself an approximately holomorphic
      section, so that the transversality problem falls in the right class, something which
       fails to hold due to the uniform positivity along $D$ of the sequence $L^{\otimes k}$.
         This is solved by introducing a new connection in the bundles of pseudo-holomorphic
          jets. Secondly, we need to define a stratification
           $\mathcal{S}$ of the right kind. This is done in
            section \ref{TBA} by  introducing the  bundles of pseudo-holomorphic jets
             for maps to projective spaces, and defining there $\mathbb{P}\mathcal{S}$
               -a ``linear'' analog of the Thom-Boardman stratification-; $\mathcal{S}$
                is then constructed  by pulling back  $\mathbb{P}\mathcal{S}$ by the corresponding
                 jet extension of the projectivization map
                 $\pi\colon \mathbb{C}^{m+1}\backslash\{0\}\rightarrow \mathbb{C}\mathbb{P}^m$.
                  The properties of both the map and of  $\mathbb{P}\mathcal{S}$ are used
                   to conclude that $\mathcal{S}$ is indeed of the right kind, and thus
                    the transversality problem falls in the right class.
                     The necessary modifications for the relative theory are also described.

In section \ref{sec:mainthm} we give the main strong transversality result.

The proofs of the theorems stated in this introduction are given in section \ref{sec:applications}.

Our results are based on the existence of plenty of approximately holomorphic sections of very ample
 line bundles.  In the integrable setting the existence of enough meromorphic functions/holomorphic sections has been used to prove results of similar nature to ours:

\begin{enumerate}
\item[(i)] In  \cite{Gh99} E. Ghys gave conditions on a compact space laminated by Riemann
 surfaces for the existence of plenty of meromorphic functions. More generally, B. Deroin
  has extended those results to laminations by complex leaves without vanishing cycle, and
   endowed with positive Hermitian line bundles \cite{De03}. The  work of Ghys and Deroin
    proves the existence of leafwise holomorphic embeddings into projective spaces of the
     aforementioned laminated spaces (compare with corollary \ref{cor:embedding}), though
      the maps -even in the case of smooth foliations- are in general only continuous
       in the transverse directions. The strategy they follow is working in the universal
        cover of the leaves of the lamination. Interestingly enough, Deroin's results are
         obtained by extending some techniques of approximately holomorphic geometry to
         the leaves,  which are open K\"ahler manifolds with bounded geometry.

\item[(ii)] In \cite{OS00}  Ohsawa and Sibony gave a solution to the tangential Cauchy-Riemann equation with $L^2$-estimates for sections of a positive CR line bundle over a Levi-flat
  compact manifold. As a consequence they were able to produce CR embeddings into projective
   space of any prescribed order of regularity (though in general non-smooth).
\end{enumerate}

Part of the results of the present paper were announced in \cite{IM03,IM04}
(proposition \ref{pro:proposition1}, corollary \ref{cor:embedding}, corollary
 \ref{cor:Tischler}, theorem \ref{thm:pencil} and  theorem \ref{thm:main1}),
 where an account of the results available through an intrinsic approximately holomorphic theory was presented.

While a more detailed study of  2-calibrated structures is feasible, we do not think
 the results that could be obtained would be relevant enough to justify its undertaking.

There are two main reasons to develop an approximately holomorphic theory  for 2-calibrated structures.
The first one is because they contain  contact structures and 2-calibrated foliations.
 Approximately holomorphic geometry has already been introduced in the contact setting
  \cite{IMP01,Pr02,Moh01,MPS01}. Its most important application has been the construction of compatible  open book
   decompositions for contact manifolds of arbitrary dimension \cite{GM02}. Our contribution
    in this paper to contact geometry is the construction of a large class of contact submanifolds
     and the determination of their homology class (corollary \ref{cor:determinantal2}).
We want to propose 2-calibrated foliations as an interesting higher dimensional generalization
 of 3-dimensional taut foliations. In \cite{Ma05} -and building on the results of this paper- it is shown that any such foliation
  $(M,\mathcal{D},\omega)$ contains a 3-dimensional taut foliation
  $(W^3,\mathcal{D}_W)\hookrightarrow (M,\mathcal{D})$ so that the inclusion
   descends to a homeomorphism between  leaf spaces. This is done by showing that $W^3$
   can be chosen to intersect each leaf of $(M,\mathcal{D})$ in a unique connected
   component; this is somehow surprising since often the leaves are immersed submanifolds dense in $M$.

The second reason to develop an approximately holomorphic theory for 2-calibrated structures
 is that sometimes they appear as auxiliary structures. If $M$ is an odd dimensional manifold and $\omega$
  a maximally non-degenerate closed 2-form,
 any distribution $D$ complementary to the kernel of $\omega$ endows $M$ with a 2-calibrated
  structure. In \cite{MMP03} this idea was applied to almost contact manifolds to construct
   (via approximately holomorphic theory) open book decomposition with control on the
   topology of the leaves (see also \cite{Qu79}).

If $(M,D,J)\hookrightarrow (\mathbb{C}\mathbb{P}^N,\omega_{FS})$ is a CR manifold of hypersurface type which has a CR embedding in projective space, then in \cite{Ma05b} we show that the constructions of this paper can be performed in the CR category. In particular CR Lefschetz pencils are constructed, yielding CR Morse functions defined away from a CR submanifold of base points.

All the applications outlined so far for contact manifolds, 2-calibrated foliations, and
 projective CR manifolds  use at most pseudo-holomorphic 1-jets. If the CR manifold is Levi-flat then it
  makes sense to speak about $r$-generic CR
functions. These are defined to be leafwise $r$-generic holomorphic functions, i.e. functions whose leafwise holomorphic $r$-jet is transverse to the Thom-Boardman stratification of the bundle of holomorphic $r$-jets over each leaf. In \cite{Ma05b} we show that Levi-flat CR manifolds embedded in projective space admit  for all $k\gg 1$ $r$-generic linear systems. These are (holomorphic) linear systems of
 $\mathcal{O}(k)\rightarrow \mathbb{C}\mathbb{P}^N$ of rank $m(r)$ whose restriction to
 $M$ define $r$-generic CR functions away from base points  (definition \ref{def:rgenseccr}). Briefly,  such functions are
 easily seen to be CR functions whose CR $r$-jet prolongation solve $\mathcal{P}_{int}$ a transversality problem  over the leaves of the foliation $\mathcal{D}$ in the bundle of CR  $r$-jets of CR maps from $M$ to $\mathbb{C}\mathbb{P}^{m(r)}$.
   One has to show that it can be ``linearized" to a transversality problem $\mathcal{P}_{lin}$,
    (the bundle, the stratification, and the notion of CR $r$-jet all have to be replaced by
     ``linear" analogs) that fits into the ones solved in theorem \ref{thm:main2}; solutions are shown to exist among restrictions of holomorphic
       sections
       $\mathcal{O}(k)$. Finally, it has to be checked that the CR solution to  $\mathcal{P}_{lin}$
        is also a solution of $\mathcal{P}_{int}$.

We think that the existence of $r$-generic linear systems for projective Levi-flat CR manifolds
 is a relevant result  by itself and justifies the  extension of the approximately holomorphic
  theory to higher order jet bundles, which  is technically awkward. We expect it to be useful to analyze
   such manifolds. For example one can use it to define $r$-generic functions
   $f\colon (M^{2n+1},\mathcal{D},J)\rightarrow \mathbb{C}\mathbb{P}^n$ (with no base points)
    for which the regular level sets are unions of circles (with variable number of components),
     and using the analysis of the singularities define a dynamical system transverse to
     $\mathcal{D}$ (at least for low values of $n\geq 2$); by iterating the Lefschetz pencil
      construction (the dimensional induction of \cite{Au02b}, section 5) one can also define
       maps to $\mathbb{C}\mathbb{P}^{n-1}$ whose fibers (by \cite{Ma05}) are 3-manifolds
       intersecting each leaf of $\mathcal{D}$ in a connected Riemann surface.

       We point out that the results in  \cite{Ma05b} do not include those of Ghys and Deroin \cite{Gh99,De03} and those of Ohsawa and Sibony \cite{OS00}. Our results require starting with a CR embedding into projective space (\cite{OS00} gives sufficient conditions to produce it).

\subsection{Acknowledgements}
I wish to express my gratitude to Alberto Ibort for the numerous fruitful conversations that
 helped to improve many aspects of this work. It is also my pleasure to thank Denis Auroux, Vicente Mu\~noz,
  Fran Presas and Ignacio Sols for their helpful comments, and \'Etienne Ghys
   for kindly pointing me to very relevant literature on the subject. I am also grateful to the referees for their valuable comments and suggestions.

\section{Ample bundles and approximately holomorphic sections}\label{sec:ampleb}

Let $(M,D,\omega)$ be  an integral  2-calibrated manifold. Let us fix once and for
 all a compatible almost complex structure $J\colon D\rightarrow D$, and a metric
 $g$ so that $g_{\mid D}=\omega(\cdot,J)$. The kernel of $\omega$  is required to
  be $g$-orthogonal to $D$, so as to make some of the computations
   in the local theory simpler. Notice that for any such metric the closed $2n$-form
    $\omega^{n}$ is a calibration for $D$ \cite{HL82}.

If we forget about the 2-form what remains is the following structure.
\begin{definition}\label{def:almostcr} An almost CR structure is a tuple $(M,D,J,g)$ where
 $D$ is a codimension one distribution, $J\colon D\rightarrow D$  an almost complex
  structure, and $g$ a metric whose restriction to $D$ is compatible with $J$ ($J$
   is $g$-orthogonal and $g$-antisymmetric).
\end{definition}

Let $(L,\nabla)\rightarrow M$ be any Hermitian line bundle -or more generally vector
 bundle- with compatible connection. Let $\hat{D}$ denote  the pullback to $L$ of $D$;
  let  $\hat{J}$ and $\hat{g}$ be the almost complex structure and metric on $L$,
   which extend the Hermitian structure on the fibers and are  defined on the horizontal
    distribution associated to $\nabla$ by pulling back $J$ and $g$ respectively.
      Then $(L,\hat{D},\hat{J},\hat{g})$ is an almost CR manifold.

Our goal is to be able to construct sections $\tau\colon M\rightarrow L$ which (i)
 are close enough to satisfying  $\tau_*J=\hat{J}\tau_*$ (for which we
  use the adjective almost holomorphic instead of almost CR to be consistent with
   the terminology of \cite{IMP01} and \cite{Pr02}), and (ii)
 transverse to suitable submanifolds of the total space of $L$. In the almost complex setting we know that what ensures their existence is roughly speaking asking the curvature of the connection to be of type (1,1)
   and positive.

\begin{definition}\label{def:ampleb} (see \cite{Au01}, definition 2.1) Given $c>0$, $\delta\geq 0$,
 a Hermitian line bundle with compatible connection
 $(L,\nabla)\rightarrow (M,D,J,g)$ is $(c,\delta)$-ample (or just ample)
 if its curvature $F$ verifies  $iF(v,Jv)\geq
cg(v,v), \forall v\in D$, and $|F_{\mid D}-F^{1,1}_{\mid D}|_g\leq \delta $,
 where we use the supremum norm.

A sequence $(L_k,\nabla_k)$ of Hermitian line bundles with compatible connections
   is asymptotically very ample (or just very ample) if fixed constants
$c>0$, $\delta, (C_j)_{j \geq 0}\geq 0$  exist, so that for all $k \gg 1$ the following inequalities
 for the curvatures  $F_k$ hold:
 \begin{enumerate}
\item $iF_k(v,Jv)\geq ckg(v,v), \forall v\in D$.
\item $|{F_k}_{\mid D}-{F_k}_{\mid D}^{1,1}|_{g}\leq \delta k^{1/2}$.
\item $|\nabla^jF_k|_g\leq C_jk$.
\end{enumerate}
\end{definition}

Another motivation for the previous definition is the case of Levi-flat CR manifolds,
 where according to the results of Ohsawa and Sibony  \cite{OS00} leafwise
  positivity grants the existence of plenty of CR sections (with an appropriate twisting by a line bundle).

The fundamental example of an ample bundle is the pre-quantum line bundle $L$ of an integral 2-calibrated manifold $(M,D,\omega)$  (with  $c=2\pi$, $\delta=0$). Its tensor powers $L^{\otimes k}$ define a very
 ample sequence of line bundles.

From now on we will only consider almost CR structures on 2-calibrated manifolds
 defined by compatible almost complex structures and metrics. Similarly, we will
  only consider the very ample sequence $L^{\otimes k}$.

 For any $\tau_k\in\Gamma(L^{\otimes k})$ we use $J$ to split the restriction of $\nabla\tau_k$
  to $D$
  \[\nabla_D\tau_k=\partial\tau_k+\bar{\partial}\tau_k,\; \partial\tau_k\in \Gamma(D^{*1,0}\otimes L^{\otimes k}),\;
\bar\partial\tau_k\in \Gamma(D^{*0,1}\otimes L^{\otimes k}).\]

We can see $\bar{\partial}\tau_k$ as a section of $T^*M\otimes L^{\otimes k}$ by
 declaring it to vanish on $D^{\perp}$, and then use the Levi-Civita connection
  on $T^*M$ to define $\nabla^{r-1}\bar{\partial}\tau_k\in \Gamma(T^*M^{\otimes r}\otimes L^{\otimes k})$.

Let us denote the rescaled metric $kg$ by $g_k$.

\begin{definition}\label{def:sucah}  A sequence of sections   $\tau_k$ of $L^{\otimes k}$
 is approximately $J$-holomorphic (or approximately holomorphic or simply A.H.)
  if positive constants ${(C_j)}_{j\geq 0}$ exist such that

\[  |\nabla^j\tau_k|_{g_k}\leq  C_j,\;\; |\nabla^{j-1} \bar{\partial}\tau_k |_{g_k} \leq C_j k^{-1/2}. \]

If we want to make the bounds explicit speak of an A.H.($C_j$) sequence.
\end{definition}

\begin{remark} The original notion of A.H. sequence introduced in  \cite{IMP01,Pr02}
 is a bit more general  than definition \ref{def:sucah}.   The  difference -as well as the fact that
  only a finite number of derivatives were taken into account-
   is that the direction orthogonal to $D$ had a different treatment.
    The main theorem of \cite{IMP01} produced appropriate A.H. sequences
     of sections with good control on any finite number of  derivatives
     along $D$, but little  along $D^{\perp}$. Using the relative theory
      one can obtain solutions with control in all directions,
       so we can avoid using the technically more complicated definition of \cite{IMP01,Pr02}.
\end{remark}

\section{The local approximately holomorphic theory}\label{sec:locah}

 Maybe the most important idea on Donaldson's  work \cite{Do96} was the construction
  of localized A.H. sections (inspired in the work of Tian \cite{Ti89}) by adopting
   a unitary point of view instead of a holomorphic one.
The use of a unitary connection in a Darboux chart allowed him to
find a model for the coupled Cauchy-Riemann equation invariant under
rescaling -provided one
 worked in the appropriate tensor power of the pre-quantum line bundle-
  and explicitly write concentrated solutions giving rise to the so called reference sections.

The local approximately holomorphic theory, both using an intrinsic construction or the symplectization
 to be introduced in subsection \ref{subsec:symp}, is based on the choice
  of appropriate families of charts. In the intrinsic local theory we need
  as well a local model for the coupled Cauchy-Riemann equations and a good
  choice of explicit solution.

For 2-calibrated manifolds the local model for the intrinsic approximately holomorphic theory -that can only be achieved asymptotically when $k\rightarrow\infty$-  is the following:

\begin{itemize}
\item The domain is $\mathbb{C}^n\times\mathbb{R}$, with coordinates $z^1,\dots,z^n,s$
(sometimes we write them as $x^1,\dots,x^{2n+1}$ or $x^1,\dots,x^{2n},s$).
\item The distribution $D_h$ is the tangent space to the level hyperplanes of the vertical or real coordinate $s$.
\item The identification of each leaf with $\mathbb{C}^n$ means that we have fixed
 the leafwise standard almost complex structure $J_0$.
\item The metric is the  Euclidean one $g_0$ with Levi-Civita connection $\mathrm{d}$
 (usual partial derivatives), and  the distance is the Euclidean norm $|\cdot|$.

\item The 2-form in the fixed coordinates is required to be
\begin{equation}\label{eqn:sympform}
\omega_\mathrm{std}=\frac{i}{2}\sum_{i=1}^ndz^i\wedge d\bar{z}^i.
\end{equation}

\item We ask for a choice of unitary trivialization of the line bundle whose connection form is
\begin{equation}\label{eqn:connform}
A=\frac{1}{4}\sum_{i=1}^nz^id\bar{z}^i-\bar{z}^idz^i.
\end{equation}
\end{itemize}

In $\mathbb{R}^N$ with coordinates $x^1,\dots,x^N$ let $\mathbb{R}^p$ denote
the distribution by $p$-planes $\textrm{span}<\partial/\partial x^{i_1},\dots,\partial/\partial x^{i_p}>$,
 $1\leq i_1<\cdots<i_p\leq N$; its Euclidean orthogonal is denoted by $\mathbb{R}^{N-p}$.
   If we have a distribution $D'$ of dimension $p$ in $\mathbb{R}^N$ which is
   transverse to $\mathbb{R}^{N-p}$, we can measure its distance to $\mathbb{R}^p$
   to order $j$ with respect to the flat connection d as follows: $D'$ can be identified with an element of $\mathrm{Hom}(\mathbb{R}^p,\mathbb{R}^{N-p})$. We let $v^{i_l}$, $l=1,\dots,p$,
    be the vector field in $\mathbb{R}^{N-p}$ such that $\partial/\partial x^{i_l}+v^{i_l}\in D'$. Then we define

\[|\textrm{d}^j(\mathbb{R}^p-D')|_{g_0}=\textrm{max}\{|\textrm{d}^jv^{i_1}|_{g_0},\dots,|\textrm{d}^jv^{i_p}|_{g_0}\},\]
which by definition is coordinate dependent.

In the previous local model let us denote the line field spanned by  $\partial/\partial s$
 by $D_v$. According to the  previous paragraph we can measure the distance in
 $\mathbb{C}^n\times \mathbb{R}$ to $D_h$ (resp. $D_v$) of any codimension one
 (resp. dimension one) distribution transverse to $D_v$ (resp. $D_h$).

\begin{definition}\label{def:dcharts} Let $\varphi_{k,x}\colon (\mathbb{C}^n\times\mathbb{R},0)\rightarrow (U_{k,x},x)$,
 for all  $x\in M$ and all  $k \gg 1$, be a family of charts with coordinates $z_k^1,\dots,z_k^n,s_k$.
  We call them a family of approximately holomorphic coordinates  if  there exist  constants
  independent of $k,x$ (uniform) so that  the following estimates hold for
  all $k\gg 1$ at the points of $B(0, \rho k^{1/2})$, $\rho>0$:
\begin{enumerate}\item The Euclidean and the induced metric  are comparable to any order, i.e.
\[\frac{1}{\gamma} g_0\leq g_k\leq \gamma g_0, \,\gamma>0, \mathrm{and} \;|\nabla^j\varphi_{k,x}^{-1}|_{g_0}\leq O(k^{-1/2}), \forall j\geq 2,\]
where $\nabla$ denotes the Levi-Civita connection with respect to $g$.
\item The kernel of $\omega$, which is $D^{\perp}$, is sent to a line field
 $\varphi_{k,x}^*D^{\perp}$ transverse to $D_h$ and such that
\begin{align*}
|\varphi_{k,x}^*D^{\perp}-D_v|_{g_0}& \leq  |(z_k,s_k)|O(k^{-1/2}),\\
|\textrm{d}^j(\varphi_{k,x}^*D^{\perp}-D_v)|_{g_0} &\leq  O(k^{-1/2}), \;\forall j\geq 1.
\end{align*}
The pullback of $D$ is transverse to $D_v$ and
\begin{align*}
|\varphi_{k,x}^*D-D_h|_{g_0} &\leq |(z_k,s_k)|O(k^{-1/2}),\\
|\textrm{d}^j(\varphi_{k,x}^*D-D_h)|_{g_0} &\leq   O(k^{-1/2}), \;\forall j\geq 1.
\end{align*}
 \item Regarding the antiholomorphic components,
\begin{align*}
|\bar{\partial}\varphi^{-1}_{k,x}(z_k,s_k)|_{g_0}&\leq |(z_k,s_k)|O(k^{-1/2}),\\
|\nabla^j \bar{\partial}\varphi^{-1}_{k,x}(z_k,s_k)|_{g_0}&\leq O(k^{-1/2}), \forall j\geq 1,
\end{align*}
where  $\bar{\partial} \varphi^{-1}_{k,x} $ is the antiholomorphic component of
 $\nabla_D(\pi_{D_h}\circ \varphi^{-1}_{k,x})$, with $\pi_{D_h}\colon\mathbb{C}^n\times \mathbb{R}\rightarrow \mathbb{C}^n$
  the projection onto the first factor.
\end{enumerate}

   We speak of Darboux coordinates when the additional condition  $\varphi_{k,x}^*k\omega=\omega_0$ holds.
\end{definition}

\begin{remark}\label{rem:splitting}
According to condition (2) (resp. (3)) we have $\varphi_{k,x}^*D=D_h$, $\varphi^*_{k,x}D^{\perp}=D_v$ (resp. $\varphi_{k,x}^*J=J_0$)
 at the origin. For most of our constructions it is enough to require the equality up to a summand
  of size $O(k^{-1/2})$ at most, but since these equalities are needed to prove results
   concerning pseudo-holomorphic jets (in particular the identities concerning local
    representations and subsets of transverse holonomy of lemma \ref{lem:locrep}) we
    choose to ask for them from the very beginning.
\end{remark}

\begin{remark}\label{rem:ahcoord}
If we are in an almost complex manifold then conditions (1) and (3) ((2) makes no sense)
 recover the notion of approximately holomorphic charts (resp. Darboux charts if we
  add the Darboux condition on the 2-form).
\end{remark}

 A chart centred at a point for which the Darboux
condition holds can always be obtained: $(M,D,\omega)$ is a coisotropic submanifold
of its symplectization, as defined in lemma \ref{lem:sympl}. The local normal form
theorem for coisotropic submanifolds (\cite{Va87}, theorem  3.4.10) provides such a chart.  Families of Darboux charts can be constructed using the same local normal form. Since this would fall into the relative theory we prefer to give a different proof.

\begin{lemma}\label{lem:dchart} Let $(M,D,\omega)$ be a (compact) 2-calibrated manifold
(with $J$, $g$ already fixed). Then a family of Darboux charts can always be constructed.
\end{lemma}
\begin{proof}
Let us  fix a family of charts $\psi_x\colon \mathbb{R}^{2n+1}\rightarrow U_x$ depending
 smoothly on $x$, where $x\in M_1$  a small enough subset of $M$, so that $\psi_x^*D=D_h$, $\psi_x^*D^{\perp}=D_v$
  at the origin. Denote by $x^1,\dots,x^{2n},s$ the coordinates on $\mathbb{R}^{2n+1}$. We
  compose $\psi_x$ with the diffeomorphism $\Theta_x\colon \mathbb{R}^{2n+1}\rightarrow \mathbb{R}^{2n+1}$
    which is the identity on $\mathbb{R}^{2n}\times \{0\}$, preserves setwise the horizontal
    foliation $D_h$ and sends $\mathrm{Ker} \psi_x^*\omega$ to $D_v$. The diffeomorphisms  $\Theta_x$ depend smoothly on $x$.

Now  we fix $J_0$ to identify $\mathbb{R}^{2n+1}$ with $\mathbb{C}^n\times \mathbb{R}$
 and compose with an element of $\textrm{Gl}(2n,\mathbb{R})\subset \textrm{Gl}(2n+1,\mathbb{R})$
  (again depending smoothly on $x\in M_1$), so that we obtain charts $\varphi_x$  for
   which the pullback of $J$ at the origin equals $J_0$.

By compactness  $M$ can be covered with a finite number of subsets $M_1,\dots,M_h$ in
 which the above charts  can be constructed. In this way we obtain charts centred at every
  $x\in M$ (we might have more than one chart for each $x\in M$, but that is not relevant)
   so that the bounds on tensors pulled back from $M$ to a ball of fixed radius in the domain of the charts will not depend  on $x$.

 We define $\varphi_{k,x}$ to be the composition $\varphi_x\circ\gamma_{k^{-1/2}}$, where
   $\gamma_{k^{-1/2}}\colon \mathbb{C}^n\times \mathbb{R}\rightarrow \mathbb{C}^n\times \mathbb{R}$
    is the homothety by factor $k^{-1/2}$.
The equalities at the origin together with the smooth dependence on $x$ of the constructions
 previous to the rescaling, imply that we have obtained approximately holomorphic coordinates.

To obtain Darboux charts we need to modify $\varphi_{k,x}$ as
follows: we apply  Darboux' lemma with estimates  (lemma 2.2 in
\cite{Au01}) to the almost complex manifolds
$(\mathbb{C}^n\times\{0\},\varphi_{k,x}^*J_{\mid
\mathbb{C}^n\times\{0\}},\varphi_{k,x}^*g_{\mid
\mathbb{C}^n\times\{0\}})$ and the 2-forms
$\varphi_{k,x}^*\omega_{\mid \mathbb{C}^n\times\{0\}}$. We get
diffeomorphisms $\Psi_{k,x} $ on this leaf that are extended to
$\mathbb{C}^n\times\mathbb{R}$ independently of the vertical
coordinate $s_k$.  The bounds on $\Psi_{k,x}$ and their derivatives
coming from  lemma 2.2 in \cite{Au01} imply that the compositions $
\varphi_{k,x}\circ\Psi_{k,x}\colon
(\mathbb{C}^n\times\mathbb{R},0)\rightarrow (U_{k,x},x)$ still
define approximately holomorphic coordinates.  Moreover, we can
assume $(\varphi_{k,x}\circ\Psi_{k,x})^*J=J_0$ at the origin.

 Since $\partial/\partial s_k$ generates the kernel of
  \[{(\varphi_{k,x}\circ\Psi_{k,x})}^*\omega=\sum_{1\leq i <l \leq 2n}\omega_{il}dx^i_k\wedge dx^l_k+\sum_{1\leq i \leq 2n}\omega_idx^i_k\wedge ds_k,\]
all $\omega_i$ vanish. Closedness implies that each function $\omega_{il}$
 is independent of $s_k$. Therefore ${(\varphi_{k,x}\circ\Psi_{k,x})}^*\omega$
  is determined by its restriction  to $\mathbb{C}^n\times\{0\}$, which by construction
   is  ${\omega_\mathrm{std}}_{\mid \mathbb{C}^n\times\{0\}}$.
Thus, $\omega$ is sent to $\omega_\mathrm{std}$.
\end{proof}

Darboux charts are useful because there local computations become simpler.

Let $d_k$ denote the distance defined by the metric $g_k$.

Recall that in the domain of a Darboux chart we can always fix $\xi_{k,x}$ a
unitary trivialization of $L^{\otimes k}$ whose connection form is $A$ (equation (\ref{eqn:connform})).

\begin{lemma}\label{lem:localcheck1} Let $\varphi_{k,x}\colon (\mathbb{C}^{n}\times\mathbb{R},0)\rightarrow (U_{k,x},x)$
 be a family of Darboux charts with coordinates $x_k^1,\dots,x^{2n}_k,x^{2n+1}_k$. Let
  $F$ be a bundle associated to either  $TM$ or $D$ and let
   $F_{k,x}\rightarrow B(0,\rho k^{1/2})\subset \mathbb{C}^{n}\times\mathbb{R}$
    denote the pullback of  $F$ by $\varphi_{k,x}$. Associated
    to the Darboux coordinates there is a canonical trivialization  $\zeta_{k,x,j}$ of $F_{k,x}$.
     Let $T_k$ be  a sequence of sections of $F\otimes L^{\otimes k}$ and use the frames
      $\zeta_{k,x,j}\otimes \xi_{k,x}$  to write $\varphi_{k,x}^*T_k$ locally as a function
       $T'_{k,x}$. Let $P_j$ be a polynomial such that for any multi-index $\alpha$ of
       length $j=0,\dots,r$,  at the points of $B(0,\rho k^{1/2})$ and for all $k \gg 1$ we have:

\[{\left|\frac{\partial}{\partial x^\alpha_k}T'_{k,x}\right|}_{g_0}\leq P_j(|(z_k,s_k)|)O(k^{-1/2}).\]
Then $|\nabla^rT_k(y)|_{g_k}\leq Q_r(d_k(x,y))O(k^{-1/2})$, where the polynomial $Q_r$ depends only on $P_1,\dots,P_r$.
Conversely, from bounds using the global  metric elements $g_k,d_k,\nabla$ we obtain similar bounds
for the local Euclidean elements.
\end{lemma}
\begin{proof} This is a simple calculation based on items (1) and (2), and in the Darboux condition in
 definition \ref{def:dcharts}. Also notice that the presence of the connection form  and
 its derivatives is absorbed by the polynomial, since $|A|\leq O(|(z_k,s_k)|)$ and  its derivatives are of order $O(1)$.
\end{proof}

\begin{remark}\label{rem:modiflem} Lemma \ref{lem:localcheck1} admits different modifications.
  It holds in a similar fashion for bounds of order $O(1)$ instead of order $O(k^{-1/2})$;
   also for sections $T_k$ of $F$ instead of $F\otimes L^{\otimes k}$ (with $F_{k,x}$ locally
    trivialized by $\zeta_{k,x,j}$); it is also possible to consider the inequalities in
     the ball of (uniform) radius $\rho>0$, rather than $\rho k^{1/2}$. There is also a version for symplectic manifolds.
\end{remark}

Let $\bar{\partial}_0$  denote the (0,1)-component with respect to
$J_0\colon \mathbb{C}^n\times \mathbb{R}\rightarrow  \mathbb{C}^n\times \mathbb{R}$ of the
leafwise derivation operator $\textrm{d}_{D_h}$.

\begin{lemma}\label{lem:localcheck2} Let $\varphi_{k,x}\colon (\mathbb{C}^{n}\times\mathbb{R},0)\rightarrow (U_{k,x},x)$
 be a family of Darboux charts with coordinates $x_k^1,\dots,x^{2n}_k,s_k$. Let
   $L_{k,x}\rightarrow B(0,\rho k^{1/2})\subset \mathbb{C}^{n}\times\mathbb{R}$
   denote the pullback of  $L^{\otimes k}$ by $\varphi_{k,x}$.   Let $\tau_k$ be
    a sequence of sections of $L^{\otimes k}$ such that
     $\varphi_{k,x}^*\tau_k=f_{k,x}\xi_{k,x}$. Let $P_j, P_{j'}$ be  polynomials
     such that for any multiindices $\alpha$, $\beta$  of length
      $j=0,\ldots,r-1$, and $j'=0,\dots,r$, respectively,  at the points
       of $B(0,\rho k^{1/2})$ and for all $k \gg 1$ the following inequalities hold:

\begin{align}\label{eqn:localcheck1bis}
{\left|\frac{\partial}{\partial x^\beta_k}f_{k,x}\right|}_{g_0} &\leq P_j'(|(z_k,s_k)|)O(1).\\
\label{eqn:localcheck2}\left|\frac{\partial}{\partial x^\alpha_k}(\bar{\partial}_0+A^{0,1})f_{k,x}\right|_{g_0} &\leq P_j(|(z_k,s_k)|)O(k^{-1/2}).
\end{align}
Then we have
\begin{align}\label{eqn:localcheck1bisa}
|\nabla^{r}\tau_k(y)|_{g_k}&\leq Q'_r(d_k(x,y))O(1),\\
\label{eqn:localcheck2a}
|\nabla^{r-1}\bar{\partial}\tau_k(y)|_{g_k}&\leq Q_{r-1}(d_k(x,y))O(k^{-1/2}),
\end{align}
where the polynomial $Q_{r-1}$ (resp. $Q'_r$) depends only on
 $P_1,\dots,P_{r-1},P_1',\dots,P_r'$ (resp. $P_1',\dots,P_r'$).
 Conversely, from bounds using  $g_k,d_k,\nabla,J$ we obtain similar bounds for $g_0, |\cdot|, \mathrm{d}+A, J_0$.
\end{lemma}
\begin{proof}
The equivalence between equations (\ref{eqn:localcheck1bis}) and (\ref{eqn:localcheck1bisa}) is the content
 of lemma \ref{lem:localcheck1}, but for bounds of order O(1) (see remark \ref{rem:modiflem}). The equivalence of equations
(\ref{eqn:localcheck1bis}), (\ref{eqn:localcheck2}) and equations (\ref{eqn:localcheck1bisa}), (\ref{eqn:localcheck2a})
 follows again easily from the properties of Darboux charts. We sketch the case $r=1$.

 From now on $\varphi_{k,x}^*J, \varphi_{k,x}^*D,\varphi_{k,x}^*g_k$, and all the tensors and sections
  pulled back to the domain of the charts will be  denoted by $J, D, g_k,\dots$ whenever there is no risk of confusion.

Let $e_i$ be any of the local vector fields associated to the first $2n$ coordinates. By condition
 (2) in definition \ref{def:dcharts} there exists $u_i$ a local vector field such that  $e_i+u_i$ is tangent to $D$ and
\begin{equation}\label{eq:cond1}
 |u_i|_{g_0}\leq |(z_k,s_k)|O(k^{-1/2}),\; |\textrm{d}^ju_i|_{g_0}\leq  O(k^{-1/2}),\; j\geq 1.
 \end{equation}
The endomorphism $J$ is defined on $D$. We can use the orthogonal  projection w.r.t $g_0$ onto
 $D_h$ to induce out of $J$ another almost complex structure $J_{D_h}\colon D_h\rightarrow D_h$.

Condition (3) in definition \ref{def:dcharts} implies that

\begin{equation}\label{eq:cond2}
|J_0-J_{D_h}|_{g_0}\leq |(z_k,s_k)|O(k^{-1/2}),\;
 |\textrm{d}^j(J_0-J_{D_h})|_{g_0}\leq  O(k^{-1/2}),\; j\geq 1.
 \end{equation}

By definition $\bar{\partial}_{e_i+u_i}\tau_k=1/2\nabla_{e_i+u_i}\tau_k+i/2\nabla_{J(e_i+u_i)}\tau_k$.

Equation (\ref{eq:cond1}) combined with lemma \ref{lem:localcheck1} implies
\[|\nabla_{u_i}\tau_k|_{g_k}\leq P'_1(d_k(x,y))O(k^{-1/2}).\]
Again equations (\ref{eq:cond1}) and (\ref{eqn:localcheck1bisa}), condition (3) in definition
\ref{def:dcharts}, and lemma \ref{lem:localcheck1} imply
\[|\nabla_{J(e_i+u_i)}\tau_k-\nabla_{J_he_i}\tau_k|_{g_k}\leq P''_1(d_k(x,y))O(k^{-1/2}).\]
Therefore  the  bounds in equation (\ref{eqn:localcheck2a}) we want for $\bar{\partial}_{e_i+u_i}\tau_k$
 are equivalent to the  same kind of bounds for
\[1/2\nabla_{e_i}\tau_k+i/2\nabla_{J_he_i}\tau_k,\]
 and  by equation (\ref{eq:cond2}) for
 \[1/2\nabla_{e_i}\tau_k+i/2\nabla_{J_0e_i}\tau_k,\]
and by definition

\[1/2\nabla_{e_i}\tau_k+i/2\nabla_{J_0e_i}\tau_k=((\bar{\partial}_0+A^{0,1})_{e_i}f_{k,x})\xi_{k,x}.\]

Bounds for higher order derivatives are proven similarly.
 \end{proof}

 \begin{definition}(see \cite{Au01}, definition 2.2) A sequence of sections of $L^{\otimes k}$ has
  Gaussian decay with respect to $x$ if  there exist polynomials $(P_j)_{j\geq 0}$  and a constant $\lambda>0$,
    so that $\;\forall y \in M$  and $\forall j\geq 0$
\[  |\nabla^j \tau_k(y)|_{g_k} \leq  P_j(d_k(x,y))\mathrm{e}^{-\lambda d_k(x,y)^2}.\]

  \end{definition}

The main purpose of the use of Darboux charts is the construction of reference sections $\tau_{k,x}^{\mathrm{ref}}$.

\begin{corollary}\label{cor:refsections} Let $(M,D,\omega)$ be a compact 2-calibrated manifold. Then for
all $x\in M$ A.H. sections $\tau_{k,x}^{\mathrm{ref}}$ with Gaussian decay with respect to $x$ can be constructed.
 The bounds are uniform on $k,x$ and these sections have norm greater than some  constant $\kappa$ in
  $B_{g_k}(x,\rho)$, where $\kappa,\rho>0$ are uniform on $k,x$.
\end{corollary}

\begin{proof} We follow Donaldson's ideas in \cite{Do96}, section 2. Let us  fix Darboux charts and
 $\xi_{k,x}$  trivializations of $L^{\otimes k}$ for which the connection form is $A$. Let $\beta$ be
  a standard  cut-off function of a single variable, with $\beta(t)=1$ when $|t|\leq 1/2$ and $\beta(t)=0$ when $|t|\geq 1$.

Define $\beta_k(z_k,s_k)=\beta(k^{-1/6}|(z_k,s_k)|)$.

In the points where the derivatives of $\beta_k$ do not vanish we have $|(z_k,s_k)|\geq Ck^{1/6}$,
$C$ uniform (on $k, x$). Using this inequality we deduce
\begin{align}\label{eq:cutoff}\nonumber
|\textrm{d}\beta_k|_{g_0}&\leq |(z_k,s_k)|^2O(k^{-1/2}),\\\nonumber
 |\textrm{d}^2\beta_k|_{g_0}&\leq |(z_k,s_k)|O(k^{-1/2}), \\ |\textrm{d}^j\beta_k|_{g_0}&\leq O(k^{-1/2}),\; j\geq 3.
\end{align}

Consider the function $f(z_k,s_k)=\mathrm{e}^{-|(z_k,s_k)|^2/4}$. We have

\begin{equation}\label{eq:cr}
\bar{\partial}_0f+A^{0,1}f=0.
\end{equation}

The reference sections are

\begin{equation}\label{eq:refsec}
\tau_{k,x}^{\mathrm{ref}}:=\beta_k f \xi_{k,x}.
\end{equation}

Equation (\ref{eq:cutoff}) implies that for any multi-index $\alpha$ of length $j\leq r$,

\[\left|\frac{\partial}{\partial x^\alpha}\beta_k f\right|_{g_0}\leq P_j(|(z_k,s_k)|)|f|O(1).\]
 Therefore, lemma \ref{lem:localcheck1} for bounds of type $\mathrm{e}^{-\lambda' |(x,y)|^2}O(1)$, $\lambda'>0$,
  gives the Gaussian decay with respect to $x$:

\[|\nabla^r\tau_{k,x}^{\mathrm{ref}}(y)|_{g_k}\leq Q_r(d_k(x,y))\mathrm{e}^{-\lambda d_k(x,y)^2}O(1), \; \lambda>0,\]
where $\lambda$  appears when relating the distance induced by $g$ and $g_0$.
The Gaussian decay also implies  \[|\nabla^r\tau_{k,x}^{\mathrm{ref}}|_{g_k}\leq  O(1).\]

The bound for $|\nabla^{r-1}\bar{\partial}\tau_{k,x}^{\mathrm{ref}}|_{g_k}$ is obtained using the same ideas:
from equations (\ref{eq:cutoff}) and (\ref{eq:cr}) it follows that  for any multi-index $\alpha$ of length $j\leq r-1$

\[\left|\frac{\partial}{\partial x^\alpha}(\bar{\partial}_0+A^{0,1})\beta_k f \right|_{g_0}\leq P_j(|(z_k,s_k)|)|f|O(k^{-1/2}).\]

Lemma \ref{lem:localcheck2}  for bounds of type
$\mathrm{e}^{-\lambda' |(x,y)|^2}O(1)$, $\mathrm{e}^{-\lambda'
|(x,y)|^2}O(k^{-1/2})$, $\lambda'>0$ (in equations
(\ref{eqn:localcheck1bis}) and (\ref{eqn:localcheck2}) resp.), gives for
some $\lambda>0$

\[|\nabla^{r-1}\bar{\partial}\tau_{k,x}^{\mathrm{ref}}|_{g_k}\leq Q_r(d_k(x,y))\mathrm{e}^{-\lambda d_k(x,y)^2}O(k^{-1/2})\leq O(k^{-1/2}).\]

The existence of constants $\kappa,\rho>0$ such that $|\tau_{k,x}^{\mathrm{ref}}|\geq \kappa $ in $B_{g_k}(x,\rho)$,
 can be easily checked.
\end{proof}

We observe that many of the inequalities we are using  (for global tensors) have the same pattern. We will introduce
 a definition that will avoid the excessive appearance in the notation of  such inequalities.

Let $E$ be a Hermitian bundle with connection, $F$ a bundle associated either to $TM$ or to $D$,
 and let $E_k$ denote  the sequence $F\otimes E\otimes L^{\otimes k}$.

\begin{definition}\label{def:apeq} Let $T_{k,x}$, $x\in M$, be a family of sequences of sections of $E_k$.
 We say that $T_{k,x}$ is $C^r$-approximately vanishing (or that the sequence vanishes in the $C^r$-approximate sense)
  and denote it by $T_{k,x}\approxeq_r 0$, if  positive constants $C_0,\dots,C_r$ exist so that
\begin{equation}\label{eq:inequality0}
|\nabla^jT_{k,x}|_{g_k}\leq C_jk^{-1/2}, \; j=0,\dots,r.
\end{equation}

There is an analogous  definition for sequences $T_k$ of sections of  $E_k$  (i.e. without extra dependence on the point $x\in M$).
\end{definition}

Using the above language one of the conditions for a sequence $\tau_k$ of $L^{\otimes k}$ to be A.H.
(definition \ref{def:sucah}) is that $\bar{\partial} \tau_k\in \Gamma(D^{*0,1}\otimes L^{\otimes k})$
 has to be approximately vanishing.

\begin{remark}\label{rem:proj} Given $\tau_k$ an approximately holomorphic sequence of sections
 of $L^{\otimes k}$, we have defined $\nabla^{r-1}\bar{\partial}\tau_k\in T^*M^{\otimes r}\otimes L^{\otimes k}$ by
 taking covariant derivatives of $\bar{\partial}\tau_k$ thought of as a section of $T^*M\otimes L^{\otimes k}$.
  We might have equally defined $\nabla^{r-1}\bar{\partial}\tau_k$ as the image of $\nabla^r\tau_k$ by the projection
$\bar{p}_r\colon T^*M^{\otimes r}\otimes L^{\otimes k}\rightarrow T^*M^{\otimes r-1}\otimes D^{*0,1}\otimes L^{\otimes k}$,
 for using Darboux charts and lemmas \ref{lem:localcheck1} and \ref{lem:localcheck2} (with
 the inequalities $|\nabla^j\tau_k|_{g_k}\leq O(1), \;j\geq 0$), one checks that $\bar{\partial}\tau_k\approxeq 0$ if and only if
$|\bar{p}_j(\nabla^j\tau_k)|_{g_k}\leq O(k^{-1/2}), \;j\geq 1$.
\end{remark}

\subsection{Relative approximately holomorphic theory and symplectizations}\label{subsec:symp}

\begin{definition} Let $(P,\Omega)$ be a symplectic manifold and $(M,D,\omega)$ a 2-calibrated manifold.
 We say that $l\colon M\hookrightarrow P$ embeds $M$ as a 2-calibrated submanifold of $P$ if $l^*\Omega=\omega$.
\end{definition}

\begin{lemma}\label{lem:sympl} Let  $(M,D,\omega)$ be a compact co-oriented 2-calibrated manifold. Then it
 is possible to define a symplectization so that $(M,D,\omega)$ embeds as a 2-calibrated  submanifold.
  Any fixed compatible almost complex structure and metric can be extended to a compatible
   almost complex structure and metric in the symplectization.
\end{lemma}
\begin{proof} Let $J$ and $g$ be fixed compatible almost complex structure and metric. The symplectization  $(M\times[-\epsilon,\epsilon],J,g,\Omega)$
 is constructed as follows: let $t$ be the coordinate of the interval. Let $\alpha$ be the unique 1-form
  of pointwise norm 1 (and positively oriented) whose kernel is $D$. The closed $2$-form $\Omega$ is defined
  to be $\omega+d(t\alpha)$, where $\alpha$ and $\omega$ represent the pullback of the corresponding forms
   to $M\times [-\epsilon,\epsilon]$. If $\epsilon $ is chosen small enough then $\Omega$ is symplectic.

In the points of  $M$ the almost complex structure is extended by sending the positively oriented $g$-unitary
 vector in $D^\perp $ to $\partial/\partial t$; in those points $\partial/\partial t$ is also defined to
  have norm 1 and to be orthogonal to $TM$. It is routine to further extend $J$ to a compatible almost complex structure on the
  symplectization. The metric defined by $\Omega$ and the almost complex structure also extends $g$.
  We will not use different notation for the extension of the almost complex structure and metric if
  there is no risk of confusion.

We also fix $G$ a $J$-complex distribution  on the symplectization restricting to $D$ at the points of $M$.
 To do that we choose any line field  that at the points of $M$ contains $\partial/\partial t$; this
 line field spans a complex line field. Its orthogonal with respect to $g$ is by construction $J$-complex and extends $D$.
\end{proof}

\begin{remark}\label{rem:2calembedding}  We want to work out a relative theory for embeddings in arbitrary symplectic manifolds -not just in symplectizations- because of our applications to CR manifolds, where we need an ambient complex manifold with plenty of holomorphic sections.
 \end{remark}

Let $(M,D,\omega)$ be a 2-calibrated submanifold of $(P,\Omega)$. Let us fix $J$ a compatible almost complex structure on $(P,\Omega)$
 so that $D$ is $J$-invariant, and let us define $g=\Omega(\cdot,J\cdot)$. The restriction of $(J,g)$ to $(M,D)$
  induces an almost CR structure. We also choose $G$ a $J$ complex distribution
   that coincides with $D$ at the points of $M$.  The main example to have in mind is the symplectization
    of $(M,D,\omega)$ with an almost complex structure as defined in lemma \ref{lem:sympl}.

We have at our disposal the approximately holomorphic theory for symplectic manifolds \cite{Au01}. At
 this point we pause to warn the reader that throughout this subsection and the rest of the paper we
  will be using  A.H. sequences of sections  defined in both 2-calibrated (definition \ref{def:sucah})
   and symplectic manifolds (see  definitions in \cite{Au01} or definition \ref{def:sucah} for an
   almost complex base space). Whenever there is no risk of confusion about the base space we
    will just speak about A.H. sequences of sections.

Let  $(L_\Omega,\nabla)\rightarrow (P,\Omega)$ be the pre-quantum line bundle. Its powers $(L^{\otimes k}_{\Omega},\nabla_k)$ define a very ample sequence of line bundles (in the sense
  of  \cite{Au01}), which restricts   to a very ample sequence of line bundles
  $(L^{\otimes k},\nabla_k)\rightarrow (M,D,J,g_k)$ (definition \ref{def:ampleb}).

One expects that if  $\tau_k\in \Gamma(L^{\otimes k}_{\Omega})$ is a (symplectic) A.H.
 sequence of sections, then ${\tau_k}_{\mid M}\colon M\rightarrow L^{\otimes k}$ is also an A.H.
  sequence of sections (definition \ref{def:sucah}). Even more, we will see that it is possible
   to construct reference sections by restricting (symplectic) reference sections centred
    at points of $M$. The key point to prove these results is the choice of appropriate charts.

Recall that in  $\mathbb{C}^p=\mathbb{R}^{2p}$ we denote the foliation whose leaves are associated
to $g$ distinguished complex coordinates (resp. $d$ distinguished real coordinates) by $\mathbb{C}^g$
 (resp. $\mathbb{R}^{d}$); its Euclidean orthogonal is denoted by $\mathbb{C}^{p-g}$ (resp. $\mathbb{R}^{2p-d}$).
    From now on if we compare the distance of $\mathbb{C}^g$ to any distribution of the same dimension,
     we will assume the latter to be transverse to $\mathbb{C}^{p-g}$.

\begin{definition}\label{def:adaptedAHcoord} Let $(P,\Omega)$ be a compact symplectic manifold  and $G$ a
$J$-complex distribution of complex dimension $g$. A family of (symplectic) approximately
holomorphic coordinates (resp. Darboux charts)
$\varphi_{k,x}\colon(\mathbb{C}^p,0)\rightarrow (U_{k,x},x)$ is said
to be adapted to $G$ if
\begin{eqnarray*}
|\mathbb{C}^g-G|_{g_0}\leq |(z_k,s_k)|O(k^{-1/2}),&&|\mathrm{d}^j(\mathbb{C}^g-G)|_{g_0} \leq O(k^{-1/2}), \;\forall j\geq 1.\\
|\mathbb{C}^{p-g} -G^{\perp}|_{g_0}\leq |(z_k,s_k)|O(k^{-1/2}),&&|\mathrm{d}^j(\mathbb{C}^{p-g} -G^{\perp})|_{g_0} \leq O(k^{-1/2}), \;\forall j\geq 1.
\end{eqnarray*}
\end{definition}

The existence of approximately holomorphic (resp. Darboux) charts
adapted to $G$ is straightforward:
 once we have approximately holomorphic (resp. Darboux) charts, we compose with a unitary
  transformation sending $G$ to $\mathbb{C}^g$ at the origin.

Given a 2-calibrated submanifold $(M,D)\hookrightarrow (P,\Omega)$, in order to select coordinate
 charts adapted to $M$ we fix a distribution $T^{||}M$ defined in a tubular neighborhood
 of $M$ as  follows: the neighborhood is defined by flowing a little bit the geodesics normal to
  $M$. For each point $y$ in the neighborhood, let $x\in M$ be the starting point of the unique
   geodesic normal to $M$ through $y$. Then  $T^{||}_yM$ is the result of parallel transport of $T_xM$ along that geodesic.

\begin{definition}\label{def:ahm} Let $(M,D)\hookrightarrow (P,\Omega)$ be a 2-calibrated submanifold, $G$ a $J$-complex distribution which extends $D$ (perhaps defined in a tubular neighborhood of $M$), and  $T^{||}M$ a distribution constructed as above. A family of (symplectic) A.H. coordinates
   $\varphi_{k,x}\colon (\mathbb{C}^p,0)\rightarrow (U_{k,x},x)$ (centred at every point of $P$)
    is adapted to $(M,G)$, if it is  adapted to $G$ and for the charts centred at points of $M$ the following conditions hold:
\begin{enumerate}
\item $M$ sits in each chart as a fixed linear subspace $\mathbb{R}^{2n+1}\times\{0\}\subset \mathbb{C}^p$ and
 at the origin $D=\mathbb{R}^{2n}\times \{0\}\subset \mathbb{R}^{2n+1}\times\{0\}$,
 $D^\perp=\{0\}\times\mathbb{R}\subset \mathbb{R}^{2n+1}\times\{0\}$.
\item $|\mathbb{R}^{2n+1}-T^{||}M|_{g_0}\leq |(z_k,s_k)|O(k^{-1/2})$,
$|\mathrm{d}^j(\mathbb{R}^{2n+1}-T^{||}M)|_{g_0} \leq O(k^{-1/2})$, $\forall j\geq 1$.
\end{enumerate}

We speak of A.H. charts adapted to $(M,G)$ and Darboux over $M$ if
\begin{equation}\label{eq:darbcond}
\varphi_{k,x}^*\omega_{\mid M}=\omega_0.
\end{equation}

\end{definition}

\begin{lemma}\label{lem:reldar} Let $(M,D)\hookrightarrow (P,\Omega)$ be a 2-calibrated submanifold.
 Then   approximately holomorphic  charts adapted to $(M,G)$ and Darboux over $M$ can always be constructed.
\end{lemma}

\begin{proof}
We start by fixing approximately holomorphic coordinates adapted to $G$. Then we forget about
 the ones centred at points of $M$, that are going to be substituted by new ones.
For every $x\in M$  we fix initial charts $\varphi_{x}$  depending smoothly on the center
 -at least in a small neighborhood about each point- with $(\varphi_{x}J^*,\varphi_{x}^*g)=(J_0,g_0)$
  at the origin. Then  we compose with maps $\Theta_x\colon (\mathbb{C}^p,0)\rightarrow (\mathbb{C}^p,0)$
   that are tangent to the identity map at the origin and send $M$ to a vector space  in $\mathbb{C}^{p}$.
    The image of the distribution $D$ is $J_0$-complex at the origin.  By composing with a unitary
    transformation $(D_x,T_xM)$ can be assumed to be sent to
    $(\mathbb{C}^{n}\times\{0\},\mathbb{R}^{2n+1}\times \{0\})\subset \mathbb{R}^{2p}$.

Next we essentially apply lemma \ref{lem:dchart} on the leaf $\mathbb{R}^{2n+1}\times\{0\}\subset \mathbb{R}^{2p}$
to get Darboux charts for $M$:  let $\Xi_x\colon \mathbb{R}^{2n+1}\rightarrow \mathbb{R}^{2n+1}$ be
 the map which is the identity on $\mathbb{C}^n\times \{0\}$, preserves the foliation by complex
  hyperplanes, and sends the kernel of $\omega$ to the ``vertical'' or real line field in
  $\mathbb{R}^{2n+1}\times\{0\}$.  We extend it to a diffeomorphism of $\mathbb{R}^{2p}$
   independently of the coordinates $x^{2n+2},\dots,x^{2p}$. Since the map is by construction
    tangent to the identity at the origin, we keep the properties at the origin described in the previous paragraph.

We now  apply Darboux' lemma on $\mathbb{R}^{2n}\times\{0\}$ for each $x$. The result  is a diffeomorphism
 on $\mathbb{R}^{2n}$ that can be assumed to preserve $J_0$ at the origin. We extend it independently
 of $x^{2n+1},\dots,x^{2p}$ to a diffeomorphism of $\mathbb{C}^p$.
Notice that  $(D_x,T_xM)$ goes to $(\mathbb{R}^{2n}\times\{0\},\mathbb{R}^{2n+1}\times\{0\})$, $J_x$
 to ${J_0}$, $G_x\oplus G^\perp_x$ to $\mathbb{C}^{n}\oplus \mathbb{C}^{p-n}$, and $\mathrm{Ker}{\omega}_{\mid D_x}$
  to the Euclidean orthogonal of $\mathbb{R}^{2n}\times\{0\}\subset\mathbb{R}^{2n+1}\times\{0\}$.
   Hence if we   apply the homothety $\gamma_{k^{-1/2}}\colon \mathbb{R}^{2p}\rightarrow \mathbb{R}^{2p}$
    we obtain a family of charts with the desired properties.
\end{proof}

\begin{lemma}\label{lem:ahm} A family of A.H. charts  $\varphi_{k,x}\colon(\mathbb{C}^p,0)\rightarrow (U_{k,x},x)$
   adapted to $(M,G)$ and Darboux over $M$ constructed as in lemma \ref{lem:reldar} restricts to $M$ to  Darboux charts.
\end{lemma}
\begin{proof}
It follows because   the charts in lemma  \ref{lem:reldar}   are obtained by applying a construction depending smoothly
on the center of the chart to obtain a number of equalities  for tensors and distributions at the origin, and then rescaling.
 Therefore when we restrict the charts to $M$ condition (1) in definition \ref{def:dcharts} holds. Conditions (2) and (3) follow
 because before rescaling  $D_x\oplus D^\perp_x$ is sent to $\mathbb{R}^{2n}\oplus \mathbb{R}$ and  $J_x$ to $J_0$.
  The Darboux  condition (equation (\ref{eq:darbcond})) holds by construction.
\end{proof}
\begin{lemma}\label{lem:triv} Let $\varphi_{k,x}\colon (\mathbb{C}^p,0)\rightarrow (P,x)$ be charts coming from lemma \ref{lem:reldar}. Then in $B(0,\rho k^{1/2})\subset \mathbb{C}^p$
 it is possible to fix a family of unitary trivializations of $\varphi^*_{k,x}L_{\Omega}^{\otimes k}$ with
 connection forms $A_{k,x}$, such that for all $k\gg 1$
\begin{enumerate}
\item $|A_{k,x}|_{g_0}\leq O(|z_k|)$, $|\mathrm{d}A_{k,x}|_{g_0}\leq O(1)$,  $|\mathrm{d}^jA_{k,x}|_{g_0}\leq O(k^{-1/2}),\;j\geq 2$.
\item ${A_{k,x}}_{\mid M}=\frac{1}{2}\sum_{i=1}^n(x_k^{2i-1}\wedge dx_k^{2i}-x_k^{2i}\wedge dx_k^{2i-1})$.
\end{enumerate}
\end{lemma}
\begin{proof}
By construction  $|\varphi_{k,x}^*k\omega|_{g_0}\leq O(1)$,
$|\mathrm{d}^j\varphi_{k,x}^*k\omega|_{g_0}\leq O(k^{-1/2})$, $j\geq
1$, on $B(0,\rho k^{1/2})$. Hence, we deduce the existence unitary
trivializations with connection forms $A'_{k,x}$ satisfying the
bounds of condition (1).

When we restrict the connection forms to $M$ they coincide with $A$ up to a exact 1-form $dF_{k,x}$
 defined on $\mathbb{R}^{2n+1}\times\{0\}$; its bounds are as in item (1) above, but on $\mathbb{R}^{2n+1}\times\{0\}$
  instead of on $\mathbb{C}^p$. We  extend it to $\mathbb{C}^p$ independently of the remaining coordinates
  and still denote it by $F_{k,x}$. It is always possible to find a unitary trivialization $\xi_{k,x}$ of
    $\varphi^*_{k,x}L_{\Omega}^{\otimes k}$ whose connection form is $A'_{k,x}+dF_{k,x}$. These
     trivializations give the desired result.
For simplicity we will denote the family by $A$ when there is no risk of confusion.
\end{proof}

Let $G$ be the $J$-complex distribution on $P$ that extends $D$. Given $\tau_k\in \Gamma(L^{\otimes k}_\Omega)$,
 the restriction of the covariant derivative of $\tau_k$ to $G$ will be denoted by
 $\nabla_G\tau_k\in \Gamma(G^*\otimes L^{\otimes k}_\Omega)$. Since $G$ is $J$-complex we can write
\[\nabla_G\tau_k=\bar{\partial}_G\tau_k+\partial_G\tau_k,\;  \bar{\partial}_G\tau_k\in \Gamma(G^{*0,1}\otimes L^{\otimes k}_\Omega),\;  \partial_G\tau_k\in \Gamma(G^{*1,0}\otimes L^{\otimes k}_\Omega).\]

\begin{lemma}\label{lem:relref}\quad
\begin{enumerate}
\item If $\tau_k\colon P\rightarrow  L^{\otimes k}_\Omega$ is an A.H. sequence then
 ${\tau_k}_{\mid M}\colon M\rightarrow L^{\otimes k}$ is also an A.H. sequence.
\item Moreover, the restriction of  a family of reference sections of $(L^{\otimes k}_\Omega,\nabla_k)\rightarrow (P,\Omega)$
 centred at the points of $M$ (as defined in \cite{Au01}) is a family of reference sections
 of $(L^{\otimes k},\nabla_k)\rightarrow (M,D,\omega)$.
\item If $\tau_k\colon P\rightarrow  L^{\otimes k}_\Omega$ is an A.H. sequence then $\bar{\partial}_G\tau_k\approxeq 0$.

\end{enumerate}
\end{lemma}
\begin{proof} We fix a family of A.H. charts adapted to $(M,G)$ and Darboux over $M$, and trivialize the bundles
 $L^{\otimes k}_\Omega$ as in lemma \ref{lem:triv}. Let $x^1_k,\dots,x_k^{2p}$ be the coordinates and write
  $\tau_{k,x}=f_{k,x}\xi_{k,x}$.

We  first observe that lemmas \ref{lem:localcheck1} and \ref{lem:localcheck2} for symplectic manifolds
also hold for the connection forms $A_{k,x}$ provided by lemma \ref{lem:triv}.
By  lemma \ref{lem:ahm} the restriction of the coordinates to $M$ are Darboux charts.
 We can apply lemma \ref{lem:localcheck1} for almost complex manifolds,
  bounds of order O(1), and the connection forms provided by lemma \ref{lem:triv}, to conclude that
   the partial derivatives of $f_{k,x}$ are bounded by $O(1)$ in the ball $B(0,\rho)\subset \mathbb{R}^{2p}$.
    In particular we get the same bounds if we only take into account the partial derivatives with respect to
    the variables $x_k^1,\dots,x_k^{2n+1}$ and restrict our attention to $B(0,\rho)\subset \mathbb{R}^{2n+1}$.
     Now if we apply back lemma \ref{lem:localcheck1} (this time for almost CR manifolds) we conclude
     that $|\nabla^j({\tau_k}_{\mid M})|_{g_0}\leq O(1)$, $\forall j\geq 0$, in $B(0,\rho)\subset \mathbb{R}^{2n+1}$
      and for all $x\in M$, the constants being independent of $x$. Therefore $|\nabla^j({\tau_k}_{\mid M})|_{g_k}\leq O(1)$,
      $\forall j\geq 0$, in all the points of $M$.

Lemma \ref{lem:localcheck2} for symplectic manifolds and the connections of lemma \ref{lem:triv} gives

\begin{equation}\label{eq:localineq}\left|\frac{\partial}{\partial x^\alpha_k}(\bar{\partial}_0+A_{k,x}^{0,1})f_{k,x}\right|_{g_0}\leq O(k^{-1/2})
\end{equation}
in $B(0,\rho)\subset \mathbb{R}^{2p}$. Let us consider the splitting $\mathbb{C}^n\times \mathbb{C}^{p-n}$.
 The operator  $\bar{\partial}_0+A_{k,x}^{0,1}$ and its derivatives  can be split into two pieces using it.
  We consider  the part involving $d\bar{z}_k^{1},\dots, d\bar{z}_k^n$, for which the above inequalities
   also hold, but now in $B(0,\rho)\subset \mathbb{R}^{2n+1}$.  Since the restriction of $A_{k,x}$ to
    $\mathbb{C}^n\times\mathbb{R}$ is $A$, the restriction to $M$ of the piece of
    $\bar{\partial}_0+A_{k,x}^{0,1}$ involving $d\bar{z}_k^1,\dots, d\bar{z}_k^n$ is
    the operator $\bar{\partial}_0+A^{0,1}$ of lemma \ref{lem:localcheck2}. Thus  we can
    apply this lemma (we already have the required bounds for the partial derivatives of $f_{k,x}$)
     to conclude $\bar{\partial}({\tau_k}_{\mid M})\approxeq 0$, and this proves item (1).

It is also easy to check that reference sections for $L^{\otimes k}_\Omega$  centred at the points
 of $M$ restrict to reference sections for $L^{\otimes k}$, and hence item (2) also holds.

To prove  $\bar{\partial}_G{\tau_k}\approxeq 0$ we use the previous ideas: equation  (\ref{eq:localineq}) and lemma \ref{lem:localcheck2} give

\[|\nabla^j\bar{\partial}_{\mathbb{C}^n}\tau_k|_{g_0}\leq O(k^{-1/2}),\; \forall j\geq 0\]
in $B(0,\rho)\subset \mathbb{R}^{2p}$, where $\bar{\partial}_{\mathbb{C}^n}$ is the part of $\bar{\partial}_0+A_{k,x}$
 involving $d\bar{z}_k^1,\dots, d\bar{z}_k^n$. The choice of A.H. charts adapted to $G$ and the bounds
  $|\nabla^j\tau_k|_{g_0}\leq O(1)$, $\forall j\geq 0$, easily imply
\[|\nabla^j(\bar{\partial}_{\mathbb{C}^n}\tau_k-\bar{\partial}_G{\tau_k})|_{g_0}\leq O(k^{-1/2}),\; \forall j\geq 0,\]
and therefore  $\bar{\partial}_G{\tau_k}\approxeq 0$.
\end{proof}

\begin{remark} Notice that item (3) in lemma \ref{lem:relref} is an assertion about a section defined on
$P$, and not on $M$ unlike in item (1).
\end{remark}

\subsection{Higher rank ample bundles} So far we have only considered approximately holomorphic theory
 for the sequence of line bundles $(L^{\otimes k},\nabla_k)\rightarrow (M,D,\omega)$, but there are
  obvious extensions for sequences of the form $E\otimes L^{\otimes k}$, where $E$ is any Hermitian
   bundle of rank $m$ with compatible connection. Regarding the local theory  the role of the reference
    sections is played by the reference frames $\tau_{k,x,1}^{\mathrm{ref}},\dots,\tau_{k,x,m}^{\mathrm{ref}}$,
     where each $\tau_{k,x,j}^{\mathrm{ref}}$ is an A.H. sequence with Gaussian decay with respect to $x$ and
     they are a frame of $E$ comparable to a unitary one in $B_{g_k}(x,\rho)$, $\rho>0$. Reference frames
      are constructed by tensoring  reference sections for $L^{\otimes k}$ with local unitary frames of $E$.

\section{Estimated transversality and finite, Whitney (A), approximate holomorphic stratifications}\label{sec:esttrans}

Let  $\tau_k$ be an A.H. sequence of sections of $L^{\otimes k}\rightarrow (M,D,\omega)$.
 Proposition \ref{pro:proposition1} for codimension two submanifolds is proved by pulling
 back the ${\bf 0}$ section of $L^{\otimes k}$. To obtain $W_k$ a 2-calibrated submanifold
  $\tau_k$ has to be transverse along $D$, so that $TW_k\cap D$ defines a codimension
  one distribution on $W_k$. Next, to make sure that $W_k\cap D$ is a  symplectic
  distribution the ratio $|\bar{\partial}\tau_k(x)|/|\partial\tau_k(x)|$ has to be smaller
   than 1; since $\nabla_D=\bar{\partial}+\partial$,  $\nabla_D\tau_k(x)$ has to be asked
    to be not only to be surjective but also to have norm greater than $O(k^{-1/2})$ (estimated transversality).

For each point $x$ we can use the reference sections to turn the local estimated  transversality
 problem along $D$ on $B_{g_k}(x,\rho)$, into an estimated transversality
  problem along $D_h$ for an A.H. sequence of functions
 \[F_{k,x}\colon B(0,\rho')\subset \mathbb{C}^n\times\mathbb{R}\rightarrow \mathbb{C},\]
  where $\tau_k\circ \varphi_{k,x}=F_{k,x}\cdot(\tau_{k,x}^{\mathrm{ref}}\circ \varphi_{k,x})$
  (more generally $\mathbb{C}^m$-valued functions for bundles of rank $m$).  Equivalently,
   we have to solve  an estimated transversality problem for a  1 real parameter family of A.H. functions
  \[F_{k,x}(\cdot, s_k)\colon B(0,\rho')\subset \mathbb{C}^n\rightarrow \mathbb{C}.\]
    This problem is known to have a solution  \cite{Au02,IMP01}.
The solution of the local transversality problem along $D_h$ will produce a new function
$F_{k,x}-u_{k,x}$, and therefore  a perturbation
\[\chi_{k,x}:=(-u_{k,x}\circ \varphi_{k,x}^{-1})\cdot\tau_{k,x}^{\mathrm{ref}}\] so that we obtain
 estimated transversality along $D$ for $\tau_k+\chi_{k,x}$ over the ball $B_{g_k}(x,\rho)$.
  But the reference section is supported in $B_{g_k}(x,\rho''k^{1/6})$, being the consequence
   that there will be interference among different local solutions.  However, unlike transversality,
    estimated transversality does behave well under addition, and in the presence of ``enough''
     local estimated transversality, Donaldson's globalization procedure gives global
     estimated transversality (see the proof of theorem \ref{thm:main2}).

\begin{definition}\label{def:esttrans} Let $(P,g)$ be a Riemannian manifold, $(E,\nabla)$
 a Hermitian bundle over it, and $Q_x$ a subspace of $T_xP$. We say that
 $\tau\colon P\rightarrow E$ is $\eta$-transverse to {\bf 0} at $x$ along
 $Q_x$, if either $|\tau(x)|\geq \eta$ or $\nabla_{Q_x}\tau(x)$ has a right inverse
 with norm bounded by $\eta^{-1}$.

If $Q$ is a distribution we say that $\tau$ is $\eta$-transverse along $Q$ to {\bf 0}
  if the above condition holds at all the points where $Q$ is defined. When $Q$ is the tangent bundle of a submanifold we also say that $\tau$ is $\eta$-transverse over the submanifold to {\bf 0}.

Let $(M,D,\omega)$ be a 2-calibrated manifold, $E_k:=E\otimes L^{\otimes k}$,  and
 $\tau_k\colon (M,g_k)\rightarrow (E_k,\nabla_k)$ a sequence of sections. We say
 that the sequence $\tau_k$ is uniformly transverse along $D$ to {\bf 0}  if
  $\eta>0$ exist such that $\tau_k$ is  $\eta$-transverse along $D$
   to {\bf 0}  for all $k \gg 1$.

For a symplectic manifold the definition of uniform transversality along a distribution $Q$
 (possibly the tangent bundle to a 2-calibrated submanifold) is analogous.
\end{definition}

It is possible to attain estimated transversality along $D$ using both the intrinsic
 and the relative point of view. Using the former what we do is (locally) solving
  transversality problems for 1-parameter families of  A.H. functions from
   $\mathbb{C}^n$ to $\mathbb{C}^m$. Regarding the latter we follow the ideas
   of J.-P Mohsen developed for contact manifolds, working in the symplectization
     $(M\times[-\epsilon,\epsilon],\Omega)$ and solving the estimated transversality
      problem for A.H. sections, but this time over $M$. Then we can use the following

\begin{lemma}\label{lem:contmohs} [\cite{Moh01}, second lemma in subsection 6.1]
 Let $(M,D,\omega)$ be a 2-calibrated manifold.  If in the symplectization
 $(M\times[-\epsilon,\epsilon],\Omega)$ we are able to find an A.H. sequence
 $\tau_k$ $\;\eta$-transverse over $M$ to {\bf 0}, then for any  constant
 $C$, $0<C<\sqrt{2}/2$, there exists $k_0(C)$ such that for any $k\geq k_0$
 the section ${\tau_k}_{\mid M}$ is  $C\eta$-transverse  along $D$ to {\bf 0}.
\end{lemma}

The proof is just an estimated version of the following elementary fact: if a $J_0$-complex linear
function $l\colon \mathbb{C}^n\times \mathbb{R}\rightarrow \mathbb{C}^m$ is surjective,
 then it has a surjective restriction to each complex hyperplane. Otherwise the kernel
  of the restriction -being complex- would have real dimension bigger than $2(n-m)+2$, and $l$ could not be surjective.

\subsection{Geometric reformulation of estimated transversality}

 We recall that in this section we deal with estimated  transversality along
 $D$ in a 2-calibrated manifold (intrinsic theory), or with estimated transversality
  over a 2-calibrated submanifold $M$ inside a symplectic manifold $P$ (relative theory).
   Sometimes we might refer to both situations as transversality along a distribution
    $Q$ in the Riemannian manifold $P$.

As remarked in the previous subsection, for sequences of 1-parameter families of
A.H. functions $F_{k,x}(\cdot, s_k)\colon B(0,\rho)\subset \mathbb{C}^n\rightarrow \mathbb{C}^m$
 one can achieve estimated transversality, and thus the use of reference frames allows us
  to get local estimated transversality along $D$ to the ${\bf 0}$ section of very ample vector
   bundles $E_k$. More generally, one expects to be able to attain estimated transversality
    along $D$ to sequences of submanifolds  $S_k\subset E_k$ of very ample vector bundles,
    where the $S_k$ locally look like the zero section of a trivial vector bundle: more precisely,
     the sequence of submanifolds should be locally defined by functions $f_k\colon U_k\subset E_k\rightarrow \mathbb{C}^l$,
       $S_k\cap U_k=f_k^{-1}({\bf 0})$, which are approximately holomorphic with respect to to the
        almost CR structure in the total space of  the bundles $(E_k,\nabla_k)\rightarrow (M,D,J,g_k)$
         induced by the one on $M$, the connection, and the Hermitian metric on $E_k$, so that
\[f_k\circ \tau_k\circ \varphi_{k,x} \colon B(0,\rho)\subset \mathbb{C}^n\times \mathbb{R}\rightarrow \mathbb{C}^m\]
 is an  A.H. sequence of functions (or a weaker  property that ensures this last condition).
  That should allow us  to find an A.H. perturbation
 \[\chi_{k,x}\colon B_{g_k}(x,\rho''k^{1/6})\rightarrow E_k\]
 so that the A.H. sequence
 \[f_k\circ (\tau_k+\chi_{k,x})\colon B_{g_k}(x,\rho)\rightarrow \mathbb{C}^m\] is uniformly transverse  along $D$
 to ${\bf 0}$. Finally, we should make sure that this implies enough estimated
 transversality along $D$ to $S_k$ for the sequence of sections
 \[\tau_k+\chi_{k,x}\colon B_{g_k}(x,\rho)\rightarrow E_k\] to make Donaldson's globalization procedure work.

In the relative context $\tau_k\colon P\rightarrow E_k$ the  estimated transversality
 problem over $M\subset P$  to the $\bf 0$ section  has the same difficulty as the usual
  estimated transversality problem to the ${\bf 0}$ section (this is the work of J.-P.
   Mohsen \cite{Moh01}, section 5). Thus, one expects this principle to  be valid in
   the case of relative estimated transversality to more complicated strata $S_k$.

To give a global definition of what transversality to a submanifold $S\subset E$ is,
 we need to recall  a more geometric definition of estimated transversality along
 a distribution $Q$, together with the following concepts.

 \begin{definition}\label{def:Mangle} Let $W$ be a vector space with non-degenerate inner
 product so that for any $u,v \in W$ we can compute the (unoriented) angle  $\angle(u,v)$. Given  $U \in
 Gr(p,W)$ and $V\in Gr(q,W)$, $p,q>0,$ the maximal angle  of  $U$  and $V$,  $\angle_\mathrm{M}(U,V)$,
is defined  as follows:
 \begin{eqnarray*}
 \angle_\mathrm{M}(U,V)&:=& \mathrm{max}_{u\in U\backslash \{0\}}\mathrm{min}_{v\in V\backslash \{0\}}\angle(u,v).
 \end{eqnarray*}
 \end{definition}

 In general the maximal angle is not symmetric, but when  $p=q$ it has symmetry and defines
 a distance in the corresponding Grassmannian (see \cite{MPS02}).

The minimum angle between transverse complementary subspaces is defined as the minimum
angle between two non-zero vectors, one on each subspace. An extension of this notion for
 transverse subspaces with non-trivial intersection is:

\begin{definition}\label{def:mangle}{(Definition 3.3. in \cite{MPS02})} Using the notation of definition
\ref{def:Mangle}, $\angle_\mathrm{m}(U,V)$ -the minimum angle
between $U$ and $V$ non-void  subspaces of $W$- is defined as follows:
\begin{itemize}
 \item  If $\dim U+\dim V<\dim W$, then  $\angle_\mathrm{m}(U,V):=0$.
 \item If the intersection is non-transverse, then $\angle_\mathrm{m}(U,V):=0$.
 \item If the intersection is transverse, we consider the orthogonal to the
  intersection and its intersections  $U_c$ and $V_c$ with $U$ and $V$ respectively. We define
  $\angle_\mathrm{m}(U,V):=\mathrm{min}_{u\in U_c\backslash \{0\}}\mathrm{min}_{v \in V_c \backslash \{0\}}\angle(u,v)$.
 \end{itemize}

 The minimum angle is symmetric.
\end{definition}

The most important property relating maximal and minimal angle is:

\begin{proposition}{(Proposition 3.5 in \cite{MPS02})}\label{pro:minmax}
For non-void subspaces  $U,V,W$ of $\mathbb{R}^n$ the following inequality holds:
\[ \angle_\mathrm{m}(U,V)\leq  \angle_\mathrm{M}(U,W)+\angle_\mathrm{m}(W,V).\]
\end{proposition}

We will also be using the following
\begin{lemma}[lemma 3.8 in \cite{MPS02}\label{lem:functang}] Let $U,V$ be non-zero subspaces
 of $\mathbb{R}^n$ and let  $h\colon U\rightarrow V^\perp$ be the
  projection from $U$ with respect to the decomposition $\mathbb{R}^n=V\oplus V^\perp$. If
  $h$ has a right inverse $\theta$ satisfying $|\theta|<\eta^{-1}$ then $\angle_m (U,V)>\eta$.
\end{lemma}

Let  $\tau\colon P\rightarrow E$ be a section of a Hermitian bundle with connection and $Q$
 a distribution on $P$. Let us denote  the pullback of $Q$ to $E$ by $\hat{Q}$. Let
 $\mathcal{H}$ be the horizontal distribution associated to the linear connection and
 let $\mathcal{H}_Q$ denote its intersection with $\hat{Q}$. Finally let $T_Q\tau$ denote
  the intersection of the tangent bundle of the graph of $\tau$ with $\hat{Q}$.

\begin{lemma}\label{lem:comang}
 There exists a constant $C>0$ determined by upper bounds on $|\nabla_Q\tau(x)|,|\tau(x)|$ such
 that:
\begin{enumerate}
\item  If $\nabla_Q\tau(x)$ has a right inverse with norm bounded by $\eta^{-1}$ then
$\angle_\mathrm{m}(\mathcal{H}_Q, T_Q\tau)\geq C^{-1}\eta $ (the angle measured in $\hat{Q}_{\tau(x)}$).
\item If $\angle_\mathrm{m}(\mathcal{H}_Q, T_Q\tau)\geq \eta$ then $\nabla_Q\tau(x)$
has a right inverse with norm bounded by  $
(C\mathrm{sin}\eta)^{-1}$.
\end{enumerate}
\end{lemma}
\begin{proof}
Let us assume $Q=TP$. The vector space $T_{\tau(x)}E=\mathcal{H}_{\tau(x)}\oplus T^vE_x$ is endowed
 with the direct sum metric. We compose with an isometry preserving the direct sum structure  so
 that  $\mathcal{H}_{\tau(x)}\oplus T^vE_x$ becomes $\mathbb{R}^a\oplus\mathbb{R}^b$ with the Euclidean metric.
Let  $h\colon T\tau(x)\rightarrow  \mathbb{R}^b$ be the orthogonal projection. By lemma
\ref{lem:functang} applied to $U=T\tau(x)$  and $V= \mathbb{R}^a\times\{0\}=\mathcal{H}_{\tau(x)}$,
 if $h$ has a right inverse $\theta$ with $|\theta|\leq \eta^{-1}$ then  $\angle_\mathrm{m}(\mathcal{H}_{\tau(x)},
  T\tau(x))\geq \eta$.

By definition  $\nabla\tau(x)\colon T_xP\rightarrow T^vE_x=E_x$ is the composition
$h\circ  d\tau(x)$, with the differential $d\tau(x)\colon T_xP\rightarrow T\tau(x)$, which is an isomorphism.
Now if $\theta'$ is a right inverse for $\nabla\tau(x)$, $|\theta'|\leq \eta^{-1}$,  then
$d\tau(x)\circ \theta'$ is a right inverse for $h$ with norm bounded by $|d\tau(x)|\eta^{-1}$.
 Thus, by lemma \ref{lem:functang} $\angle_\mathrm{m}(\mathcal{H}_{\tau(x)}, T\tau(x))\geq |d\tau(x)|^{-1}\eta$.

Conversely, the projection $h$ has always a right inverse $\theta$ of minimum norm. Let us define
$W:= T\tau(x)\cap \mathcal{H}_{\tau(x)}$ and $U_c:=T\tau(x)\cap W^\perp$. If we compose
$\theta$  with the orthogonal projection $ T\tau(x)\rightarrow U_c$, we obtain a right
inverse $\hat{\theta}$ for $h_{\mid U_c}$ such that $|\theta|=|\hat{\theta}|$. If now
 $\angle_\mathrm{m}(\mathcal{H}_{\tau(x)}, T\tau(x))\geq \eta$ then the equation involving
  inequalities of lemma 3.8 in \cite{MPS02} implies
  \begin{equation}\label{eq:inverseprojbound}
   |\hat{\theta}|\leq (\mathrm{sin}\eta)^{-1},
   \end{equation}
   and therefore $d\tau(x)^{-1}\circ\theta$ is a right inverse for $\nabla\tau(x)$ with norm
    bounded by $|d\tau(x)|^{-1}(\mathrm{sin}\eta)^{-1}$.

In the case $Q\neq TP$ we fix an isometry sending  $(\mathcal{H}_Q,\mathcal{H})$ at $\tau(x)$ to
 $(\mathbb{R}^{a'}\times \{0\},\mathbb{R}^{a})$ with the Euclidean metric, and apply the above
  arguments to $\mathbb{R}^{a'}\oplus\mathbb{R}^b$.

Note that we have $C=|d_Q\tau(x)|$, with $d_Q\tau(x)$ the restriction of $d\tau(x)$ to $Q_x$.
 Observe that a bound for $|d_Q\tau(x)|$ can be obtained from upper bounds for $|\tau(x)|$ and $|\nabla_Q\tau(x)|$.
\end{proof}

\begin{remark}\label{rem:compmimang} In the definition of minimum angle $\angle_\mathrm{m}(U,V)$, when $U,V$ are not
 complementary we work with the intersections in $(U\cap V)^\perp$ where we can apply the
  usual notion of minimum angle for complementary subspaces. Instead of $(U\cap V)^\perp$
  one might choose any other  subspace $W$ complementary to $U\cap V$ to give a different
  notion of minimum angle. In certain situations this is a good strategy because  there are
   natural complementary subspaces  available. It is easy to see that the new notion of
   minimum angle is comparable to the one of definition \ref{def:mangle}, and the comparison is
   given by multiplying by a constant depending only on  $\angle_\mathrm{m}(U\cap V,W)$
   (there is no ambiguity since these are complementary subspaces). Actually, those new
    notions depending on the complementary coincide with the one given in \ref{def:mangle},
     but for a new metric, which is comparable to the Euclidean one in terms of $\angle_\mathrm{m}(U\cap V,W)$
      (very much as it happened with the isomorphism $d_Q\tau(x)$ in the previous lemma).
\end{remark}

We need a second result relating angles and intersections.

\begin{lemma}\label{lem:minmax2}
Let $U,V,W$ be linear subspaces of $\mathbb{R}^n$ such that  $\angle_\mathrm{m}(V,W)\geq \gamma>0$.
 Let $\angle_\mathrm{M}(U,V)\leq \delta$. Then there exists $C(\gamma,\mathrm{dim}V,n)>0$ such that
\[ \angle_\mathrm{M}(U\cap W,V\cap W)\leq C\delta.\]
\end{lemma}
\begin{proof}
For each $u\in U\backslash \{0\}$, we have  $\angle(u,V)=\angle(u,h(u))$, where
$h\colon \mathbb{R}^n\rightarrow V$ is the orthogonal projection. We consider a complementary space to
 $V$ possibly different from $V^\perp$: because  $\angle_\mathrm{m}(V,W)\geq \gamma>0$ the dimension of $W$ is greater or equal than the codimension of $V$, and the intersection of $V$ and $W$ is transverse. As a consequence any subspace (of $W$) complementary to $V\cap W$ in $W$ is also complementary to $V$ in $\mathbb{R}^n$. We let $V_W$ be the orthogonal to $V\cap W$ in $W$,
 and we define $h_W\colon \mathbb{R}^n\rightarrow V$ to be the projection along $V_W$ (whose restriction to $W$ is the orthogonal projection onto $V\cap W$).
  It follows that $\angle(u,h_W(u))\leq C\angle(u,h(u))=C\angle(u,V)$, and by construction
   if $u\in U\cap W$ then $\angle(u,h_W(u))=\angle(u,V\cap W)$.
\end{proof}

Let $S\subset E$ be a submanifold in the total space of the vector bundle $E$ over either
a 2-calibrated or a symplectic manifold, transverse to the fibers.  Let $\hat{g}$ be the
 metric in $E$ induced by the connection, the bundle metric, and the metric $g$ in the base. The submanifold  might
  not have a tubular neighborhood of positive radius. If we assume $S$ to be in a compact region -as it will be the case in our applications- then the problem comes from the behavior near its boundary $\partial S=\bar{S}\backslash S$. Thus a reasonable extension of definition \ref{def:esttrans} to our non-linear setting must deal separately with points close to $\partial S$ and with the other points of $S$.

 \begin{definition}\label{def:farfrombound} Given  $\bar{\eta}>0$ the points of $S$ $\bar{\eta}$-far
  from (resp. $\bar{\eta}$-close to) the boundary are those points in $S$ at $\hat{g}$-distance
   of $\partial S$ greater or equal (resp. smaller) than $\bar{\eta}>0$.  For any $\eta>0$
    -typically much smaller than $\bar{\eta}$- we define $\mathcal{N}_S(\eta,\bar{\eta})$
    to be those points that can be joined to a point  $\bar{\eta}$-far from the boundary
    by a geodesic arc normal to $S$ and of length smaller or equal than $\eta$.
 \end{definition}

 We now define the distribution   $T^{||}S$ at the points of $\mathcal{N}_S(\eta,\bar{\eta})$
  by parallel transport of $TS$ along the geodesics normal to $S$, starting at the points
  $\bar{\eta}$-far from the boundary of $S$.

$T^{||}S$ plays the role of $\mathcal{H}$.  We use the notation $T^{||}_QS:=T^{||}S\cap \hat{Q}$.

\begin{definition}\label{def:esttrans2} $\tau$ is  $(\eta,\bar{\eta})$-transverse along $Q$  to $S$
 at $x$ if either (i) $\tau(x)$ misses the union of $S$ with $\mathcal{N}_S(\eta,\bar{\eta})$,
  or (ii) $\tau(x)$ enters in $\mathcal{N}_S(\eta,\bar{\eta})$ so that $\angle_\mathrm{m}(T_Q\tau,T^{||}_QS)\geq \eta$
   at $\tau(x)$, or (iii) $\tau(x)$ intersects $S$ at the points $\bar{\eta}$-close to the boundary with
    $\angle_\mathrm{m}(T_Q\tau,T_QS)\geq \bar{\eta}$ at $\tau(x)$.

Uniform transversality of $\tau_k$ along $Q$ to $S_k$ is defined as $(\eta,\bar{\eta})$-transversality
 for some $\eta,\bar{\eta}> 0$ and for all $k \gg 1$.

\end{definition}

Conditions on a sequence of submanifolds $S_k$ of complex codimension $l$ (or more generally
 on stratifications) can be imposed, so that local estimated transversality along $Q$ of
 $\tau_{k,x}$ at the points of $B_{g_k}(x,\rho)$ to the points of $S_k$  far from $\partial S_k$,
 is equivalent to  estimated  transversality along $Q$  of a related $\mathbb{C}^l$-valued
 function to ${\bf 0}$ (lemma \ref{lem:localchartrans}).

 We will consider stratifications $\mathcal{S}=(S_k^a), a\in
A_k$,  which are (i) finite in the sense that $\#(A_k)$ must
 be bounded independently of $k$, and (ii) the boundary of each strata $\partial{S}_k^b=\bar{S}_k^b\backslash S_k^b$
  will be the union of the strata of smaller dimension

 \[\partial S_k^b=\bigcup_{a<b} S_k^b.\]

\begin{definition}\label{def:stratification}
Let $E_k=E\otimes L^{\otimes k}\rightarrow (M,D,J,g_k)$  and let
$(S_k^a)_{a \in A_k}$ be finite  stratifications of  $E_k$ whose strata are transverse to the
 fibers. Let   $r\in \mathbb{N}$, $r\geq 2$. The sequence of strata is  Whitney $C^r$-approximately holomorphic
$(C^r$-A.H.) if for any bounded open set $U_k$ of the total space of $E_k$  and any $\epsilon>0$,
  constants  $C_{\epsilon},\rho_{\epsilon}>0$ only depending on
$\epsilon$ and on the size of   $U_k$  -but not on  $k$- can be found, so that for any point
 $y \in U_k$ in a strata $S_k^a$ for which $d_{\hat{g}_k}(y,\partial S_k^a)>\epsilon$, there
  exist complex valued functions   $f_1,\dots,f_l$ such that  $B_{\hat{g}_k}(y,\rho_{\epsilon})\cap
S_k^a$ is given $f_1=\dots=f_l=0$, and the following properties hold:
\begin{enumerate}
\item (Uniform transversality to the fibers + transverse comparison) The restriction of
$df_1\wedge \cdots\wedge df_l$ to $T^vE_k$  is bounded from below by  $\rho_\epsilon$.
\item (Approximate holomorphicity along the fibers) The restriction of the function
$f=(f_1,\dots,f_l)$ to each fiber is  $C^r$-A.H.$(C_{\epsilon})$.
\item (Horizontal approximate holomorphicity + holomorphic variation of the restriction to
 the fiber + estimated variation of the restriction to the fiber) For any $\lambda, k$, and $\tau$
$C^r$-A.H.$(\lambda)$ local section of  $E_k$ with image cutting
 $B_{\hat{g}_k}(y,\rho_{\epsilon})$, $f_j\circ \tau$ is
$C^r$-A.H.$( \lambda C_{\epsilon})$. Moreover, if $\theta$ is a local $C^r$-A.H.$(\lambda)$
 section of  $\tau^*T^vE_k$, $df_\tau(\theta) $ is
$C^r$-A.H.$(\lambda C_{\epsilon})$.
\item (Estimated Whitney's condition (A)) For each $\eta>0$ small enough, there exists $\delta(\eta)>0$
 such that $\forall y \in S_k^b$ at distance smaller than $\delta$ of
$S_k^a\subset \partial S_k^b$,  $\angle_\mathrm{M}(T^{||}S_k^a, TS_k^b)$ at $y$ is bounded by
 $\eta$.
\end{enumerate}
\end{definition}


\begin{remark} If we give the corresponding definition using as base space an almost complex
manifold instead of an almost CR manifold,  we almost recover the definition 3.2 in \cite{Au01}
 (our condition (4) is a bit weaker).
\end{remark}

Condition (1) is equivalent to the strata have minimum angle with the fibers bounded from below.
 We just try to mimic the picture of the {\bf 0} section with respect to the fibers of a vector bundle,
 in which case we even have orthogonality.

Conditions (2) and (3) guarantee that if $\tau_k\colon M\rightarrow E_k$ is A.H., then the corresponding
 $\mathbb{C}^l$-valued function  to be made transverse to {\bf 0} is A.H.

Recall that for a stratification $\mathcal{S}$ of some $\mathbb{R}^N$, a stratum $S^b$ satisfies Whitney's
condition (A) if for every converging sequence $x_n\rightarrow x$, $x_n\in S^b$,
 $x\in S^a\subset \partial S^b$, so that $T_{x_n}S^b$ is converging,
 the limit contains $T_xS^a$. Condition (4) is an estimated Whitney's condition (A).

 \begin{definition}\label{def:transstr} Let  $\mathcal{S}$ be  as in definition  \ref{def:stratification}
  (over either a 2-calibrated or a symplectic manifold). Then $\tau_k$ is uniformly transverse
   along $Q$ to  $\mathcal{S}$ if there exists strictly positive numbers  $(\eta_a,\bar{\eta}_a)$  for all $a\in A_k$ such that:
\begin{enumerate}
\item  For all  $a\in A_k$ and for all $k \gg 1$  $\tau_k$ is $(\eta_a,\bar{\eta}_a)$-transverse along $Q$ to $S^a_k$.
\item For each $b$,  $\bigcup_{a<b}\mathcal{N}_{S_k^a}(\eta_a,\bar{\eta}_a)$  contains
 the points of ${S}^b_k$ $\bar{\eta_b}$-close to $\partial S_k^b$.
\end{enumerate}
\end{definition}

Now that we have the notion of uniform transversality of a sequence of sections to an appropriate
 stratification, we need tools to relate it with local uniform transversality for sequences of (related) functions.

\begin{lemma}\label{lem:localchartrans} Let $S^a_k$ be a sequence of strata as those in the stratifications
 of definition \ref{def:stratification} for the base space $P$ either an almost CR manifold
(intrinsic theory) or an almost complex manifold (relative theory).
Let $\epsilon >0$ and $0<\eta\ll \epsilon$. Let $y\in E_k$ be a
point in the stratum   $\epsilon$-far from the boundary, and let
$f=(f_1,\dots,f_l)$ be the corresponding local $\mathbb{C}^l$-valued
function defining the stratum in $B_{\hat{g}_k}(y,\rho_\epsilon)$.
Let $\tau_k$ be a section of $E_k$ whose graph enters in
$B_{\hat{g}_k}(y,\rho_\epsilon)$. Then there exist  constants
$\rho'(\epsilon,\eta,|\tau_k|)$,
$C(\epsilon,|\nabla_Q\tau_k|,|\tau_k|),C'(\epsilon,|\nabla_Q\tau_k|,|\tau_k|)>0$
such that:
\begin{enumerate}
\item  If $\angle_\mathrm{m}(T_Q\tau,T^{||}_{Q}{S^a})\geq \eta$ in $B_{\hat{g}_k}(y,\rho_\epsilon)$,
 then $d_Q(f\circ\tau)$ has a right inverse with norm bounded by $(C\mathrm{sin}(\eta/2))^{-1}$
in $B_{\hat{g}_k}(y,\rho')$.

\item If $d_Q(f\circ\tau)$ has a right inverse with norm bounded by $\eta^{-1}$ in $B_{\hat{g}_k}(y,\rho_\epsilon)$,
 then $\angle_\mathrm{m}(T_Q\tau,T^{||}_QS^a)\geq {C'}^{-1}\eta$ in $B_{\hat{g}_k}(y,\rho')$.
\end{enumerate}
\end{lemma}
\begin{proof}
By simplicity we omit the subindices for the sections $\tau_k$, the bundles, and the strata.

Let us assume $\angle_\mathrm{m}(T_Q\tau,T^{||}_{Q}{S^a})\geq \eta$.

\emph{Step1:} Show the existence of $\rho'(\epsilon,\eta,|\tau|)>0$ such
that $\angle_\mathrm{m}(T_Q\tau,\mathrm{Ker} df\cap \hat{Q})\geq \eta/2$ in $B_{\hat{g}_k}(y,\rho')$.

 According to proposition \ref{pro:minmax} (proposition 3.5 in \cite{MPS02})
\[ \angle_\mathrm{m}(T^{||}_{Q}{S^a},T_Q\tau)\leq \angle_\mathrm{M}(T^{||}_{Q}{S^a},\mathrm{Ker} df\cap \hat{Q})+
\angle_\mathrm{m}(\mathrm{Ker} df\cap \hat{Q},T_Q\tau),\]
and therefore  we need to prove the existence of $\rho'>0$  so that  in $B_{\hat{g}_k}(y,\rho')$
\begin{equation}\label{eq:angle1}
\angle_\mathrm{M}(T^{||}_{Q}{S^a},\mathrm{Ker} df\cap \hat{Q})\leq \eta/2.
\end{equation}
Condition (1) in definition \ref{def:stratification} implies
$\angle_\mathrm{m}(\mathrm{Ker} df,\hat{Q})\geq \gamma(\epsilon)$. If we find $\rho'>0$
 such that  in $B_{\hat{g}_k}(y,\rho')$
\begin{equation}\label{eq:angle2}
\angle_\mathrm{M}(T^{||}{S^a},\mathrm{Ker} df)\leq C(\gamma(\epsilon))^{-1}\eta/2,
\end{equation}
we can  apply  lemma \ref{lem:minmax2}, where $U=T^{||}{S^a}$, $V=\mathrm{Ker} df$,
 $W=\hat{Q}$, to conclude that equation  (\ref{eq:angle1}) holds.

Equation  (\ref{eq:angle2}) is proven using  appropriate charts.
The situation we are trying to mimic is that of a locally trivialized vector bundle and
we measure the maximal angle between the parallel copies of the {\bf 0} section
(here the leaves of $\ker df$) and $\mathcal{H}$ (here $T^{||}S^a$).

Due to the bounds in definition \ref{def:stratification} we can find
a chart $\Phi_y\colon \mathbb{R}^a\rightarrow
B_{\hat{g}_k}(y,\rho_\epsilon)$ such that in  $B(0,\rho'')\subset
\mathbb{R}^a$ (i) the metrics $g_0$ and $\Phi_y^*\hat{g}_k$ (that we
write $\hat{g}_k$ if it is clear that we work in the chart) are
comparable,  and  the Christoffel symbols of $\hat{g}_k$ are bounded
by $O(1)$  (the bounds being uniform on $k,y$), and (ii) the
foliation  $\mathrm{Ker} df$ is sent  to the foliation  $\mathbb{R}^{a-2l}$.
In $B(0,\rho'')\subset \mathbb{R}^a$  the stratum $S$ becomes
$\mathbb{R}^{a-2l}\times\{0\}$ and tubular neighborhoods for
$\hat{g}_k$ and $g_0$ are comparable. At any point $q$ in the
neighborhood, a vector in $u\in T^{||}S$ is the result of parallel
translating (with $\hat{g}_k$) a vector $v$ in
$\mathbb{R}^{a-2l}\times\{0\}$ over $y'\in
\mathbb{R}^{a-2l}\times\{0\}$ along the corresponding
$\hat{g}_k$-geodesic. Since the Christoffel symbols are bounded,
$\angle(u,v)$ is bounded by $e^{\Gamma t}-1$, $\Gamma>0$. So by
decreasing $t$, the distance of $q$ to $S$, we bound the maximal
angle by $C(\gamma)^{-1}\eta/2$. Therefore the final radius $\rho'$
depends on $\eta$, on $\epsilon$ (because $C(\gamma)$ depends on
$\epsilon$), and on how $g_0$ and $\hat{g}_k$ are related (to order
one). This final relation depends on $f$ (and hence on $\epsilon$)
and on the metric $\hat{g}_k$ (and hence on  $|\tau|$).

\emph{Step 2:} Show that $\angle_\mathrm{m}(T_Q\tau, \mathrm{Ker}
df\cap \hat{Q})\geq \eta/2$ implies that $d_Q(f\circ\tau)$ has a
right inverse with norm bounded by
$(C(\epsilon,|\nabla_Q\tau|,|\tau|)\mathrm{sin}(\eta/2))^{-1}$.

 The proof of item (2) in lemma \ref{lem:comang} implies that the orthogonal projection
 $h\colon T_Q\tau \rightarrow {(\mathrm{Ker}df\cap \hat{Q})}^{\perp}$
  has a right inverse with norm bounded by ${(\textrm{sin}(\eta/2))}^{-1}$ (equation (\ref{eq:inverseprojbound})).
   Let $V_E$ denote the orthogonal in the fiber $T^vE$ of $(\mathrm{Ker}df\cap \hat{Q})\cap T^vE$.
Due to  condition (1)  in definition \ref{def:stratification}, this is a subspace complementary to $\mathrm{Ker}df\cap \hat{Q}$
 and such that $\angle_\mathrm{m}(V_E, \mathrm{Ker}df\cap \hat{Q})$ is
  bounded from below in terms of $\rho_\epsilon$, and hence in terms of $\epsilon$.

Let $h_E\colon T_Q\tau \rightarrow V_E$ be the  projection along $\mathrm{Ker}df\cap \hat{Q}$.
It follows that there is a  constant $C_1(\epsilon)^{-1}>0$ and a right inverse for $h_E$
 with norm bounded by $C_1(\epsilon)^{-1}{(\textrm{sin}\eta/2)}^{-1}$. We now define

\[h''=df\circ h_E\circ d_Q\tau \colon Q\rightarrow \mathbb{C}^l.\]
By construction $h''=d_Q(f\circ \tau)$. Condition (1) about the restriction of $df$
 to the fiber implies the existence of  a right inverse for $h''$  with norm bounded by
$|d_Q\tau|^{-1}C_2(\epsilon)^{-1}C_1(\epsilon)^{-1}{(\textrm{sin}\eta/2)}^{-1}$.
Therefore, $d_Q(f\circ \tau)$ has a right inverse with norm bounded
by $C(\epsilon,|d_Q\tau|)\textrm{sin}(\eta/2)^{-1}$ in
$B_{\hat{g}_k}(y,\rho'(\epsilon,\eta,|\tau|))$, this proving item (1).

Conversely, if $d_Q(f\circ \tau)$ has a right inverse
$B_{\hat{g}_k}(y,\rho_\epsilon)$ with norm bounded by $\eta^{-1}$,
step 2 above implies that  $h_E\circ d_Q\tau$ has a right inverse
with norm bounded by $({C'_1}(\epsilon)2\eta)^{-1}$.

Item (1) in  lemma \ref{lem:comang} gives

\[\angle_\mathrm{m}(T_Q\tau,\ker df\cap \hat{Q})\geq C'(\epsilon,|\nabla_Q\tau|,|\tau|)2\eta,\]
and combined with step 1 we conclude
\[\angle_\mathrm{m}(T_Q\tau,T^{||}_QS^a)\geq C'(\epsilon,|d_Q\tau|)\eta \;\; \mathrm{in} \;\; B_{\hat{g}_k}(y,\rho(\epsilon,\eta,|\tau|)).\]

Observe that the constants $C,C'$ grow very large as $\epsilon$ and $\eta$ tend to zero.
\end{proof}

\begin{remark}\label{rem:weakstrat} The previous lemma does not involve almost
complex structures at all. Hence it also holds for arbitrary Hermitian bundles, sections,
and strata which fulfill condition (1) in definition  \ref{def:stratification}.
\end{remark}

Using appropriate choices of complementary subspaces to get a bound
from below for certain minimal angles -as noticed in remark
\ref{rem:compmimang}- we can prove the following

 \begin{lemma}\label{lem:reltransv} Let  $\mathcal{S}=(S_k^a)_{a \in A}$ be a sequence of approximately
   holomorphic stratifications as in definition  \ref{def:stratification}. Assume that the sequence
    $\tau_k$ is uniformly transverse to  $\mathcal{S}$ along a distribution
     $Q$ whose dimension is greater of equal than the codimension of the strata, and that the
      uniform bounds $|\tau_k|,|\nabla\tau_k|_{g_k}\leq O(1)$ hold. Then for each $a\in A$, $\tau_k^{-1}(S_k^a)$
       is a subvariety of $M$ uniformly transverse to $Q$.
\end{lemma}
 \begin{proof} We must prove that for a sequence of points $x(k)$ in $\tau_k^{-1}(S_k^a)$ we have
 \begin{equation}\label{eq:minanginvim0}
 \angle_m(T_{x}\tau_k^{-1}(S_k^a),Q)\geq \gamma>0
 \end{equation}
for all $k \gg 1$ independently of the points.

Denote $\tau_k(x)=q$. We claim that equation (\ref{eq:minanginvim0})
would follow from
\begin{equation}\label{eq:minanginvim}
\angle_m(\tau_{k*}T_{x}\tau_k^{-1}(S_k^a),\tau_{k*}Q)\geq \gamma'>0,
\end{equation}
where the angle is measured in $T\tau_k(x)$ with the induced metric. The reason is
 that the bound on $|\tau_k|$ implies that the metric in $E_k$ is comparable to
  the product metric given by any trivialization by reference frames (and using
   on each factor the Hermitian metric in the fiber and $g_k$ coming from the base).
    Then we use the bound on $|\nabla \tau_k|$ to conclude that in this product
    metric $\angle_m(T\tau_k(x),T^vE_k(q))\geq \delta_1>0$, where $T^vE_k$ is
    the tangent space to the fiber. Hence, our claim follows.

We can rewrite equation (\ref{eq:minanginvim}) as
\begin{equation}\label{eq:minanginvim2}
\angle_m(T_Q\tau_k(x),T\tau_k\cap TS_k^a(q))\geq \gamma'>0.
\end{equation}
Our second claim is that
\begin{equation}\label{eq:minangstrat}
\angle_m(TS_k^a(q),T^vE_k(q))\geq \delta_2>0.
\end{equation}
Indeed, this follows from condition (1) in definition
\ref{def:stratification} if we are in a point $\bar{\eta}$-far from the
boundary of $S^a_k$. For points $\bar{\eta}$-close, we use the estimated
Whitney's condition (A) together with proposition \ref{pro:minmax}
to prove  equation (\ref{eq:minangstrat}). Since $T^vE_k\subset
\hat{Q}$, we also conclude
\begin{equation}\label{eq:minangstrat2}
\angle_m(TS_k^a(q),\hat{Q})\geq \delta_3>0.
\end{equation}

We will reinterpret equation (\ref{eq:minangstrat2}) by
choosing a suitable complementary space to $TS_k^a\cap \hat{Q}(q)$
which is not its orthogonal $W$ (see remark \ref{rem:compmimang}).
Let $W_1\subset \hat{Q}$ (resp. $W_2\subset \hat{Q}$)  be the intersection
 of $T\tau_k(x)$ (resp. $TS_k^a(q)$) with the orthogonal of
  $TS_k^a\cap T\tau_k\cap \hat{Q}(q)$ inside $\hat{Q}$, and
  let $W_3$ be the intersection of $T\tau_k(x)$ with the orthogonal of $\hat{Q}$.
From  $\angle_m(T\tau_k(x),T^vE_k(q))\geq \delta_1$ we obtain
  $\angle_m(T\tau_k(x),\hat{Q})\geq \delta_1$, and by hypothesis
  $\angle_m(T_QT\tau_k(x),T_QS^a_k(q))\geq \delta_4>0$. Both inequalities
   imply that $W':=W_1 \oplus W_3$ can be used instead of $W$. By construction
    $W'\cap \hat{Q}=W_1$,  so from equation (\ref{eq:minangstrat2}) we conclude
\begin{equation}\label{eq:minangstrat3}
\angle_m(W'\cap TS_k^a(q),W_1)\geq \delta_5>0.
\end{equation}

Notice as well that to compute equation (\ref{eq:minanginvim2}) we have to intersect
 the corresponding vector subspaces with the orthogonal of $TS_k^a\cap T\tau_k\cap \hat{Q}(q)$
  inside $T\tau_k(x)$. From what we have seen, we can rather choose as complementary space
   $W'$. Since $W'\cap T_Q\tau_k(x)=W_1$ and $W'\cap (T\tau_k\cap TS_k^a(q))=W'\cap TS_k^a(x)$,
    we have to compute the left hand side of  equation (\ref{eq:minangstrat2}), so the result follows.
\end{proof}

In particular the following corollary is deduced:

\begin{corollary}\label{cor:transD}  Let  $\mathcal{S}=(S_k^a)_{a \in A}$ be a sequence of A.H.
 stratifications over the 2-calibrated manifold $(M,D,\omega)$ as in definition
  \ref{def:stratification}. Assume that the A.H. sequence  $\tau_k$ is uniformly
  transverse to $\mathcal{S}$ along $D$. Then for each  $a\in A_k$, $\tau_k^{-1}(S_k^a)$
   is either empty -if the codimension of  $S_k^a$ is bigger than the dimension of
    $D$ (or $M$)- or a subvariety uniformly transverse to $D$.

For a symplectic manifold,  transversality along the directions of a (compact) subvariety
  $N$ implies  that  either (i)  $\tau_k^{-1}(S_k^a)$ is at $g_k$-distance of $N$ bounded
   from below or (ii) it is a subvariety (at least defined in a $g_k$-neighborhood of $N$) uniformly transverse to $N$.
\end{corollary}

If we analyze the proof of lemma \ref{lem:reltransv}, corollary
\ref{cor:transD} for $2$-calibrated manifolds is equivalent to
saying that uniform transversality along $D$ implies uniform
transversality over $M$ (along $TM$). The converse is also true, extending
therefore Mohsen's relative transversality result to appropriate
sequences of stratifications.

\begin{corollary}\label{cor:fulltrans} Let  $\mathcal{S}=(S_k^a)_{a \in A_k}$ be a sequence
 of A.H. stratifications over the 2-calibrated manifold $(M,D,\omega)$ as in definition
  \ref{def:stratification}.
Assume that the A.H. sequence $\tau_k$ is uniformly transverse to $\mathcal{S}$ (over  $M$),
 for suitable constants $(\eta_a,\bar{\eta}_a)$, $a\in A_k$. Then $\tau_k$ is also uniformly
 transverse  along $D$ to $\mathcal{S}$.
\end{corollary}
\begin{proof} By induction we can assume that  $\tau_k$ is uniformly transverse along $D$ to $S_k^a$, for every $a<b$.
Let $q\in S^b_k$, with $\tau_k(x)=q$, $\bar{\eta}'$-close to $\partial S_k^b$.  We want to show
\[\angle_\mathrm{m}(T_D\tau_k(x),T_DS_k^b(q),)\geq \bar{\eta}',\]
and we will do it by applying for some index $a\in A_k$ the inequality

\begin{equation}\label{eq:ineq1}
\angle_\mathrm{m}(T_D\tau_k(x),T^{||}_DS^a_k(q))\leq  \angle_\mathrm{M}(T^{||}_DS_k^a(q),T_DS_k^b(q))+\angle_\mathrm{m}(T_D\tau_k(x),T_DS_k^b(q)).
\end{equation}
If $\bar{\eta}'$ is small enough condition (2) in definition \ref{def:transstr} implies
 the existence of an index $a\in A_k$  such that $q\in \mathcal{N}_{S_k^a}(\eta_a,\bar{\eta}_a)$.
  If we apply induction we conclude $\angle_\mathrm{m}(T_D\tau_k(x),T^{||}_DS^a_k(q))\geq \eta_a$,
so we only need to make
\[\angle_\mathrm{M}(T^{||}_DS_k^a(q),T_DS_k^b(q))\ll\eta_a.\]
This is done using
lemma \ref{lem:minmax2} with $U=T^{||}{S^a}(q)$, $V=TS^b(q)$, $W=\hat{D}$. We need to check
\begin{align}\label{eq:whitang1}
\angle_\mathrm{M}(T^{||}S_k^a(q),TS_k^b(q))&\ll\eta_a,\\\label{eq:whitang2}
\angle_\mathrm{m}(TS^b_k(q),\hat{D})&\geq \gamma.
\end{align}
Equation (\ref{eq:whitang1}) follows by the estimated Whitney's condition by taking $\bar{\eta}'$
small enough; equation (\ref{eq:whitang2}) uses again the inequality of proposition  \ref{pro:minmax}
\[\angle_\mathrm{m}(\hat{D},T^{||}S^a_k(q))\leq  \angle_\mathrm{M}(T^{||}S_k^a(q),TS_k^b(q))+\angle_\mathrm{m}(\hat{D},TS_k^b(q)),\]
together with   $\angle_\mathrm{m}(\hat{D},T^{||}S^a_k(q))\geq 2\gamma$ (by condition (1)
in definition  \ref{def:stratification}) and  equation (\ref{eq:whitang1}).

So far we deduced some $\bar{\eta}'$-transversality only at the points $\bar{\eta}'$-close to the boundary of $S_k^b$.
Now let us assume that for some $\eta>0$, $\angle_\mathrm{m}(T\tau_k(x),T^{||}S_k^b(q))\geq \eta$
 in the tubular neighborhood $\mathcal{N}_{S_k^b}(\eta,\bar{\eta}')$ (here comes the requirement
  on the constants controlling the transversality, i.e. in those points $\bar{\eta}'$-far
   from the boundary   we need to make sure that $\angle_\mathrm{m}(T\tau_k(x),T^{||}S_k^b(q))$ is uniformly bounded from below).
If $\tau_k(x)\in \mathcal{N}_{S_k^b}(\eta,\bar{\eta}')$ then by lemma \ref{lem:localchartrans} $\eta$-transversality
implies  $\eta'$-transversality to ${\bf 0}$ of the function
$f\circ\tau_k\colon B_{g_k}(x,\rho')\rightarrow \mathbb{C}^l$. From the  approximate
 holomorphicity of the composition  $f\circ \tau_k$, for all $k \gg 1$ a result analogous
  to lemma \ref{lem:contmohs} grants $\frac{\sqrt{2}}{3}\eta'$-transversality along
  $D$, which again by lemma  \ref{lem:localchartrans} gives $\eta''$-transversality  along $D$
   to $S_k^b$ (we suppose $\eta''\leq \eta$).

Therefore, it follows that $\tau_k$ is $(\eta'',\bar{\eta}')$-transverse  along $D$
to $S_k^b$.
\end{proof}

 \section{Pseudo-holomorphic jets}\label{sec:pholjets}

 The main applications of the theory of approximately holomorphic geometry for 2-calibrated
  manifolds are deduced from the existence of generic rank $m$ linear systems.

Let us assume that $(M,\mathcal{D},J)$ is a  Levi-flat CR manifold  and $L\rightarrow M$ a
 positive CR line bundle. Let $\underline{\mathbb{C}}^m\rightarrow M$ denote the trivial
  (and trivialized) bundle of rank $m$ endowed with the trivial connection.

\begin{definition}\label{def:rgenseccr}  A CR section $\tau\colon M\rightarrow \underline{\mathbb{C}}^{m+1}\otimes L$ (or
 a rank $m$ linear system of $L$) is $r$-generic if its zero set $B$ is a CR submanifold of the
  expected dimension, and the projectivization  $\phi\colon M\backslash B\rightarrow \mathbb{C}\mathbb{P}^m$
   is a leafwise $r$-generic holomorphic map, i.e. when restricted to each leaf it is
   transverse to the Thom-Boardman stratification of the bundle of holomorphic $r$-jets
   of holomorphic maps from the leaf to $\mathbb{C}\mathbb{P}^m$.
\end{definition}

The proof of the existence of $r$-generic linear systems (possibly of large enough powers of $L$)  is the main subject of \cite{Ma05b}.

The strong transversality property for a CR function $\phi\colon
M\rightarrow \mathbb{C}\mathbb{P}^m$ to be $r$-generic is as
follows: we consider $\mathcal{J}^r_{CR}(M,\mathbb{C}\mathbb{P}^m)$
 the bundle of CR $r$-jets (of foliated holomorphic $r$-jets) of CR maps from $M$ to
 $\mathbb{C}\mathbb{P}^m$. This bundle admits a CR Thom-Boardman stratification $\mathbb{P}\Sigma$,
  which restricts to each leaf to the corresponding holomorphic Thom-Boardman stratification.
   A CR function $\phi$ is $r$-generic if and only if its CR $r$-jet
   $j^r_{CR}\phi\colon M\rightarrow \mathcal{J}^r_{CR}(M,\mathbb{C}\mathbb{P}^m)$
    (which by definition is the foliated holomorphic $r$-jet) is  transverse along $\mathcal{D}$ to $\mathbb{P}\Sigma$.

Assume that our CR submanifold embeds holomorphically in some complex
manifold $P$ and that $\mathcal{D}$ extends to a holomorphic
foliation integrating the complex distribution  $G$.  There is a canonical submersion
$p_G\colon
\mathcal{J}^r(P,\mathbb{C}\mathbb{P}^m)\rightarrow
\mathcal{J}^r_G(P,\mathbb{C}\mathbb{P}^m)$ from
holomorphic $r$-jets to foliated ones. The foliated Thom-Boardman
stratification $\mathbb{P}\Sigma\subset
\mathcal{J}^r_G(P,\mathbb{C}\mathbb{P}^m)$ restricts over
$M$ to the CR Thom-Boardman stratification $\mathbb{P}\Sigma$ of
$\mathcal{J}^r_{CR}(M,\mathbb{C}\mathbb{P}^m)$. Let us denote the
pullback  ${p_G}^{-1}(\mathbb{P}\Sigma)$  by
$\mathbb{P}\Sigma^G$.

 For any holomorphic function $\phi\colon P\rightarrow \mathbb{C}\mathbb{P}^m$  it is an elementary fact that
   $j^r_G\phi\in \Gamma(\mathcal{J}^r_G(P,\mathbb{C}\mathbb{P}^m))$ -the
    holomorphic $r$-jet along $G$- is transverse along $G$ to $\mathbb{P}\Sigma$
     at the points of $M$, if and only if $j^r\phi\in \Gamma(\mathcal{J}^r(P,\mathbb{C}\mathbb{P}^m))$
      is transverse along $G$ to $\mathbb{P}\Sigma^G$ at the points of $M$. By the results
       of the previous section, this is equivalent to being transverse  over $M$ to $\mathbb{P}\Sigma^G$.

To obtain an $r$-generic linear system there is an additional complication coming from the base locus.
 We first need to make sure that   $\tau\colon P\rightarrow \underline{\mathbb{C}}^{m+1}\otimes L$
 is transverse over $M$ to the zero section, and then solve the $r$-genericity problem for the
 projectivization (in a compact region of $P\backslash \tau^{-1}({\bf 0})$). Instead of working
 first with the section $\tau$ and then with the projectivization, following ideas of D. Auroux \cite{Au01}
  we restate the whole issue  as a unique transversality problem over $M$ for the
  pseudo-holomorphic $r$-jet extension of $\tau$, a section of  a vector bundle
  $\mathcal{J}^r(\underline{\mathbb{C}}^{m+1}\otimes L)$. The advantage is that we work with vector bundles and we can use the module structure of sections.

\subsection{The integrable case}\label{ssec:integra}

Let $E \rightarrow P$ be a Hermitian bundle over a complex manifold with compatible connection $\nabla$,   whose curvature verifies  $F^{0,2}_\nabla=0$. The total space of the bundle is a complex
  manifold (theorem $2.1.53$ in \cite{DK90}) and there is a  notion of holomorphic section and
  hence of holomorphic $r$-jet. The space of $r$-jets has natural charts obtained out of holomorphic
   coordinates in the base and a holomorphic trivialization of the bundle. They provide a local
   identification of the holomorphic $r$-jets with   $\mathcal{J}^r_{n,m}$, the usual $r$-jets for
    holomorphic maps from  $\mathbb{C}^n$ to $\mathbb{C}^m$.

Let  $\partial_0 $ be the Cauchy-Riemann operator defined (locally) using the canonical structure
 $J_0$ in the base (the chart) and the trivial connection $\mathrm{d}$ in $\underline{\mathbb{C}}^m$.
  The connection on the fiber bundle can be used to give a different notion of local holomorphic
  $r$-jet (in principle chart dependent) by just considering the operator $\partial_\nabla$:  if the
   connection matrix in the trivialization is $A_x=A_x^{1,0}$, then the coupled 1-jet of a
   holomorphic section $\tau$ is defined to be $(\tau,\partial_0\tau+A_x\tau)$). Higher order coupled
    jets are constructed by induction using the connection  induced by the flat metric and $\nabla$.

Observe that locally for the above choice of coordinates and trivialization of the bundle, both the
 usual $r$-jets and coupled $r$-jets fill the bundle
 \[(\sum_{j=0}^{r}(T^{*1,0}{\mathbb{C}^n})^{\odot j})\otimes\underline{\mathbb{C}}^m=\mathcal{J}^{r}_{m,n},\]
  where $\odot $ stands for the symmetric part of the tensor product and
  $(T^{*1,0}{\mathbb{C}^n})^{\odot 0}\otimes\underline{\mathbb{C}}^m$ for $\underline{\mathbb{C}}^m$.
   This is due to the existence through any point of $E$ of  holomorphic frames tangent to the
    horizontal distribution of the connection, together  with the vanishing  $F^{2,0}_\nabla$
    (the latter implying that $\textrm{d}A$ and its derivatives are symmetric tensors when evaluated on (1,0)-vectors).

For Levi-flat CR manifolds the local model for the pseudo-holomorphic jets to be introduced is
the following: the base space is  $(\mathbb{C}^n\times \mathbb{R},J_0,g_0)$ (or rather a ball
of Euclidean radius $\rho>0$), the bundle is assumed to be trivialized by a CR frame and
 the curvature is of type (1,1).  The bundle of CR $r$-jets is denoted by $\mathcal{J}^r_{D_h,n,m}$
  (foliated holomorphic $r$-jets along $D_h$); its fiber over each point is that of $\mathcal{J}^r_{n,m}$.
   There is an obvious notion of  CR coupled $r$-jet. The hypothesis on the trivialization and on
   the curvature  imply that they are also symmetric, so they fill the bundle
   $\mathcal{J}^r_{D_h,n,m}=\mathcal{J}^r_{n,m}\times \mathbb{R}$.

Using Darboux charts and suitable trivializations this model will be achieved in an approximate
way in the theory for 2-calibrated manifolds.

There is a final local model we wish to introduce that would appear in K\"ahler manifolds $P$
 with a holomorphic foliation integrating a complex distribution $G$. Locally, we have holomorphic  coordinates
 $\mathbb{C}^g\times \mathbb{C}^{p-g}$ with $G$ sent to $\mathbb{C}^g$ (which integrates into the foliation
  with leaves $\mathbb{C}^g\times\{\cdot\}$), and we  work with foliated coupled jets along the
   leaves of  $\mathbb{C}^g$. The corresponding bundle of coupled foliated $r$-jets is denoted by
     $\mathcal{J}^r_{\mathbb{C}^g,p,m}$. It coincides  with  $\mathcal{J}^r_{g,m}\times \mathbb{C}^{p-g}$.
       Transversality problems for this bundle will be transferred to transversality problems
        in $\mathcal{J}^r_{p,m}$, so we  need no further analysis of its properties, though we
         will be interested at some point in studying the natural submersion
         $\mathcal{J}^r_{p,m}\rightarrow \mathcal{J}^r_{\mathbb{C}^g,p,m}$.  This local model
          is achieved in an approximate way in a symplectic manifold (with compatible almost complex structure and metric)
           with a $J$-complex distribution $G$ -not necessarily integrable-, by using approximate holomorphic charts adapted to $G$.

\subsection{Pseudo-holomorphic jets}

Denote sequence $E\otimes L^{\otimes k}\rightarrow (M,D,\omega)$ by $E_k$.  We define the bundles
 \[\mathcal{J}^r_D E_k:=(\sum_{j=0}^{r}{(D^{*1,0})}^{\odot j})\otimes E_k, \]
 where $\odot $ stands for the symmetric part of the tensor product of complex vector bundles.
They carry Hermitian vector bundle metrics induced by ${g_k}_{\mid D}$, the one on $E_k$, and the symmetrization map
 \begin{equation}\label{eq:symmap}{\mathrm{sym}}_j\colon {(D^{*1,0})}^{\otimes j}\rightarrow {(D^{*1,0})}^{\odot j}.
\end{equation}
The Levi-Civita connection induces a connection on $D^*$
 (using the metric to see $D^*\hookrightarrow T^*M$ and then projecting  $T^*M\rightarrow D^*$)
  and therefore in $D^{*1,0}$ (using the splitting $D^{*1,0}+D^{*0,1}$); combined with the
  connection on $E_k$ and the symmetrization map  they define connections $\nabla_{k,r}$.
 The total spaces $\mathcal{J}^r_D E_k$ also carry metrics constructed in the usual fashion out of the metric in the base, the connection, and the vector bundle Hermitian metric.

The definition of pseudo-holomorphic $r$-jets along $D$ (or just pseudo-holomorphic $r$-jets) for a
 sequence  $E_k$ of Hermitian vector bundles is given by induction (see \cite{Au01}). Let $\tau_k$ be a sequence of A.H. sections of $E_k$. By definition $j^0_D\tau_k=\tau_k$.  Let   $j^{r-1}_D\tau_k \in \mathcal{J}^{r-1}_DE_k$ be the ($r-1$)-jet of $\tau_k$. It has homogeneous
 components of degrees $0,1,\dots,r-1$. We will denote the homogeneous component of degree
  $j\in \{0,\dots,r-1\}$ by $\partial^j_{\mathrm{sym}}\tau_k\in \Gamma(({D^{*1,0}})^{\odot j}\otimes E_k)$.
The connection $\nabla_{k,r-1}$ is actually a direct sum of connections defined on the direct
summands $({D^{*1,0}})^{\odot j}\otimes E_k$, $j=0,\dots,r-1$. For simplicity and if there is
 no risk of confusion we will use the same notation for the restriction of $\nabla_{k,r-1}$ to each of the summands.
The restriction of $\nabla_{k,r-1}\partial^{r-1}_{\mathrm{sym}}\tau_k$ to $D$ defines a section
 $\nabla_{k,r-1,D}\partial^{r-1}_{\mathrm{sym}}\tau_k\in \Gamma(D^*\otimes({D^{*1,0}})^{\odot r-1}\otimes E_k)$.
  For each $x\in M$ it is a form on $D$ with values in the complex vector space
   $({D^{*1,0}})^{\odot r-1}\otimes E_k$. Therefore we can consider its (1,0)-component
    $\partial \partial^{r-1}_{\mathrm{sym}}\tau_k\in \Gamma(D^{*1,0}\otimes({D^{*1,0}})^{\odot r-1}\otimes E_k)$.
     By applying the symmetrization map $\mathrm{sym}_r$ of equation (\ref{eq:symmap}) we obtain
      $\partial^r_\mathrm{sym}\tau_k\in \Gamma(({D^{*1,0}})^{\odot r}\otimes E_k)$.

 \begin{definition}\label{def:psjets} Let  $\tau_k$ be a section of  $(E_k,\nabla_k)$.
 The pseudo-holomorphic $r$-jet  $j^r_D\tau_k$ is a section of the bundle
 $\mathcal{J}^r_DE_k=(\sum^{j=0}_{r}({D^{*1,0}})^{\odot j})\otimes E_k$ defined out of the ($r-1$)-jet by the formula
 $j^r_D\tau_k:=(j^{r-1}_D\tau_k,\partial^r_\mathrm{sym}\tau_k)$.
 \end{definition}

\begin{remark} The previous definition incorporates the fact that the degree $r$ and ($r-1$) homogeneous
 components of the $r$-jet are  symmetrization of the pseudo-holomorphic  1-jet of $\partial^{r-1}_{\mathrm{sym}}\tau_k$;
  then we have to add the homogeneous components of lower degree. Actually, we could have equally
  defined $j^r_D\tau_k$ by taking the symmetrization of the pseudo-holomorphic 1-jet of $j^{r-1}_D\tau_k$
   (because this gives the homogeneous components of degree $1,\dots,r$) and then adding $\tau_k$,
   the degree zero homogeneous component.
 \end{remark}

\begin{remark} The pseudo-holomorphic $r$-jets are useless for our purposes for low values of $k$. We are interested in having a notion of $r$-jet of an A.H. sequence which in approximately holomorphic coordinates and for suitable local trivializations of $E_k$, is  as close as possible to the local coupled
  holomorphic $r$-jet  defined in $\mathbb{C}^n \times \mathbb{R}$ using  $J_0$ and the flat metric
   (introduced in subsection \ref{ssec:integra}). As $k$ grows large and
    due to the proximity between $g_k,J$ and $J_0,g_0$,  in $B(0,\rho)\subset \mathbb{C}^n \times \mathbb{R}$ we will see that
     the norm of the difference at any order between the two notions of $r$-jet is bounded by $O(k^{-1/2})$.
\end{remark}

For a symplectic  manifold  $(P,\Omega)$ with a $J$-complex distribution $G$ the bundle of
 pseudo-holomorphic $r$-jets along $G$ will be defined to be
 \[\mathcal{J}^r_G E_k:=(\sum_{j=0}^{r}({G^{*1,0})}^{\odot j})\otimes E_k.\]
 We have a canonical projection $p_G\colon\mathcal{J}^rE_k\rightarrow \mathcal{J}^r_G E_k$.
 We also use the splitting   $TP=G\oplus G^{\perp}$ to see  $\mathcal{J}^r_G E_k$ as a subbundle of
   $\mathcal{J}^rE_k$; hence every section of $\mathcal{J}^r_G E_k$ can be seen as a section of
     $\mathcal{J}^rE_k$.
      To define the pseudo-holomorphic $r$-jet along $G$ we use the same induction procedure as in
      the definition of pseudo-holomorphic $r$-jets along $D$, but either before or after symmetrizing
       we project   $T^{*1,0}P\rightarrow G^{*1,0}$ (or even before taking the (1,0)-component
        we project   $T^*P_{\mathbb{C}}\rightarrow G^*_{\mathbb{C}}$); the result of either choice is the same.

 Once approximately holomorphic coordinates have been fixed we have a canonical pointwise $(J_0-J)$-complex
  linear identification
 \begin{align}\nonumber
 T\mathbb{C}^n & \rightarrow D\\ \nonumber
 \frac{\partial}{\partial x_k^i} &\mapsto  \frac{\partial}{\partial x_k^i}+a_i \frac{\partial}{\partial s_k}\\\label{eq:linidentf}
   \frac{\partial}{\partial y_k^i} &\mapsto J\left( \frac{\partial}{\partial x_k^i}+a_i\frac{\partial}{\partial s_k}\right).
  \end{align}
The inverse of its dual is a $(J_0-J)$-complex bundle map
\begin{equation}\label{eq:identif}
\varpi_{k,x}\colon T^{*1,0}\mathbb{C}^n\rightarrow D^{*1,0}.
\end{equation}
It should be stressed that this identification  is only important  in the ball of  some $g_k$ radius  $\rho>0$,
 the  region where our computations have to be more accurate (in order to obtain local estimated transversality).
   There, for some constant $\gamma>0$
\begin{equation}\label{eq:bundlemap}
|\varpi_{k,x}|_{g_0}\leq \gamma,\;\;|\varpi_{k,x}^{-1}|_{g_0}\leq \gamma\;\;\mathrm{and}\;\;|\mathrm{d}^j\varpi_{k,x}|_{g_0}\leq  O(k^{-1/2}),\; \forall j\geq 1.
 \end{equation}
The Gaussian decay of the reference sections will take care of what happens out of these balls. We also notice
 that by  writing   $dz_k^i$ we will mean   $\varpi_{k,x}(dz_k^i)$.

 Let us assume that we have also fixed a family of reference sections of $\tau_{k,x}^{\mathrm{ref}}\in \Gamma(L^{\otimes k})$.
  Using any local unitary basis of $E$ (with  bounds uniform on $x$) together with the reference sections,
   we have a family of trivializations $\tau_{k,x,j}^{\mathrm{ref}},j=1,\dots,m$, of $E_k$ in the balls
   $B_{g_k}(x,\rho)$ for all $x$ and for all $k$ large enough. The A.H. coordinates and the associated
    bundle maps $\varpi_{k,x}$ provide a local basis  $dz_k^1,\dots,dz_k^n$  of $D^{*1,0}$.  We  obtain
     a family of trivializations of $\mathcal{J}^r_DE_k$ about any point as follows:  for
     $I=(i_0,i_1,\dots,i_n)$, with  $1\leq i_0\leq m$,
 $0\leq i_1+\cdots +i_n \leq r$, we set
  \begin{equation}\label{eqn:canframe}
  \mu_{k,x,I}:={dz_k^1}^{\odot i_1}\odot \cdots \odot {dz_k^n}^{\odot i_n}\otimes\tau^{\mathrm{ref}}_{k,x,i_0}.
  \end{equation}

\begin{definition}\label{def:frame} A family of sequences $\tau_{k,x,I}\colon M\rightarrow E_k$,   is called
 a family of holonomic frames if:
\begin{enumerate}
\item They are A.H. sections with Gaussian decay w.r.t to $x$.
\item There exist $\rho,\gamma >0$ such that in the balls $B_{g_k}(x,\rho)$ and for all point and all
 $k$ large enough the sequences $j^r_D\tau_{k,x,I}\colon M\rightarrow \mathcal{J}^r_DE_k$ define a
  frame which is $\gamma$-comparable to $\mu_{k,x,I}$ in the following sense: if we write $j^r_D\tau_{k,x,I}$
  in the basis $\mu_{k,x,I}$, for  the corresponding matrix $M_{k,x}$ we have
\[ |M_{k,x}|_{g_0}\leq \gamma, \; |M_{k,x}^{-1}|_{g_0}\leq \gamma.\]
\end{enumerate}
\end{definition}

One checks that the notion of holonomic reference frame does not depend either on the fixed approximately
 holomorphic coordinates, or in the chosen reference sections of $E_k$ to define $\mu_{k,x,I}$. Only
  the constants involved in the definition change.

In this situation there is still a weak point. The main goal is to construct sections whose
 pseudo-holomorphic $r$-jets are transverse to certain stratifications. For that we need the
  pseudo-holomorphic $r$-jets to be A.H. sections of the bundles $\mathcal{J}^r_DE_k$ (resp.
  $\mathcal{J}^rE_k$  for symplectic manifolds with $J$-complex distribution $G$), so that we
   can apply the transversality results from approximately holomorphic theory (to be proved in
    section \ref{sec:mainthm}).  We intend to use holonomic reference frames  defined as follows:
     if  $I$ is one of the  ($n+1$)-tuples introduced before we set
\begin{equation}\label{eqn:holonframe}
\nu_{k,x,I}:=j^r_D\tau_{k,x,I}^{\mathrm{ref}},\;\; \mathrm{where}\;\; \tau_{k,x,I}^{\mathrm{ref}}:={(z_k^1)}^{i_1}\cdots {(z_k^n)}^{i_n}\tau_{k,x,i_0}^{\mathrm{ref}}\in \Gamma(E_k).
\end{equation}
In the K\"ahler case and due to the presence of curvature (see \cite{Au02}),  the coupled jets
  are not anymore holomorphic sections of $ \mathcal{J}^r_{n,m}$ with respect to the  complex structure
   induced by the connection. Similarly, the frames $\nu_{k,x,I}$  fail to be families of holonomic
    frames because the sections are not approximately holomorphic if $r\geq 1$.
This difficulty is overcome by introducing a new almost complex structure (a new connection)
in  $\mathcal{J}^r_DE_k$ (resp. $\mathcal{J}^rE_k$). This is the content of the following proposition
 whose proof is given in \ref{sec:newconn}.

\begin{proposition}\label{pro:perthol} The sequence  $\mathcal{J}_D^r E_k\rightarrow (M,D,J,g_k)$ -which
is very ample for the connections  $\nabla_{k,r}$ previously described- admits new connections  $\nabla_{k,H_r}$ such that:
\begin{enumerate}
\item $\nabla_{k,r}-\nabla_{k,H_r} \in D^{*0,1}\otimes \mathrm{End}(\mathcal{J}_D^r E_k)$. Hence,
if in order to compute the pseudo-holomorphic jets (definition \ref{def:psjets}) we use the
connections $\nabla_{k,H_r}$ instead of $\nabla_{k,r}$, then the result is the same.
\item Let us denote the curvatures of $\nabla_{k,H_r}$ and $\nabla_{k,r}$  by $F_{k,H_r}$ and
$F_{k,r}$ respectively. Then  $F_{k,H_r}\approxeq F_{k,r}$ and hence  $(\mathcal{J}_D^r E_k,\nabla_{k,H_r})$
 is a very ample sequence.

\item If $\tau_k \colon M\rightarrow E_k$ is a $C^{r+h}$-A.H. sequence of sections, then
$j^r_D\tau_k\colon M \rightarrow \mathcal{J}_D^rE_k$ is a  $C^h$-A.H. sequence of sections
 for the connections  $\nabla_{k,H_r}$.
\end{enumerate}

In the integrable model $(E,\nabla)\rightarrow (\mathbb{C}^n\times\mathbb{R},D_h,J_0,g_0)$, with
$E=L_1\oplus\cdots \oplus L_m$, we can introduce new connections  $\nabla_{H_r}$ (here there is no dependence in $k$, since distribution, (almost) complex structure, and metric are the standard ones). If the curvature  $F_i$ of each line bundle   $L_i$, $i=1,\dots,m$,
 restricted to the leaves is of type (1,1) and has constant components with respect to the coordinates
 $z_1,\dots,z_n$, then the restrictions to each leaf of the curvatures $F_{H_r}$ and $F_{r}$
 (item (2) above)  coincide. As a consequence the new almost CR structure in the total space of
  $\mathcal{J}^r_{D_h,n,m}$ induced by $\nabla_{H_r}$ is also integrable (the foliation does
   not vary, just the leafwise complex structure). Also if  $\tau$ is a CR section ($\mathbb{C}^m$-valued
   function), then the coupled CR jet  is a CR section of $(\mathcal{J}^r_{D_h,n,m},\nabla_{H_r})$.

In the case of $(P,\Omega)$ symplectic with a $J$-complex distribution $G$, analogous results hold
for $\mathcal{J}^rE_k$ and for the integrable model.

\end{proposition}

As we said we postpone the proof until \ref{sec:newconn}, but we introduce
the formula for the connection.

Let $\sigma_k=(\sigma_{k,0},\sigma_{k,1})$ be a section (maybe local) of $\mathcal{J}_D^1 E_k$. We
  define \[\nabla_{H_1}(\sigma_{k,0},\sigma_{k,1})= (\nabla\sigma_{k,0},\nabla
\sigma_{k,1})+ (0,-F^{1,1}_D \sigma_{k,0}),\]
 where $F^{1,1}_D
\sigma_{k,0}\in D^{*0,1}\otimes D^{*1,0}\otimes E_k$ (see \cite{Au02}).

\begin{remark}\label{rem:metriccomparison} The approximate equality  $F_{H_1,k}\approxeq
F_{k}$ has useful consequences. Assume for simplicity $E_k=L^{\otimes k}$. Fix approximately holomorphic coordinates and trivialize the line bundle so that the connection form is $A$ (equation
  (\ref{eqn:connform})). Then in the local frame $(1,0)\otimes \tau_k,(0,dz_k^1)\otimes \tau_k,\dots, (0,dz_k^n)\otimes \tau_k$ $\mathcal{J}^1_DL_k$ and over $B(0,\rho)\subset \mathbb{C}^n\times\mathbb{R}$,
   the connection matrix of $\nabla_{k,H_1}$ is,
    up to summands bounded (at any order) by $O(k^{-1/2})$

\begin{equation}\nonumber
\left| \begin{array}{cccc}
 A & -\frac{1}{2}d\bar{z}_k^1 & \cdots & -\frac{1}{2}d\bar{z}_k^n \\
 0 & A & \cdots & 0\\
 & & \ddots & \\
 0& 0& \cdots & A
 \end{array}
 \right|
 \label{delta}
\end{equation}
In particular we have a uniform control on the new metric of the total space of the bundles
 $\mathcal{J}^1_DL_k$ (resp. $\mathcal{J}^1L_k$). In a similar manner this uniform control also
  holds for the bundles  $\mathcal{J}^r_DE_k$ (resp. $\mathcal{J}^rE_k$).
A useful outcome is that if we have  a sequence of  stratifications $\mathcal{S}$ such
 that for a choice of approximate holomorphic coordinates and reference frames, in the associated
 local basis  $\mu_{k,x,I}$ of equation (\ref{eqn:canframe}) the strata  $S_k^a$ are given by equations
  (functions) that do not depend neither on $k$ nor on $x$, then the different bounds associated
   to the strata (basically those of the local functions defining them) will not depend on $k$
    and $x$ (because we can compute them for the corresponding model with the Euclidean metric elements).
\end{remark}

\section{The linearized Thom-Boardman stratification }\label{TBA}

For the very ample sequences $E_k$ there is an easy sufficient condition for a sequence of
stratifications to be finite, Whitney (A), and approximately holomorphic.

Let us denote by $\mathsf{T}$ the group of translations of $\mathbb{C}^n\times \mathbb{R}$
 (resp. $\mathbb{C}^p$ in the relative case).

 \begin{lemma}\label{lem:trivializationweak}
 Let $(S_k^a)_{a \in A}$ be a sequence of stratifications of $E_k\rightarrow (M,D,\omega)$
 such that, for a choice of approximately holomorphic coordinates and approximately holomorphic
  trivialization, it is sent to $(S^a)_{a\in A}$, a fixed CR finite, Whitney (A) stratification
   of $\underline{\mathbb{C}}^m\rightarrow \mathbb{C}^n\times\mathbb{R}$ transverse to the fibers.
 Then the sequence  $(S^a_k)_{a \in A}$ is as in definition \ref{def:stratification}.

Conversely, from a Whitney (A) CR stratification of  $\underline{\mathbb{C}}^m \rightarrow \mathbb{C}^n\times\mathbb{R}$
 transverse to the fibers and invariant under  the action of  $\mathsf{T}\times Gl(m,\mathbb{C})$
  (or  $\mathsf{T}\times\mathbb{C}^*$), using the local identifications of $E_k$ with
   $\underline{\mathbb{C}}^m$ furnished by A.H. coordinates and A.H. trivializations, it
   is possible to induce an approximately holomorphic sequence of finite, Whitney (A) stratifications of  $E_k$.
\end{lemma}

\begin{proof} Recall that we are interested in constructing A.H. sequences of sections transverse
to $(S_k^a)_{a \in A}$; in particular this sections will be uniformly bounded. Therefore, for each
 $k,x$ we can work in the subset
 $B(0,\rho)\times B(0,R)\subset (\mathbb{C}^n\times\mathbb{R})\times \mathbb{C}^m =\underline{\mathbb{C}}^m$,
  for some $R>0$. Let  $f$ be a function defining locally a stratum $S^a$, which by hypothesis
  can be chosen to be CR.  Condition (1) in  definition  \ref{def:stratification} holds trivially
   for the model $S$ and therefore for $(S_k^a)_{a \in A}$, because when we compare the Euclidean
    metric and $\hat{g}_k$  we get the same inequalities as in condition (1) in definition  \ref{def:dcharts}.

Since the model stratification is Whitney  (A) and we work in
a compact region, Whitney's condition (A) implies the estimated
Whitney's condition (A) for the Euclidean metric and hence for
$\hat{g}_k$.

Let $\hat{J}_0$ be the leafwise holomorphic structure associated to the canonical CR structure of
$\underline{\mathbb{C}}^m = (\mathbb{C}^n\times\mathbb{R})\times \mathbb{C}^m $ and let $\hat{D}_h$
 denote the foliation by complex hyperplanes. Since the local function  $f$ defining $S^a$  is CR,
   in particular it is fiberwise holomorphic, and this proves condition (2) in  definition  \ref{def:stratification}.

Let $(\hat{D},\hat{J},\hat{g}_k)$ be the almost CR structure on $B(0,\rho)\times B(0,R)$ induced by the
 one on $E_k$. In order to prove condition (3) it suffices to check that $f$ is A.H. with respect to the this almost
  CR structure. We are going to  slightly modify the induced almost CR structure: instead of  $\hat{D}$,
   we select $\hat{D}_h$. By using the Euclidean orthogonal projection, we can push
    $\hat{J}\colon \hat{D}\rightarrow \hat{D}$ into an almost complex structure $J'\colon  \hat{D}_h\rightarrow \hat{D}_h$.

Since $|\mathrm{d}^j(\hat{D}-\hat{D}_h)|_{g_0}\leq O(k^{-1/2})$, for all $j\geq 0$, then $f$ is A.H. with respect to
 $(\hat{D},\hat{J},\hat{g}_k)$ if and only if it is A.H. with respect to $(\hat{D}_h,J',g_0)$ (this appears
  in the proof of lemma \ref{lem:localcheck2}).

In $\underline{\mathbb{C}}^m = (\mathbb{C}^n\times\mathbb{R})\times \mathbb{C}^m $ we have canonical
 coordinates $z_k^1,\dots,z_k^n,s_k,u_k^1,\dots,u_k^m$. These are CR coordinates w.r.t $(\hat{D}_h,\hat{J}_0)$.
  By hypothesis
\[\frac{\partial f}{{\partial}\bar{z}^1_k}=\cdots \frac{\partial f}{{\partial}\bar{z}^n_k}=\frac{\partial f}{{\partial}\bar{u}^1_k}=\cdots\frac{\partial f}{{\partial}\bar{u}^m_k}=0.\]
If we show that  $z_k^1,\dots,z_k^n,s_k,u_k^1,\dots,u_k^m$ are A.H. coordinates for $(\hat{D}_h,J',g_0)$
then we are done (this is again lemma \ref{lem:localcheck2} in the absence of connection form). But this
 follows from the fact that the trivialization of $\underline{\mathbb{C}}^m$ is given by an A.H. frame
  and therefore the induced distribution (by the connection form) $\mathcal{H}$ on $\hat{D}_h$ is  such
  that $|\mathrm{d}^j(\mathcal{H}-\hat{J}_0\mathcal{H})|_{g_0}\leq O(k^{-1/2})$, for all $j\geq 0$.

To prove the result in the other direction we fix A.H. coordinates and A.H. frames for $E_k$. The
 $\mathsf{T}\times Gl(m,\mathbb{C})$-invariance of  $(S^a)_{a\in A}\subset\underline{\mathbb{C}}^m$
  means that the local identifications define a sequence of global stratifications, and that these
  do not depend either on the A.H. coordinates or on the A.H. trivializations. It is an approximately
   holomorphic sequence of finite, Whitney (A) stratifications by the first  part of the proof.
\end{proof}

 In contrast to what happens for $0$-jets, it is not easy to find non-trivial approximately holomorphic
  stratifications for higher order jets.  The difficulty comes from the fact that the modification of
  the connection of proposition \ref{pro:perthol} that makes the $r$-jets of A.H. sequences of sections
  of $E_k$  into  A.H. sequences of sections of $\mathcal{J}_{D}^rE_k$, makes it very complicated to
   guarantee that the strata are given by  functions whose composition with an A.H. section is an A.H. function.

\begin{example}\label{ex:example1} Let $L^{\otimes k}_\Omega$ be the sequence of powers of the pre-quantum
line bundle of a  symplectic manifold of dimension $2p$. Let us consider the following sequence of strata in
$\mathcal{J}^1L^{\otimes k}_\Omega$:

\[\Sigma_{k,p}=\{(\sigma_0,\sigma_1)|\sigma_1=0\}.\]
The second subindex in our notation indicates the complex dimension of the kernel of the degree one homogeneous component of the 1-jet (see equation (\ref{eq:order1sing})).
Using the local sections $\mu_{k,x,I}$ of equation (\ref{eqn:canframe}), where $I=1,\dots,p$, and taking reference
 sections in Darboux charts, we get coordinates $z_k^1,\dots,z_k^p,v_k^0,v_k^1,\dots,v_k^p$ for the total
  space.  $\Sigma_{k,p}$ is then defined by the zeros of the function
  $f=(v_k^1,\dots,v_k^p)\colon \mathbb{C}^{2p+1}\rightarrow \mathbb{C}^p$, which is not  holomorphic
    (or A.H.) with respect to the modified almost complex structure of the total space. Otherwise, the
     composition  $f\circ j^1(z_k^1\tau_{k,x}^{\mathrm{ref}})$ would be  A.H., but that composition
      is $(1+z_k^1\bar{z}_k^1,z_k^1\bar{z}_k^2,\dots,z_k^1\bar{z}_k^p)$.

Actually, we cannot find A.H. functions $f$ defining  $\Sigma_{k,p}$: let us work in Darboux
 coordinates with the canonical complex structure  $J_0$ in the base.  Assume that $\mu_{k,x,I}$
  is built out of the reference section $\mathrm{e}^{-|z_k|^2/4}\xi$, where  $\xi$ is a unitary
  trivialization of $L_\Omega$ whose connection form is $A$ in equation (\ref{eqn:connform}).
  Then $\mathcal{J}^1L^{\otimes k}_\Omega$ becomes locally $\underline{\mathbb{C}}^{p+1}$
  with diagonal connection matrix $A\mathrm{I}_{p+1\times p+1}$. Proposition \ref{pro:perthol}
   for complex manifolds implies that the  modified   almost complex structure on
   $\underline{\mathbb{C}}^{p+1}$ is integrable. The submanifold $z_k^2=\cdots =z_k^p=v_k^2=\cdots =v_k^p=0$
    is complex with respect to the modified almost complex structure. Therefore, we can restrict our
     attention to the case $p=1$. The sections
     $j^1_{\mathrm{hol}}\mathrm{e}^{-|z_k|^2/4}\xi$, $j^1_{\mathrm{hol}}z_k\mathrm{e}^{-|z_k|^2/4}\xi$
     are by  proposition  \ref{pro:perthol} holomorphic. If we use them to trivialize
      $\mathcal{J}^1L^{\otimes k}_\Omega$ in a neighborhood of the origin, then we obtain
      a new identification with $\underline{\mathbb{C}}^{3}$ with its canonical complex structure.
       Let  $z_k$,$t_k$,$s_k$ be the new complex coordinates. A short computation shows that

\begin{eqnarray*}
 v_k^0&= & t_k+z_ks_k,\\
 v_k^1&= & -\bar{z}_k/2t_k+( 1-z_k\bar{z}_k/2)s_k.
\end{eqnarray*}
Hence  away from the origin $\Sigma_{k,p}$ admits the parametrization

\[(z_k,s_k)\mapsto (z_k,s_k,s_k(2/\bar{z}_k-z_k)).\]

Therefore, $\Sigma_{k,p}$ is not holomorphic with respect to the modified almost complex structure,
 and it follows that we cannot find $f$ A.H. defining $\Sigma_{k,p}$ locally.

\end{example}

\subsection{Quasi-stratifications}

For the applications we have in mind the notion of stratification has to be weakened. We start
doing it for the local model (endowed with the trivial connection).

Let $\sigma\in S$, $S$ a submanifold of $\mathcal{J}^{r+1}_{D_h,n,m}$. We say that
 $\alpha\in \Gamma(\mathcal{J}^r_{D_h,n,m})$ is a local representation for
 $\sigma$ if (i) $\alpha(0)=\pi^{r+1}_r\sigma$, and (ii) $\sigma=j^1_{D_h}\alpha(0)\in \mathcal{J}^{r+1}_{D_h,n,m}$,
  where $\pi^{r+1}_r\colon \mathcal{J}^{r+1}_{D_h,n,m}\rightarrow \mathcal{J}^r_{D_h,n,m}$
  is the natural projection, and  $j^1_{D_h}\alpha$  denotes the CR 1-jet of $\alpha$. The equality in (ii)
   should be understood in the following sense: the degree 1 component of the 1-jet should
   give an element of $\mathcal{J}^{r+1}_{D_h,n,m}$ (with vanishing degree 0 homogeneous component)
    and whose homogeneous components of degree $1,\dots,r+1$ coincide with those of $\sigma$.

\begin{definition}\label{def:quasistratum}(see \cite{Au02}) Let $S$  be a submanifold of $\mathcal{J}^r_{D_h,n,m}$
 (resp. $\mathcal{J}^r_{\mathbb{C}^g,p,m}$).  We define $\Theta_S$ to be the set of points  $\sigma \in S$
   for which  there exists an ($r+1$)-jet $\tilde{\sigma}$ (resp. ($r+1$)-jet along $G$) such that
   $\pi^{r+1}_r\tilde{\sigma}=\sigma$ and with a local representation $\alpha$ intersecting  $S$
    at $\sigma$ transversely along  $D_h$ (resp. along $\mathbb{C}^g$). We refer to $\Theta_S$
    as the holonomic transverse subset of $S$.
\end{definition}

It  can   be checked that if $S$ is invariant under the action of $\mathsf{T}\times (Gl(n,\mathbb{C})\times Gl(m,\mathbb{C}))$
 -the second factor $Gl(n,\mathbb{C})\times Gl(m,\mathbb{C})$ acting  fiberwise- (resp.
 $\mathsf{T}\times (Gl(g,\mathbb{C})\times Gl(m,\mathbb{C}))$), then $\Theta_S$  has the same invariance property.

When  an ($r+1$)-jet $\sigma$ is represented by  a local section of $\mathcal{J}^r_{D_h,n,m}$, in order
 to check whether $\pi^{r+1}_r\sigma \in S$  belongs to $\Theta_S$ the local representation is essentially
  unique: regarding transversality, it is enough to consider the degree $1$ part of the Taylor expansion
   in the coordinates $z_k^1,\bar{z}_k^1,\dots,z_k^n,\bar{z}_k^n$ (we turn the section into a function
    using the basis $\mu_I$). The degree $0$ part is determined by the  $r$-jet,  the hypothesis implies that
     the antiholomorphic part is vanishing and the holomorphic part is determined by the ($r+1$)-jet.
     That means in particular that we can restrict our attention to CR representations if necessary.

The importance of $\Theta_S$ is two-fold: on the one hand it will be used to define the stratifications
 we are interested in. On the other hand it is a very relevant subset when we study transversality to the
  strata: indeed, if $\tau$ is a CR section of $\underline{\mathbb{C}}^m$ and $\alpha:=j^r_{D_h}\tau$ is
  such that $\alpha(0)=\sigma$ and $\sigma\notin \Theta_S$, then $\alpha$ cannot be transverse along $D_h$
   to $S$ at $\sigma$ (notice that
    $\tilde{\sigma}:=(\tau(0),\textrm{d}_{D_h}\alpha(0))=j^{r+1}_{D_h}\tau(0)\in \mathcal{J}^{r+1}_{D_h,n,m}$,
    and therefore $\alpha$ is a local representation of $\tilde{\sigma}$). The consequence is that if
 $S\backslash\Theta_S\subset\partial S'$, transversality of $\tau$ to $S$ implies that $\tau$ misses  a
  neighborhood of $S\backslash\Theta_S$ in $S'$.

Definition \ref{def:quasistratum} extends to strata $S_k\subset \mathcal{J}^r_DE_k$ (resp. $ \mathcal{J}^r_GE_k$):
  we have a notion of pseudo-holomorphic $1$-jet of a section of $\mathcal{J}^r_DE_k$ (resp. pseudo-holomorphic
   $1$-jet along $G$ of a section of $\mathcal{J}^r_GE_k$) -because we have a  connection $\nabla_{H,D}$
   (resp. a connection on $\mathcal{J}^r_GE_k$ defined out of $\nabla_H$ and the projection
   $p_G\colon \mathcal{J}^rE_k\rightarrow \mathcal{J}^r_GE_k$)- and hence the notion of local
   representation. Then  $\Theta_{S_k}$ are those points $\sigma$ with lifts $\tilde{\sigma}$
   having a local representation transverse along $D$ (resp. $G$) to $S_k$ at $\sigma$.

Recall that once a family of A.H. charts has been fixed we have identifications
$\varpi_{k,x}\colon T^{*1,0}\mathbb{C}^n\rightarrow D^{*1,0}$. If we also fix a family of A.H.
 trivializations of $E_k$ over the charts there is an induced identification

\begin{equation}\label{eq:jetidentity}
\Pi_{k,x}\colon \mathcal{J}^r_DE_k\rightarrow \mathcal{J}^r_{D_h,n,m}.
\end{equation}

\begin{lemma}\label{lem:locrep} Let $S_k$ be a sequence of strata of  $\mathcal{J}^r_DE_k$, where either
 $r=0$ and $E_k=E\otimes L^{\otimes k}$, or $E_k=\underline{\mathbb{C}}^m$ and $r\in \mathbb{N}$.

\begin{enumerate}
 \item If $E_k=\underline{\mathbb{C}}^m$ assume that  for a choice of A.H. charts  $\Pi_{k,x}(S_k)=S$,
  where $S\subset\mathcal{J}^r_{D_h,n,m}$ is invariant under the action $\mathsf{T}\times Gl(n,\mathbb{C})$.
   Then $\Pi_{k,x}(\Theta_{S_k})=\Theta_S$.
\item The same result holds for $E\otimes L^{\otimes k}$ and $r=0$; we need to fix A.H. trivializations of
 $E_k$ (so  $\Pi_{k,x}$ is defined) and require invariance of $S\subset\underline{\mathbb{C}}^m$ under
  the action of  $\mathsf{T}\times Gl(m,\mathbb{C})$.
\end{enumerate}

For jets along $G$ we have analogous results, but we need A.H. charts adapted to $G$  and  we ask for
$\mathsf{T}\times Gl(g,\mathbb{C})$-invariance of $S$ instead of $\mathsf{T}\times Gl(n,\mathbb{C})$-invariance.
\end{lemma}
\begin{proof}
Since $S$ is $Gl(n,\mathbb{C})$-invariant, so is  $\Theta_S$.
We have the local identifications $\Pi_{k,x}\colon \mathcal{J}^r_DE_k\rightarrow \mathcal{J}^r_{D_h,n,m}$.
Let $y\in M$ belong to $B(0,\rho)$ in the domain of the charts centred at $x_1$ and $x_2$, for some $k$.
 Then there is a fiber bundle isomorphism
\begin{equation}\label{eq:jetid2}
\Phi_{k,x_1,x_2}\colon \mathcal{J}^r_{D_h,n,m}\rightarrow \mathcal{J}^r_{D_h,n,m}
\end{equation}
defined as follows: for each point $y$ in the intersection of the domains of the charts, the restriction
 of the differential to $D$ is a complex $J$-linear map $L_y$. Consider the  linear map
 $\varpi_{k,x_2}\circ L^*_y\circ \varpi^{-1}_{k,x_2}\colon T^{*1,0}\mathbb{C}^n\rightarrow T^{*1,0}\mathbb{C}^n$,
  which belongs to $Gl(n,\mathbb{C})$. $\Phi_{k,x_1,x_2}$ in the fiber over $y$ (or over the origin in
   both charts due to the $\mathsf{T}$-invariance) is the vector space isomorphism induced by
   $\varpi_{k,x_2}\circ L^*_y\circ \varpi^{-1}_{k,x_2}$ (and the identity acting on the $\mathbb{C}^m$
    factor of the tensor product). Since $S$ is invariant under the $\mathsf{T}\times Gl(n,\mathbb{C})$-action, $\Phi_{k,x_1,x_2}(\Theta_S,S)=(\Theta_S,S)$.
 In particular  the pair  $(\Theta_S,S)$ does not depend on the chosen family of A.H. charts.  We
  construct an appropriate family of A.H. charts (there is no Darboux condition involved here) by
  the usual rescaling procedure, but starting from normal coordinates composed with a linear transformation
   so that $(D,J)=(D_h,J_0)$ at the origin. Recall that since $E_k=\underline{\mathbb{C}}^m$, the connection
     $\nabla_{k,r}$  on $\mathcal{J}^r_DE_k$  is just induced by the Levi-Civita connection (in the
     $\underline{\mathbb{C}}^m$ factor we use the trivial connection $\mathrm{d}$). Hence the pushforward
      of $\nabla_{k,r}$ by $\Pi_{k,x}$ to $\mathcal{J}^r_{D_h,n,m}$ has vanishing connection form at the
      origin. Since we also have  $(D\oplus D^\perp,J)=(D_h\oplus D_v,J_0)$ at the origin, for any section
       $\alpha$ of $\mathcal{J}^r_{D_h,n,m}$ we have $j^1_D\alpha(0)=j^1_{D_h}\alpha(0)$. Therefore, the
       local representations at the origin for the canonical CR structure and the induced one coincide.
       From that and $D=D_h$ at the origin, we conclude $\Pi_{k,x}(\Theta_{S_k})=\Theta_S$.

Item (2) is proven in the same fashion. The  $Gl(m,\mathbb{C})$-invariance implies that we can choose
any arbitrary family of A.H. trivializations. What we do is selecting trivializations such that the
connection form over the origin is vanishing (here we deal with the connection $\nabla_k$ on $E_k$).

Notice that we cannot state item (2) for higher order jets because  the action of
$Gl(n,\mathbb{C})\times Gl(m,\mathbb{C})$ does not allow us to kill at the origin of each chart the connection
 form of  the modified connection $\nabla_{k,H_r}$.

For the relative results we start by modifying a bit  the vector bundle isomorphism
$\varpi_{k,x}\colon T^{*1,0}\mathbb{C}^p\rightarrow T^{*1,0}P$; the original $(J_0,J)$-complex map
 $T\mathbb{C}^p\rightarrow TP$ can  be easily arranged to be compatible with the splittings
  $T\mathbb{C}^g\oplus T\mathbb{C}^{p-g}$ and $G\oplus G^\perp$. Due to the
  $\mathsf{T}\times (Gl(g,\mathbb{C}))$-invariance we are free to pick any family of A.H. charts
   adapted to $G$. The ones we need come from rescaling  normal coordinates composed with  a linear
    transformation sending $(G\oplus G^{\perp},J)$ to $(\mathbb{C}^g\oplus \mathbb{C}^{p-g},J_0)$
     at the origin. In these coordinates the connection form on $T^{*1,0}\mathbb{C}^g$ is vanishing,
      because we project the Levi-Civita connection which is already vanishing at the origin. Hence
       the 1-jets along $G$ and $\mathbb{C}^g$ at the origin coincide (also because
         $(G\oplus G^{\perp},J)=(\mathbb{C}^g\oplus \mathbb{C}^{p-g},J_0)$ at the origin), and this proves the result.
\end{proof}

The only relevant strata $S_k\subset \mathcal{J}^r_{D}E_k$ for which we have to consider the subsets
 $\Theta_{S_k}$ are the zero sections $Z_k$.  In that case (see \cite{Au00}) the subsets $\Theta_{Z_k}$
  are those $r$-jets whose degree 1 component is onto.

\begin{definition}\label{def:quasidef} An approximately holomorphic quasi-stratification of $\mathcal{J}^r_{D}E_k$
 is an approximately holomorphic stratification in which the partial order condition is relaxed in the following way: $Z_k$ are strata of the quasi-stratification, and for any other  strata $S_k\neq Z_k$ when we approach $Z_k$, it
 accumulates into points of $Z_k\backslash\Theta_{Z_k}$ (so in particular $Z_k$ is not in the closure of $S_k$).
\end{definition}

\subsection{The Thom-Boardman-Auroux stratification for maps to projective spaces}\label{ssec:TBA}

Let $E_k=\underline{\mathbb{C}}^{m+1}\otimes L^{\otimes k}$. Let $Z^0,\dots,Z^m$ be the complex coordinates
 associated to the trivialization of $\mathbb{C}^{m+1}$ (at any fiber) and let
  $\pi\colon \mathbb{C}^{m+1}\backslash \{0\}\rightarrow \mathbb{C}\mathbb{P}^m$ be the
   canonical projection. Consider the canonical affine coordinates
\begin{eqnarray}
 \nonumber\varphi_i^{-1}\colon U_i& \longrightarrow  & \mathbb{C}^m \\ \nonumber
 [Z_0:\dots:Z_m]& \longmapsto & \left(\frac{Z^0}{Z^i},\dots,\frac{Z^{i-1}}{Z^i},\frac{Z^{i+1}}{Z^i},\dots,\frac{Z^m}{Z^i}\right).
 \end{eqnarray}
For each chart $\varphi_i$ we consider the bundle
\begin{equation}\label{eqn:projjets}
\mathcal{J}_{D}^r(M,\mathbb{C}^m)_i:=(\sum_{j=0}^r {(D^{*1,0})}^{\odot j})\otimes \underline{\mathbb{C}}^m.
\end{equation}

We now bring back the discussion at the beginning of section \ref{sec:pholjets}. Assume for the
moment that $M$ is a Levi-flat CR manifold and fix a family of CR charts. Over each of the balls
$B_{g_k}(x,\rho)$ we have the bundles $\mathcal{J}^r_{D_h,n,m}$ of CR $r$-jets. Notice that if we
use the frames  $\mu_{k,x,I}$ of equation (\ref{eqn:canframe}) they are vector bundles.

The local bundles $\mathcal{J}^r_{D_h,n,m}$  glue into the non-linear bundle  $\mathcal{J}^r_{CR}(M,\mathbb{C}^m)_i$:
 let $y\in M$ be a point belonging to two different charts centred at $x_0$ and $x_1$ respectively.
 If we send  $y$ in both charts to the origin via a translation, then the change of coordinates
 restricts to the leaf through the origin to a holomorphic map fixing the origin. The fibers over
  $y$ are  related  by the action of the holomorphic $r$-jet of the bi-holomorphism.  If we only take
   the linear part of the action -which is the vector bundle map $\Phi_{k,x_1,x_2}$ of equation
   (\ref{eq:jetid2})- we are equally defining a bundle, for the cocycle condition still holds. Moreover,
   it is a vector bundle. Besides, since we only use the linear part we do not need either $D$ or $J$
    to be integrable. This bundle is  $\mathcal{J}_{D}^r(M,\mathbb{C}^m)_i$ as defined in   equation
     (\ref{eqn:projjets}) (what we defined there, it  is rather a sequence in which the metric in the
      $D^{*1,0}$ factors is induced from $g_k$). Thus for  Levi-flat manifolds the vector bundles
       $\mathcal{J}_{D}^r(M,\mathbb{C}^m)_i$ are ``linear approximations'' of the non-linear bundles
        $\mathcal{J}^r_{CR}(M,\mathbb{C}^m)_i$.

\begin{proposition}\label{pro:projjets}\quad
\begin{enumerate}
\item The  vector bundles $\mathcal{J}_{D}^r(M,\mathbb{C}^m)_i$ can be glued to define the almost complex
 fiber bundles $\mathcal{J}^r_D(M,\mathbb{C}\mathbb{P}^m)$  of pseudo-holomorphic $r$-jets of maps from
 $M$ to $\mathbb{C}\mathbb{P}^m$, so that their fibers inherit a canonical holomorphic structure.

\item Given $\phi_k\colon M\rightarrow \mathbb{C}\mathbb{P}^m$ there is a notion of pseudo-holomorphic
 $r$-jet extension $j^r_D\phi_k\colon M\rightarrow \mathcal{J}^r_D(M,\mathbb{C}\mathbb{P}^m)$, which is
  compatible with the notion of pseudo-holomorphic $r$-jet for the sections
  $\varphi_i^{-1}\circ\phi_k\colon M\rightarrow \mathbb{C}^m$ of definition \ref{def:psjets}.
   If $\phi_k\colon M\rightarrow \mathbb{C}\mathbb{P}^m$ is an A.H. sequence then
    $j^r_D\phi_k\colon M\rightarrow  \mathcal{J}^r_D(M,\mathbb{C}\mathbb{P}^m)$ is also A.H.
\end{enumerate}

Analogous results hold in the relative setting for the bundles $\mathcal{J}^r(P,\mathbb{C}\mathbb{P}^m)$
and $\mathcal{J}^r_G(P,\mathbb{C}\mathbb{P}^m)$. Also there is an approximately holomorphic sequence of
 canonical  submersions $p_G\colon \mathcal{J}^r(P,\mathbb{C}\mathbb{P}^m) \rightarrow\mathcal{J}^r_G(P,\mathbb{C}\mathbb{P}^m)$.
  These submersions are left inverses of the natural inclusions
   $l_G\colon \mathcal{J}^r_G(P,\mathbb{C}\mathbb{P}^m)\hookrightarrow \mathcal{J}^r(P,\mathbb{C}\mathbb{P}^m)$
   so that for $\phi_k\colon P\rightarrow \mathbb{C}\mathbb{P}^m$ an A.H. sequence,
    $j^r_G\phi_k\colon P\rightarrow \mathcal{J}^r_G(P,\mathbb{C}\mathbb{P}^m) \hookrightarrow \mathcal{J}^r(P,\mathbb{C}\mathbb{P}^m)$
     is A.H.
\end{proposition}
\begin{proof}
Let us denote  the change of coordinates $\varphi_j^{-1}\circ \varphi_i$  by $\Psi_{ji}$.
For any   $y\in M$ the points in  $\{y\}\times (U_i\cap U_j)\subset\mathcal{J}_{D}^r(M,\mathbb{C}^m)_i$
are identified with points in  $\{y\}\times (U_i\cap U_j)\subset\mathcal{J}_{D}^r(M,\mathbb{C}^m)_j$
using the same transformation $j^r\Psi_{ji}$ in $\mathcal{J}^r_{n,m}$ induced by the fiberwise holomorphic change
 of coordinates $\Psi_{ji}$. In other words, if we take an approximately holomorphic chart centred at
  $x$ say and containing  $y$, we get as in equation (\ref{eq:jetidentity}) a vector bundle isomorphism
  $\Pi_{k,x,i}\colon\mathcal{J}_{D}^r(M,\mathbb{C}^m)_i\rightarrow \mathcal{J}^r_{D_h,n,m}$. Thus
  for $\sigma\in \mathcal{J}_{D}^r(M,\mathbb{C}^m)_i$ there exists $F\colon \mathbb{C}^n\rightarrow \mathbb{C}^m$
   a CR function such that $\Pi_{k,x,i}(\sigma)=j^r_{D_h}F(x)$.

The bundle map we define is:
\begin{align}\nonumber
j^r\Psi_{ji}\colon \mathcal{J}_{D}^r(M,\mathbb{C}^m)_i&\longrightarrow \mathcal{J}_{D}^r(M,\mathbb{C}^m)_j\\\label{eq:jetid}
\sigma &\longmapsto  \Pi_{k,x,j}^{-1}(j^r_{D_h}(\Psi_{ji}\circ F)(x)).
\end{align}
 This map does not depend either on the charts: if we have a point $y$ in two different  charts
 centred at $x_1$ and $x_2$, then we saw in the proof of lemma \ref{lem:locrep} that the vector space
   isomorphism $\Phi_{k,x_1,x_2}\colon \mathcal{J}^r_{D_h,n,m}\rightarrow \mathcal{J}^r_{D_h,n,m}$
   was induced by $T\in Gl(n,\mathbb{C})$. The bundle map of equation (\ref{eq:jetid}) is equivariant with respect to
    this action, because in the CR setting it is equivariant with respect to the action in the base of CR transformations.
     Hence, the result follows by considering the affine CR transformation sending $y$ in the first chart
     to its image in the second and whose linear part is
     $T^*\times \mathrm{I}\colon \mathbb{C}^{n}\times\mathbb{R}\rightarrow \mathbb{C}^{n}\times\mathbb{R}$.

Equivalently,  the $r$-jet of $\Psi_{ji}\circ F$ admits a coordinate free expression only in terms of the $r$-jet of $F$.

Therefore the identifications $j^r\Psi_{ji}$  give rise to  a well defined locally trivial fiber bundle
  $\mathcal{J}^r_D(M,\mathbb{C}\mathbb{P}^m)$.

\begin{remark}\label{rem:linid} If our manifold is CR and we have $x$ belonging to two different CR charts,
 then there is a natural induced identification $\mathcal{J}^r_{D_h,n,m}\rightarrow \mathcal{J}^r_{D_h,n,m}$
  over the points belonging to both charts. This identification is just the action of the CR $r$-jet of the
   change of coordinates. We observe that this is not the action of $\Phi_{k,x_1,x_2}$, which is just the
    action induced by the 1-jet of the change of coordinates (the only one available for all almost CR structures!).
\end{remark}

The fibers of $\mathcal{J}^r_D(M,\mathbb{C}\mathbb{P}^m)$  admit a canonical holomorphic structure because
 using the local identifications $\Pi_{k,x,i}$ the fiber is some $\mathbb{C}^N$ and the change of coordinates
  is a  fiberwise holomorphic map (because it is the holomorphic $r$-jet of $\Psi_{ji}$), and this proves item (1).

Let $\phi\colon (M,J,D)\rightarrow \mathbb{C}\mathbb{P}^m$. Its pseudo-holomorphic $r$-jet $j^r_{D}\phi$
is defined as follows: the affine charts of  projective space induce maps  $\phi_i:=\varphi_i^{-1}\circ \phi
\colon M\rightarrow \mathbb{C}^m$. Using the trivial connection $\mathrm{d}$ in this trivial vector bundle
 and  the induced connection on $D^{*1,0}$, we can define the corresponding pseudo-holomorphic $r$-jet
 $j^r_{D}\phi_i$ (definition \ref{def:psjets}). We must check that
 \begin{equation}\label{eq:welldefjet}
 j^r_{D}\phi_j=j^r\Psi_{ji}(j^r_{D}\phi_i).
 \end{equation}

More generally let  $H\colon \mathbb{C}^{m_1}\rightarrow
\mathbb{C}^{m_2}$ be any holomorphic map. Then use the local identifications
$\Pi_{k,x,s}\colon\mathcal{J}_{D}^r(M,\mathbb{C}^{m_s})\rightarrow \mathcal{J}^r_{D_h,n,m_s}$, $s=1,2$,
to induce the map $j^rH\colon \mathcal{J}_{D}^r(M,\mathbb{C}^{m_1})\rightarrow
\mathcal{J}_{D}^r(M,\mathbb{C}^{m_2})$. We claim that for any function  $\phi\colon M\rightarrow \mathbb{C}^{m_1}$
 we have
\begin{equation}\label{eq:welldefjet2}
j_{D}^r(H\circ \phi)=j^rH(j^r_{D}\phi).
\end{equation}
Equation (\ref{eq:welldefjet}) follows from the claim by taking $H=\Psi_{ji}$.

The proof of the claim take the next two and a half pages, and it is by induction on $r$. Firstly we notice that from the proof of the claim for $m_2=1$, the proof for any $m_2$ follows immediately. Therefore we assume $m_2=1$. Secondly we observe that it is enough to check the equality in (\ref{eq:welldefjet2}) for the
 degree $r$ homogeneous component of the $r$-jet.

 We shall denote  the degree $r$ homogeneous component of $j^rH$
 by $\textrm{d}^rH$; recall that $\textrm{d}^rH(j^r_{D}\phi(x))$ depends on the components of every order
 of $j^r_{D}\phi(x)$. Let $F=(F^1,\dots,F^{m_1})\colon \mathbb{C}^n\times \mathbb{R}\rightarrow \mathbb{C}^{m_1}$
  be a CR function such that

\[j^r_D\phi(x)=j^r_{D_h}F(x).\]
Also the degree $r$ homogeneous component of $j^r_{D_h}F$ is denoted by $\partial^r_0F$.
 By definition

\begin{equation}\label{eq:jetequal}
\partial^j_{\mathrm{sym}}\phi(x)=\partial^j_0F(x),\; j=0,\dots,r.
\end{equation}

We start the proof of the claim for 1-jets. Once we use the identification $\partial \phi(x)=\partial_0F(x)$, we have

\begin{equation}\label{eq:1jetbis}
\textrm{d} H(\partial\phi(x)):=\textrm{d} H(\partial_0F(x))=\sum_{a=1}^{m_1}\frac{\partial_0 H}{\partial_0 z^a}\partial_0 F^a(x),
\end{equation}
and using  the identification of equation (\ref{eq:jetequal}) back we get the following formula  for the right hand side of equation (\ref{eq:welldefjet2}) for 1-jets:

\begin{equation}\label{eq:1jetbisbis}
\textrm{d} H(\partial\phi(x))=\textrm{d} H(\partial_0F(x))=\sum_{a=1}^{m_1}\frac{\partial_0 H}{\partial_0 z^a}\partial \phi^a(x),
\end{equation}
where the partial derivatives of $H$ are evaluated on $\phi(x)=F(x)$, but we omit it in the notation.

Regarding the left hand side  of equation (\ref{eq:welldefjet2}),
the computation of  $\partial (H\circ \phi)(x)$ is done by firstly taking in  $\nabla (H\circ \phi)(x)$ its
 projection over $D^*$ (or restricting the differential to $D$).  Since

\begin{equation}\label{eq:1jet}
\nabla (H\circ \phi)(x)=\sum_{a=1}^{m_1}\frac{\partial_0 H}{\partial_0 z^a}\nabla\phi^a(x)
\end{equation}
is the sum of partial derivatives of $H$ multiplied by the components $\nabla \phi_a(x)$ of $\nabla \phi(x)$,
 taking $\nabla_D (H\circ \phi)(x)$ amounts to substituting  in equation (\ref{eq:1jet}) the factors  $\nabla \phi^a(x)$
 by $\nabla_{D} \phi^a(x)$.

Next the holomorphic component is singled out; since  $H$ is holomorphic $\partial (H\circ \phi)(x)$ is
computed by taking the component   $\partial \phi^a(x)$ of $\nabla_{D}\phi^a(x)$ in equation (\ref{eq:1jet}).
Thus we obtain the same result as in equation (\ref{eq:1jetbisbis}), and this proves the claim for 1-jets.

We need to prove the claim for  2-jets before going to the induction step. The reason is that for 1-jets the symmetrization step is not present, unlike the case of higher order jets.

By definition
\begin{equation}\label{eq:2jetbisbis}
\textrm{d}^2 H(j^2_{D_h}F(x))=\sum_{b,a=1}^{m_1}\frac{\partial^2_0H}{\partial_0 z^a\partial_0 z^b}\partial_0 F^a(x)\otimes\partial_0 F^b(x)+
\sum_{c=1}^{m_1}\frac{\partial_0 H}{\partial_0 z^c}\partial^2_0 F^c(x),
\end{equation}
so using equation (\ref{eq:jetequal}) we get for the right hand side of equation (\ref{eq:welldefjet2})

\begin{equation}\label{eq:2jetbisbisbis}
\textrm{d}^2 H(j^2_{D}\phi(x))=\sum_{b,a=1}^{m_1}\frac{\partial^2_0H}{\partial_0 z^a\partial_0 z^b}\partial \phi^a(x)\otimes\partial \phi^b(x)+
\sum_{c=1}^{m_1}\frac{\partial_0 H}{\partial_0 z^c}\partial^2_{\mathrm{sym}} \phi^c(x).
\end{equation}
To compute $\partial^2_\mathrm{sym}(H\circ \phi)(x)$  we first differentiate $\partial(H\circ \phi)$ at $x$:

\begin{equation}\label{eq:2jet}
\nabla \partial(H\circ \phi)(x)=\sum_{b,a=1}^{m_1}\frac{\partial^2_0H}{\partial_0 z^a\partial_0 z^b}\nabla \phi^a(x)\otimes\partial \phi^b(x)+
\sum_{c=1}^{m_1}\frac{\partial_0 H}{\partial_0 z^c}\nabla\partial \phi^c(x).
\end{equation}
Taking the component along $D$ and then the holomorphic part amounts to substituting in equation
(\ref{eq:2jet}) $\nabla \phi^a(x)$ by $\partial \phi^a(x)$, and  $\nabla\partial\phi^c(x)$ by $\partial^2 \phi^c(x)$:

\begin{equation}\label{eq:2jetbis}
\partial^2(H\circ \phi)(x)=\sum_{b,a=1}^{m_1}\frac{\partial^2_0H}{\partial_0 z^a\partial_0 z^b}\partial \phi^a(x)\otimes\partial \phi^b(x)+
\sum_{c=1}^{m_1}\frac{\partial_0 H}{\partial_0 z^c}\partial^2 \phi^c(x).
\end{equation}
We need to show that symmetrizing equation (\ref{eq:2jetbis}) amounts to writing $\partial^2_{\mathrm{sym}}\phi^c(x)$
 instead of $\partial^2 \phi^c(x)$.

In equation (\ref{eq:2jetbis}) we have  terms of ``type'' 2 -those containing a second derivative of $\phi$-
and  terms of ``type'' (1,1) which contain the tensor product of two derivatives of $\phi$. Terms of ``type''
 (1,1) are already  symmetric (just exchange the indices $a,b$); the symmetrization -being a linear projection-
  does not alter them. Now one checks that the symmetrization of each summand
   $\frac{\partial_0 H}{\partial_0 z^c}\partial^2 \phi^c(x)$ is exactly  $\frac{\partial_0 H}{\partial_0 z^c}\partial^2_{\mathrm{sym}}\phi^c(x)$, this proving the claim for 2-jets.

We now move onto the induction step. We assume   $\mathrm{d}^{r}H(j^{r}_D\phi(x))=\partial^{r}_\mathrm{sym}(H\circ \phi)(x)$ and we want to prove the claim for $(r+1)$-jets.  By a partition of $r$ of degree $s$ we understand any (ordered) s-tuple $(r_1,\dots,r_s)$, $1\leq s\leq r$,
 $1\leq r_j\leq r$, $r_1+\cdots +r_s=r$. In the computation of $\mathrm{d}^{r}H(j^{r}_D\phi(x)):=\partial^{r}_0(H\circ F)(x)$
  we  get an algebraic expression whose summands are of the form

\begin{equation}\label{eq:induction}
\frac{\partial_0^{r_1+\cdots +r_s}H}{\partial_0^{r_1} z^{i_1}\dots \partial_0^{r_s} z^{i_s}}\partial_0^{r_1}F^{i_1}(x)\otimes\cdots\otimes \partial^{r}_0 F^{i_s}(x),
\end{equation}
each belonging to a partition $(r_1,\dots,r_s)$. Notice that to some partitions correspond summands
that are originated from different partitions of $r-1$.  For example, in degree 3 we have (1,2)-terms
coming from the derivation of the terms of ``type'' 2 and others obtained from the derivation of the
(1,1)-terms. We do not add summands of the same ``type'', but keep them distinguished. By induction we
 assume that $\partial^r_\mathrm{sym}(H\circ\phi)(x)$ is computed  by the same algebraic expression as
  $\mathrm{d}^{r}H(j^{r}_{D_h}F(x))$, but writing in the summands in equation (\ref{eq:induction})
   $\partial^{r_j}_\mathrm{sym}\phi^{i_j}$ in place of $\partial^{r_j}_0F^{i_j}(x)$, and then evaluating at $x$.

To compute $\partial^{r+1}_\mathrm{sym}(H\circ\phi)(x)$ we have to firstly differentiate the algebraic
expression that computes $\partial^r_\mathrm{sym}(H\circ\phi)(x)$. From the previous assumption  a one to
 one correspondence compatible with the partitions between the summands of $\mathrm{d}^{r+1}H(j^{r+1}_{D_h}F(x))$
  and of $\nabla\partial^{r}_\mathrm{sym}(H\circ\phi)(x)$  can be established. It is clear that restricting to
   $D$ and taking the (1,0)-component does not affect the identification.

In each summand of  $\partial \partial_\mathrm{sym}^r(H\circ \phi)(x)$ all the factors but possibly one
 in the tensor product are of the form  $\partial^{r_j}_{\mathrm{sym}}\phi^{i_j}$ and hence already symmetric;
  the different one is of the form $\partial\partial^{r_j'}_{\mathrm{sym}}\phi^{i_j'}$.
Observe  that the symmetrization of each summand in $\partial \partial_\mathrm{sym}^r(H\circ \phi)(x)$
amounts to putting  instead of $\partial\partial^{r_j'}_{\mathrm{sym}}\phi^{r_j'}$, its symmetrization
$\partial_{\mathrm{sym}}^{r_{j}'+1}\phi^{r_j'}$ and then symmetrizing the resulting expression (this is
an elementary result concerning symmetric products which is proved by suitably regrouping the permutations).
Thus we have proven that
 \[\partial^{r+1}_\mathrm{sym}(H\circ\phi)(x)=\textrm{sym}_{r+1}(\mathrm{d}^{r+1}H(j^{r+1}_{D}\phi(x))),\]
 but $\mathrm{d}^{r+1}H(j^{r+1}_{D_h}F(x))$ is already symmetric. Therefore we conclude
  \[\partial^{r+1}_\mathrm{sym}(H\circ\phi)(x)=\mathrm{d}^{r+1}H(j^{r+1}_{D}\phi(x)),\]
  where the equality also holds for each summand in the algebraic expression computing both quantities.

Therefore we conclude that the pseudo-holomorphic
 $r$-jet of a map to  $\mathbb{C}\mathbb{P}^m$ is well defined.

To be able to say when a sequence of functions of  $\mathcal{J}^r_D(M,\mathbb{C}\mathbb{P}^m)$ is A.H. we need
to introduce  an almost CR structure in the total space of the $r$-jets. This can be done using a connection
(for example out of the Levi-Civita connection associated to the Fubini-Study metric in the projective space
 and of the connection on  $D^*$). In our case we choose to do something different but equivalent: we use
  the identifications with $\mathcal{J}^r_D(M,\mathbb{C}^m)_i$.  Each of these trivial vector bundles with
  trivial connection has a natural almost CR structure.  Let $K_i\subset U_i$ be compact sets whose
  interiors cover $\mathbb{C}\mathbb{P}^m$. We have the corresponding subsets
  $ \mathcal{J}^r_D(M,\varphi_i^{-1}(K_i))\subset \mathcal{J}^r_D(M,\mathbb{C}^m)_i$.

We say that $\sigma_k\colon M\rightarrow \mathcal{J}^r_D(M,\mathbb{C}\mathbb{P}^m)$ is A.H. if there exist
 constants $(C_j)_{j\geq 0}$ such that for all $x\in M$, $j\geq 1$, and $k\in \mathbb{N}$

\begin{align*}
\mathrm{max}_{i\in\{0,\dots,m\}}|\nabla^j(j^r\varphi_i^{-1}\circ\sigma_k)(x)|_{g_k}&\leq  C_j,\\
\mathrm{max}_{i\in\{0,\dots,m\}}|\nabla^{j-1} \bar{\partial}(j^r\varphi_i^{-1}\circ\sigma_k)(x)|_{g_k} &\leq  C_j k^{-1/2},
\end{align*}
where for each $x$ we only take into account those indices for which $\sigma_k(x)$ belongs to the interior
of $ \mathcal{J}^r_D(M,K_i)$.

Notice that in the local models the identifications $j^r\Psi_{ji}$ are holomorphic, therefore when restricted
 to subsets associated to compact regions of $\mathbb{C}^m_i$ and $\mathbb{C}^m_j$ the sequence of maps
 $j^r\Psi_{ji}\colon \mathcal{J}_D^r(M,\mathbb{C}^m)_i\rightarrow  \mathcal{J}_D^r(M,\mathbb{C}^m)_j$ is
 A.H. In particular the notion of a sequence $\sigma_k\colon M\rightarrow \mathcal{J}^r_D(M,\mathbb{C}\mathbb{P}^m)$
  being A.H. does not depend on the covering $K_i$. It is also clear that if a sequence of functions $\phi_k$
  is A.H. then $j^r_D\phi_k$ is also A.H.  This proves item (2) of the proposition.

If $(P,\Omega)$ is symplectic the definition of  $\mathcal{J}^r(P,\mathbb{C}\mathbb{P}^m)$ is the same
(we just do not need to project the full derivative into the subspace $D^*$).
When we have a $J$-complex distribution  $G$ there is an analogous definition of the bundle of pseudo-holomorphic
 $r$-jets along $G$.  Using the previous affine coordinates of  projective space we consider the sub-bundles
\[\mathcal{J}_G^r(P,\mathbb{C}^m)_i:=(\sum_{j=0}^r {(G^{*1,0})}^{\odot j})\otimes \underline{\mathbb{C}}^m,\]
 where $\mathcal{J}_G^r(P,\mathbb{C}^m)_i\subset \mathcal{J}^r(P,\mathbb{C}^m)_i$ via the splitting $G\oplus G^{\perp}=TP$.

 It is easily checked using the local identification between $\mathcal{J}^r_{p,m}$ and $\mathcal{J}(P,\mathbb{C}^m)$
  coming from approximately holomorphic coordinates adapted to $G$, that the diffeomorphisms
  $j^r\Psi_{ji}\colon \mathcal{J}^r(P,\mathbb{C}^m)_i\rightarrow \mathcal{J}^r(P,\mathbb{C}^m)_j$ preserve these sub-bundles.

The  proof  that shows that the  $j^r\phi$ is well defined is exactly the same we gave for 2-calibrated
manifolds; a small modification shows that $j_G^r\phi$ is well defined (instead of keeping the component
 $\nabla_D$ of the odd dimensional case, we project over $G^*$).

Going to the models furnished by approximately holomorphic coordinates adapted to $G$, the submersion
$p_G\colon \mathcal{J}^r_{p,m}\rightarrow \mathcal{J}^r_{\mathbb{C}^g,p,m}$ is just a projection on
some of the holomorphic coordinates, and therefore it is an approximately holomorphic sequence of maps.

Using approximately holomorphic coordinates adapted to $G$ it is straightforward to check that if
$\phi_k\colon P\rightarrow \mathbb{C}\mathbb{P}^m$ is A.H., then both $j^r_G\phi_k$ and $j^r\phi_k$
are A.H. sequences of $\mathcal{J}^r(P,\mathbb{C}\mathbb{P}^m)$.
\end{proof}

 We recall that  $Z_k$ denotes the sequence of strata of  $\mathcal{J}^r_DE_k$ (resp.  $\mathcal{J}^rE_k$,
   $\mathcal{J}^r_GE_k$) of $r$-jets whose degree $0$-component vanishes. We define
    $\mathcal{J}^r_DE_k^*:=\mathcal{J}^r_DE_k\backslash Z_k$ (resp. $\mathcal{J}^rE_k^*:=\mathcal{J}^rE_k\backslash Z_k$,
     $\mathcal{J}^r_GE_k^*:=\mathcal{J}^r_GE_k\backslash Z_k$).

\begin{proposition}\label{pro:TBAsubmer}\quad
\begin{enumerate}
\item There exists a bundle map  $j^r\pi\colon \mathcal{J}^r_DE_k^*\rightarrow \mathcal{J}^r_D(M,\mathbb{C}\mathbb{P}^m)$
 which is a fiberwise holomorphic submersion.
\item Let $\tau_k$ be a section of $E_k$, and let $\phi_k=\pi\circ \tau_k\colon M\backslash Z(\tau_k)\rightarrow \mathbb{C}\mathbb{P}^m$ be its projectivization defined away from the zero subset of $\tau_k$. Then the following equation holds:
\begin{equation}\label{eq:fundeq}
j^r\pi(j^r_D\tau_k)=j^r_D\phi_k.
\end{equation}
\end{enumerate}

In the almost complex case we  have an analogous  map $j^r\pi$, and for  $\tau_k\colon P\rightarrow E_k$ and its projectivization $\phi_k$ the equality
\begin{equation}\label{eq:fundeq2}
j^r\pi(j^r\tau_k)=j^r\phi_k
\end{equation}
holds where defined.

Given  $G$ a $J$-complex distribution we have  the following commutative square of submersions:

\begin{equation}\label{cd:submersions}\begin{CD}\mathcal{J}^rE_k^* @>{p_G}>>  \mathcal{J}^r_GE_k^*\\
@VV{j^r\pi}V         @VV{j^r\pi}V     \\
\mathcal{J}^r(P,\mathbb{C}\mathbb{P}^m) @>{p_G}>>  \mathcal{J}^r_G(P,\mathbb{C}\mathbb{P}^m)
\end{CD}
\end{equation}

\vskip 3mm
If $j^r_G\tau_k$ is a section of $\mathcal{J}^r_GE_k^*$ the equality
\begin{equation}\label{eq:fundeq3}
 j^r\pi(j^r_G\tau_k)=j^r_G\phi_k
 \end{equation}
 holds where defined.
 \end{proposition}

\begin{proof}
We define  $j^r\pi$ to have the same expression as in the integrable case. That means that we fix approximately
holomorphic coordinates and a section trivializing $L^{\otimes k}$ and a local frame of $E=\underline{\mathbb{C}}^{m+1}$,
 so that the  $r$-jet $\sigma$ in question is identified with the usual CR  $r$-jet at a point $x$ of a CR function  $F$.
Then  $j^r\pi(\sigma)$ is defined  to be the  CR $r$-jet of $\pi\circ F$. Notice that for an appropriate chart
 $\varphi_i^{-1}$ of projective space,

\begin{equation}\label{eq:rproj}j^r\pi(\sigma):=\Pi^{-1}_{k,x,i}(j^r_{D_h}(\varphi^{-1}_i\circ \pi\circ F)(x))\in  \mathcal{J}^r_{D}(M,\mathbb{C}^m)_i.
\end{equation}

The arguments in proposition \ref{pro:projjets} that showed that the bundles $\mathcal{J}^r_D(M,\mathbb{C}\mathbb{P}^m)$
 are well defined,  also prove that $j^r\pi(\sigma)$ is well defined independently of the approximately holomorphic
 coordinates and of the chart of  $\mathbb{C}\mathbb{P}^m$ we used; it is as well independent of the local frame
 of $E_k$, because the map is equivariant with respect to the action of $Gl(m+1,\mathbb{C})$ on the fibers of $E_k$ and on
  $\mathbb{C}\mathbb{P}^m$.

It is clear that  $j^r\pi$ is a submersion, and it is fiberwise holomorphic because in each fiber we have a map from
 some $\mathbb{C}^{m_1}$ to some  $\mathbb{C}^{m_2}$ (after composing with a chart $\varphi_i$), whose formula is
  that of the integrable case which is  holomorphic,  so item (1) holds.

We now prove the equality  $j^r_{D}(\pi \circ\tau_k)=j^r\pi(j^r_{D}\tau_k)$: let  $\varphi_i^{-1}$ be any chart whose
 domain contains $\pi\circ\tau_k(x)$. Then by the definition given in proposition  \ref{pro:projjets}
\[j^r_D(\pi\circ\tau_k)(x):=j^r_{D}(\varphi_i^{-1}\circ \pi \circ\tau_k)(x).\]
We just defined  in equation (\ref{eq:rproj})
\[j^r\pi(j^r_D\tau_k(x)):=\Pi^{-1}_{k,x,i}(j^r_{D_h}(\varphi^{-1}_i\circ \pi\circ F)(x)).\]
By proposition \ref{pro:projjets} the right hand side of the two previous equalities coincides, i.e.
 \[\Pi^{-1}_{k,x,i}(j^r_{D_h}(\varphi^{-1}_i\circ \pi\circ F)(x))=j^r_{D}(\varphi_i^{-1}\circ \pi \circ\tau_k)(x).\]
Here the holomorphic function $\varphi^{-1}_i\circ \pi\colon \mathbb{C}^{m+1}\backslash \{0\}\rightarrow \mathbb{C}^m$
 plays the role of $H$ in  proposition \ref{pro:projjets}. Also observe that the proposition is in principle
 only valid when  $\underline{\mathbb{C}}^{m+1}$ has the trivial connection. In the current situation
 $\underline{\mathbb{C}}^{m+1}$ is endowed with a diagonal connection coming from the one in $L^{\otimes k}$.
  The key point is that the composition $\varphi^{-1}_i\circ \pi\circ\phi_k$ is  a section of
  $\underline{\mathbb{C}}^m\otimes L^{\otimes k}\otimes  L^{-\otimes k}$ and hence a $\mathbb{C}^m$-valued
   function independently of the trivialization of $L^{\otimes k}$. Therefore the flat connection
   $\mathrm{d}$ on $\underline{\mathbb{C}}^m$ is induced from $\mathrm{d}\otimes\mathrm{I}+\mathrm{I}\otimes\nabla_k$ in
   $\underline{\mathbb{C}}^{m+1}\otimes L^{\otimes k}$, where $\nabla_k$ is any Hermitian connection
    on $L^{\otimes k}$. In other words, the equations of proposition  \ref{pro:projjets} involving
    the connection $\nabla_g\otimes \mathrm{I}+\mathrm{I}\otimes\mathrm{d}$ on ${(T^{*1,0}\mathbb{C}^n}^{\odot r})\otimes \underline{\mathbb{C}}^{m+1}$
     are also valid in this setting for the connection $\nabla_g\otimes \mathrm{I}+\mathrm{I}\otimes (\mathrm{d}\otimes\mathrm{I}+\mathrm{I}\otimes\nabla_k)$,
     and this finishes the proof of item (2).

The previous ideas work word by word to show  that for symplectic manifolds
$j^r\pi\colon \mathcal{J}^rE_k^*\rightarrow \mathcal{J}^r(P,\mathbb{C}\mathbb{P}^m)$ is a well defined
 submersion and that equation (\ref{eq:fundeq2}) holds.

If we have a distribution   $G$, once we use the local identification coming from approximately holomorphic
 coordinates adapted to $G$, the commutativity of the diagram   (\ref{cd:submersions}) follows from the
 commutativity in the holomorphic case. It is also clear that
 $j^r\pi\colon\mathcal{J}^r_GE_k^*\rightarrow \mathcal{J}^r_G(P,\mathbb{C}\mathbb{P}^m)$ is a submersion
 and that equation (\ref{eq:fundeq3}) holds.
\end{proof}

In order to  describe the linearized Thom-Boardman stratification we need to define  -at least for certain
 kinds of strata $\mathbb{P}S_k^a$ of $\mathcal{J}^r_D(M,\mathbb{C}\mathbb{P}^m)$- the corresponding subsets
  of transverse holonomy $\Theta_{\mathbb{P}S_k^a}$.

\begin{definition}\label{def:holon} Let $\mathbb{P}S_k$ be a sequence of strata of $\mathcal{J}^r_D(M,\mathbb{C}\mathbb{P}^m)$
  so that in canonical affine charts of  $\mathbb{C}\mathbb{P}^m$ and approximately holomorphic coordinates it is
   identified  with a stratum  $\mathbb{P}S$ of $\mathcal{J}^r_{D_h,n,m}$ invariant under the action of
   $\mathsf{T}\times Gl(n,\mathbb{C})$.
We let $\mathbb{P}S_{k,i}:=\mathbb{P}S_k\cap \mathcal{J}_{D}^r(M,\mathbb{C}^m)_i$ and then we define

\[\Theta_{\mathbb{P}S_k}:=\bigcup_{i\in\{0,\dots,m\}}\Theta_{\mathbb{P}S_{k,i}}.\]
For  $S_k:=j^r\pi^{-1}(\mathbb{P}S_k)$, with
$j^r\pi\colon \mathcal{J}^r_DE_k^*\rightarrow \mathcal{J}^r_D(M,\mathbb{C}\mathbb{P}^m)$
 the submersion of proposition \ref{pro:TBAsubmer}, we define  $\check{\Theta}_{S_k}:=j^r\pi^{-1}(\Theta_{\mathbb{P}S_k})$.

 In the relative theory we assume that for a choice of approximately holomorphic coordinates adapted
 to  $G$ and canonical affine charts of  projective space, the sequence
  $\mathbb{P}S_{k,i}\subset \mathcal{J}^r_G(P,\mathbb{C}^m_i)$ is identified with a stratum
    $\mathbb{P}S$ of $\mathcal{J}^r_{\mathbb{C}^g,p,m}=\mathcal{J}^r_{g,m}\times \mathbb{C}^{p-g}$
     invariant under the action of  $\mathsf{T}\times Gl(g,\mathbb{C})$. Then we define
     \[\Theta_{\mathbb{P}S_k}:=\bigcup_{i\in\{0,\dots,m\}}\Theta_{\mathbb{P}S_{k,i}}.\]
      For $S_k:=j^r\pi^{-1}(\mathbb{P}S_k)\subset \mathcal{J}^r_GE_k^*$, $S_k^G:={p_G}^{-1}(S_k)\subset\mathcal{J}^rE_k^*$,  we define the subset $\check{\Theta}_{S_k^G}\subset S_k^G$ by pulling back $\Theta_{\mathbb{P}S_k}$ to $\mathcal{J}^rE_k^*$ using either of the sides of the commutative diagram (\ref{cd:submersions}).
     \end{definition}

Notice that by item (1) of lemma \ref{lem:locrep} the subsets $\Theta_{\mathbb{P}S_{k,i}}$ are well defined, so definition \ref{def:holon} makes sense. It is also satisfactory because of the following result:
\begin{lemma}\label{lem:hol2} We have
\[ \Theta_{\mathbb{P}S_k}\cap \mathcal{J}^r_D(M,\mathbb{C}^m_i)=\Theta_{\mathbb{P}S_{k,i}}.\]

\[ \Theta_{\mathbb{P}S_k}\cap \mathcal{J}^r_G(P,\mathbb{C}^m_i)=\Theta_{\mathbb{P}S_{k,i}}.\]
\end{lemma}
\begin{proof}  Fix  approximately holomorphic coordinates and canonical affine charts  of $\mathbb{C}\mathbb{P}^m$,
 so that $\Pi_{k,x,i}(\mathbb{P}S_{k,i})=\mathbb{P}S$, for all $k,x,i$. We need to show is that
  \[j^r\Psi_{ji}(\Theta_{\mathbb{P}S})=\Theta_{\mathbb{P}S}\] in the domain of definition of $j^r\Psi_{ji}$, where $\Psi_{ji}$ is a change of canonical affine coordinates.

Let $\psi$ be an $r$-jet in $\Theta_{\mathbb{P}S}$. Then we have a lift  $\tilde{\psi}$ to
$\mathcal{J}^{r+1}_{D_h,n,m}$ and a local representation $\alpha$ of the lift cutting $\mathbb{P}S$
 transversally along $D_h$ at $\psi$. As we mentioned regarding transversality the local representation
  is essentially unique. That means in particular that any other representation $\alpha'$ will also share
   the transversality property. By definition  $\tilde{\psi}$ is the ($r+1$)-jet of a local CR function
    $F$. Then  $j^r_{D_h}F(0)=\psi$ and
    $(F(0),\mathrm{d}_{D_h}j^r_{D_h}F(0))=(F(0),\partial_0 j^r_{D_h}F(0))=j^{r+1}_{D_h}F(0)=\tilde{\psi}$.
     Thus, $j^r_{D_h}F$ is a local representation of $\tilde{\psi}$ which is transverse to $\mathbb{P}S$ along $D_h$ at $\psi$.

Since $j^{r+1}\Psi_{ji}(j^{r+1}_{D_h}F)=j^{r+1}_{D_h}(\Psi_{ij}\circ F)$, we deduce that $j^{r+1}\Psi_{ji}(\tilde{\psi})$
 is a lift of $j^r\Psi_{ji}(\psi)$ with local representation $j^{r}_{D_h}(\Psi_{ij}\circ F)$, which is
  obviously transverse  along $D_h$ to  $j^{r}\Psi_{ji}(\mathbb{P}S)=\mathbb{P}S$ because $j^{r}\Psi_{ji}$
  is a diffeomorphism that preserves the pullback of $D_h$ to $\mathcal{J}^r_{D_h,n,m}$. We just checked one inclusion, but that suffices because $\Psi_{ji}$ is a diffeomorphism, thus the result for jets along $D$ follows.

An analogous proof shows the desired result for jets along $G$.
\end{proof}

 The  linearized Thom-Boardman stratification is  the pullback to  $\mathcal{J}^r_DE_k^*$ by $j^r\pi$ of
  the analog of the  Thom-Boardman stratification of $\mathcal{J}^r_{D}(M,\mathbb{C}\mathbb{P}^m)$ (see
  for example \cite{Bo67}), together with the strata  $Z_k$. The definition is the natural extension
   of the one given for symplectic  manifolds by D. Auroux in \cite{Au01}.

A first rough definition of the stratification  of $\mathcal{J}^r_{D}(M,\mathbb{C}\mathbb{P}^m)$ is the
following: we fix approximately holomorphic coordinates and canonical affine charts of projective space,
 so we have charts $\Pi_{k,x,i}^{-1}\colon\mathcal{J}^r_{D_h,n,m}\rightarrow\mathcal{J}^r_D(M,\mathbb{C}^m)_i$.
In  each $\mathcal{J}^r_{D_h,n,m}$ there is a CR Thom-Boardman stratification which is
$\mathsf{T}\times(\mathcal{H}^r_n\times \mathcal{H}^r_m)$-invariant, where $\mathcal{H}^r_l$ is
the group of $r$-jets of germs of bi-holomorphic transformations from $\mathbb{C}^l$ to $\mathbb{C}^l$;
 in particular it is $\mathsf{T}\times Gl(n,\mathbb{C})$-invariant, so it defines a stratification
  on each  $\mathcal{J}^r_D(M,\mathbb{C}^m)_i$. The $\mathcal{H}^r_m$-invariance implies that the
  identifications that define $\mathcal{J}^r_D(M,\mathbb{C}\mathbb{P}^m)$ are compatible  with the
  aforementioned stratifications  on  $\mathcal{J}^r_D(M,\mathbb{C}^m)_i$.

Once we pullback the stratification to $\mathcal{J}^r_DE_k^*$ the behavior of the strata  when they
approach $Z_k$ needs to be clarified. To do that we redefine the stratification  as follows (see \cite{Au01}):

Given $\sigma \in \mathcal{J}^r_{D}E_k^*$
let us denote its image in
$\mathcal{J}^r_{D}(M,\mathbb{C}\mathbb{P}^m$)  by  $\phi=(\phi_0,\dots,\phi_r)$. Let us define
 \begin{equation}\label{eq:order1sing} \Sigma_{k,i}=\{\sigma\in \mathcal{J}^r_{D}E_k^*\,|\,\dim_\mathbb{C}\ker\phi_1=i\}.
 \end{equation}
If $\mathrm{max}(0,n-m)<i\leq n$, the strata  $\Sigma_{k,i}$ are smooth submanifolds whose boundary is
  the union    $\bigcup_{j>i}\Sigma_{k,j}$ together with a subset of   $Z_k\backslash\Theta_{Z_k}$.

Each $\Sigma_{k,i}$   is the  pullback of a stratum  $\mathbb{P}\Sigma_{k,i}\subset\mathcal{J}^r_{D}(M,\mathbb{C}\mathbb{P}^m)$,
 and the given description of their closure is easy to check.

For $r\geq 2$,
  define $\check{\Theta}_{\Sigma_{k,i}}$ as the subset of  $r$-jets $\sigma=(\sigma_0,\dots,\sigma_r)\in \Sigma_{k,i}$
  so that

\begin{equation}\label{eqn:strkernel}
\Xi_{k,i;\sigma}=\{u\in D\,|\, (i_u\phi,0)\in T_\phi \mathbb{P}\Sigma_{k,i}\}
\end{equation}
has the expected (complex) codimension  in   $D$, which is the
(complex) codimension of $\Sigma_{k,i}$ in $\mathcal{J}^r_{D}E_k$,
which equals the  codimension of $\mathbb{P}\Sigma_{k,i}$ in
$\mathcal{J}^r_{D}(M,\mathbb{C}\mathbb{P}^m)$.

The subset $\check{\Theta}_{\Sigma_{k,i}}$ is also the one coming from definition \ref{def:holon}: observe that
 $\Theta_{\mathbb{P}\Sigma_{k,i}}$  are exactly those points of  $\mathbb{P}\Sigma_{k,i}$ which have a
  lift with a transverse local representation. Since the term that we
add to the $r$-jet to define the lift is of order  $r+1>2$, the
transversality of the local representation does not depend on the
lift, that can be chosen to have vanishing  component of order $r+1$.

Fix as in the proof of lemma \ref{lem:locrep}  A.H. coordinates so that at the origin $(D\oplus D^\perp,J)=(D_h\oplus D_v,J_0)$
 and the induced connection form (on $\mathcal{J}^r_{D_h,n,m}$) is vanishing; fix also the canonical
  affine charts of $\mathbb{C}\mathbb{P}^m$. Then the strata $\mathbb{P}\Sigma_{k,i}$ are sent to
  the Thom-Boardmman stratum $\mathbb{P}\Sigma_{i}$ of $\mathcal{J}^r_{D_h,n,m}$. The local representation
   of $(\phi,0)$ can be taken to be a CR section  $\alpha$ of $\mathcal{J}^r_{D_h,n,m}$. The stratum
   $\mathbb{P}\Sigma_{i}$ is CR, therefore
\[T_{D_h}j^1_{D_h}\alpha(0)\cap (T\mathbb{P}\Sigma_{i}\cap\hat{D}_h)\]
is a complex subspace of $T\mathbb{C}^n$. Undoing the
identifications the previous subspace  goes to
 the subspace in equation (\ref{eqn:strkernel}). By definition of transversality along $D$,
 $\Theta_{\mathbb{P}\Sigma_{k,i}}$ are exactly those  $\phi$ for which $\Xi_{k,i;\sigma}$   has
  the codimension of  $\mathbb{P}\Sigma_{k,i}$ in  $\mathcal{J}^r_D(M,\mathbb{C}\mathbb{P}^m)$. By
   construction (equation (\ref{eqn:strkernel}))
\[\check{\Theta}_{\Sigma_{k,i}}=j^r\pi^{-1}(\Theta_{\mathbb{P}\Sigma_{k,i}}).\]
Hence $\check{\Theta}_{\Sigma_{k,i}}$ is the same subset introduced in definition \ref{def:holon}.

 If $p+1\leq r$, we define inductively
\[\Sigma_{k,i_1,\dots,i_p,i_{p+1}}=\{\sigma\in \Theta_{\Sigma_{k,i_1,\dots,i_p}}\,|\,
\dim_\mathbb{C}(\ker \phi_1\cap \Xi_{k,i_1,\dots,i_p;\sigma})=i_{p+1}\},\]
\noindent with
\[\Xi_{k,I;\sigma}=\{u\in D\,|\, (i_u\phi,0)\in T_\phi\mathbb{P}\Sigma_{k,I}\}.\]
As in the previous case we define $\check{\Theta}_{\Sigma_{k,I}}$ either as the points such that
  the complex codimension of  $\Xi_{k,I;\sigma}$ in $D$ is the same as the codimension of $\Sigma_{k,I}$
   in  $\mathcal{J}^r_{D}E_k$, or as the pullback of $\Theta_{\mathbb{P}\Sigma_{k,I}}$.

If $i_1\geq\cdots\geq i_{p+1}\geq 1$, $\Sigma_{k,i_1,\dots,i_{p+1}}$ is -in the local model- a smooth CR
 submanifold whose closure in $\Sigma_{k,i_1,\dots,i_p}$ is the union of the  $\Sigma_{k,i_1,\dots,i_p,j}, j>i_{p+1}$,
  and a subset of $\Sigma_{k,i_1,\dots,i_p}\backslash\check{\Theta}_{\Sigma_{k,i_1,\dots,i_p}}$ \cite{Bo67}.
   The problem is that for large values of $r,n,m$, the closure of the strata in $\mathcal{J}^r_{D_h,n,m}$ is
    hard to understand, and what we have defined -once $Z_k$ has been added-  might very well not be a
     Whitney (A) quasi-stratification. More precisely, let
     $\Sigma_{m+1;q}:=\Sigma_{m+1,1,\overset{(q)}{\dots},1}\subset \mathcal{J}^r_{D_h,n,m}$ be a
     so called Morin stratum. Then in \cite{Wi79} it is shown that
 \[\overline{\Sigma_{m+1;q}}\cap \Sigma_{m+2,0}\neq \emptyset,\]
 but for $q$ large enough $\mathrm{dim}\Sigma_{m+1;q}<\mathrm{dim}\Sigma_{m+2,0}$, thus Whitney's
 condition (A) can never hold. It is known that $\mathcal{J}^r_{D_h,n,m}$ admits a Whitney (A) stratification
 containing the Morin strata. If the dimensions satisfy $n<4$ or $2n>3m-4$, then a generic function will
 avoid $\Sigma_{m+2,0}$ and $\Sigma_{m+1,2}$ and therefore will only intersect the Morin strata,
 so the aforementioned previous stratification suffices (also  because the strata $\Sigma_{k,I}$ do
  not accumulate in points of  $\Theta_{Z_k}$). In general one must refine the Thom-Boardman stratification.

Recall that using the local identifications  the stratification we have defined (minus  $Z_k$) is the union
 running over the affine charts of the  pullback by
 $j^r(\varphi^{-1}_i\circ \pi)\colon \mathcal{J}^{r}_{D_h,n,m+1}\backslash Z\rightarrow \mathcal{J}^r_{D_h,n,m}$
  of the CR Thom-Boardman stratification $\mathbb{P}\Sigma$ of $\mathcal{J}^r_{D_h,n,m}$. The latter is
   CR and $\mathsf{T}\times (Gl(n,\mathbb{C})\times \mathcal{H}^r_m)$-invariant.

 On the domain of each chart  $\mathcal{J}^r_{D_h,n,m}$ we can use the results of Mather  \cite{Ma75}
  to refine $\mathbb{P}\Sigma$  into  a  CR finite, Whitney (A) stratification transverse to the fibers
  and invariant under the action of  $\mathsf{T}\times (Gl(n,\mathbb{C})\times \mathcal{H}^r_m)$,
  and such that the submanifolds  $\mathbb{P}\Sigma_{I}$ are unions of strata of the refinement. Due
   to the required invariance properties for the refinements,  they can be glued to give a refinement
   of the stratification  $\mathbb{P}\Sigma_k\subset\mathcal{J}^r_D(M,\mathbb{C}\mathbb{P}^m)$, which
    is independent of the choice of approximately holomorphic coordinates. Thus, its pullback is a
    finite, Whitney (A) stratification of   $\mathcal{J}^r_DE_k^*$  and such that the $\Sigma_{k,I}$ are
    union of strata. It is by construction invariant by the action of $Gl(m+1,\mathbb{C})$ on the fiber.

It is important to notice that since all the strata are contained in the closure of  $\Sigma_{k,\mathrm{max}(0,n-m)+1}$,
 they accumulate near $Z_k$ in points of   $Z_k\backslash\Theta_{Z_k}$. Therefore, by adding  $Z_k$ we
 obtain a quasi-stratification of
$\mathcal{J}^r_{D}E_k$.

If we have a distribution $G$ we use exactly the same definitions but in the subbundles $\mathcal{J}^r_{G}E_k$
 and $\mathcal{J}^r_{G}(P,\mathbb{C}\mathbb{P}^m)$. That is, we have the strata

\[\mathbb{P}\Sigma_{k,i}=\{\phi\in \mathcal{J}^r_G(P,\mathbb{C}\mathbb{P}^m)\,|\,\dim_\mathbb{C}\ker\phi_1=i\}\]
and for $r\geq 2$,
  $\check{\Theta}_{\mathbb{P}\Sigma_{k,i}}\subset \mathbb{P}\Sigma_{k,i}$ is the subset of  $r$-jets
  along $G$, $\phi=(\phi_0,\dots,\phi_r)$
  so that
  \begin{equation}\label{eqn:strkernel2}
\Xi_{k,i;\sigma}=\{u\in G\,|\, (i_u\phi,0)\in T_\phi \mathbb{P}\Sigma_{k,i}\}
\end{equation}
has the expected (complex) codimension  in   $G$, which is the (complex) codimension of
$\mathbb{P}\Sigma_{k,i}$ in  $\mathcal{J}^r_{G}(P,\mathbb{C}\mathbb{P}^m)$.

The subsets $\mathbb{P}\Sigma_{k,I}$ are defined similarly. The result is  a stratification $\mathbb{P}\Sigma_k$
 of $\mathcal{J}^r_G(P,\mathbb{C} \mathbb{P}^{m})$. In charts adapted to $G$ as in the proof of lemma
  \ref{lem:locrep} and affine charts  -in which $\mathcal{J}^r_{\mathbb{C}^g,p,m}=\mathcal{J}^r_{g,m}\times \mathbb{C}^{p-g}$-,
   the induced stratification $\mathbb{P}\Sigma$  is seen to be the leafwise    Thom-Boardman stratification,
   i.e. the  Thom-Boardman stratification  of $\mathcal{J}^r_{g,m}$ multiplied by   $\mathbb{C}^{p-g}$.

Using the lower part of the commutative diagram (\ref{cd:submersions}), we pull back $\mathbb{P}\Sigma_k$
 to $\mathbb{P}\Sigma^G_k\subset \mathcal{J}^r(P,\mathbb{C} \mathbb{P}^{m})$.
Let  $\Sigma^G_k$ be the pullback of $\mathbb{P}\Sigma^G_k$ to $\mathcal{J}^rE_k^*$. To refine it  we first
locally refine  $\mathbb{P}\Sigma_k$ as follows: we go the leafwise Thom-Boardman stratification  furnished
by  the previous A.H. coordinates and  affine charts and construct a holomorphic
 $\mathsf{T}\times (Gl(g,\mathbb{C})\times \mathcal{H}^r_m)$-invariant  refinement in one of the
 leaves of $\mathcal{J}^r_{\mathbb{C}^g,p,m}$ (which is identified with $\mathcal{J}^r_{g,m}$).
 Next we extend it independently of the remaining $p-g$ complex coordinates $z_k^{g+1},\dots,z_k^p$.
  The local refinements of the leafwise Thom-Boardman stratification glue well and thus define a
  sequence of Whitney (A) stratifications  $\mathcal{J}^r_G(P,\mathbb{C}\mathbb{P}^m)$, which does not
   depend either on the A.H. coordinates adapted to $G$  or in the chosen affine charts of  $\mathbb{C}\mathbb{P}^m$.
     Its pullback to $\mathcal{J}^rE_k^*$ refines $\Sigma^G_k$ to a sequence of Whitney (A) stratifications.

\begin{definition}\label{def:TBA}(see \cite{Au01}). \quad
\begin{enumerate}
\item
Given $(M,D,J,g_k)$ and $E_k=\underline{\mathbb{C}}^{m+1}\otimes L^{\otimes k}$, the Thom-Boardman-Auroux
  stratification of $\mathcal{J}^r_{D}(M,\mathbb{C}\mathbb{P}^m)$, denoted by $\mathbb{P}\Sigma_k$,  is
   the stratification (or rather its refinement) built out of the pieces of the Thom-Boardman stratifications
    of $\mathcal{J}^r_{D_h,n,m}$. The Thom-Boardman-Auroux  quasi-stratification of  $\mathcal{J}^r_{D}E_k$
     is the pullback of the Thom-Boardman-Auroux stratification of  $\mathcal{J}^r_{D}(M,\mathbb{C}\mathbb{P}^m)$
      together with the zero section. We denote it by $\Sigma_k$.

\item Given $(P,J,G,g_k)$ and $E_k=\underline{\mathbb{C}}^{m+1}\otimes L^{\otimes k}$, the Thom-Boardman-Auroux
 stratification of
$\mathcal{J}^r(M,\mathbb{C}\mathbb{P}^m)$ along $G$, denoted by $\mathbb{P}\Sigma^G_k$,  is the stratification
 (or rather its refinement) built out of the pieces of the Thom-Boardman stratifications of $\mathcal{J}^r_{\mathbb{C}^g,p,m}$.
   The Thom-Boardman-Auroux  quasi-stratification of  $\mathcal{J}^rE_k$ along $G$, that we denote by $\Sigma^G_k$,
   is the pullback of the Thom-Boardman-Auroux stratification of  $\mathcal{J}^r(M,\mathbb{C}\mathbb{P}^m)$
   along $G$ together with $Z_k$.
\end{enumerate}
\end{definition}

\begin{lemma}\label{lem:TBA} The Thom-Boardman-Auroux quasi-stratification of $\mathcal{J}^r_DE_k$ and the
 Thom-Boardman-Auroux quasi-stratification of $\mathcal{J}^rE_k$ along $G$ are finite, Whitney (A), and approximately holomorphic.
\end{lemma}

\begin{proof} We start with jets along $D$. The description of the closure of the strata inside $Z_k$ implies
that the quasi-stratification condition holds.

The delicate point is checking that the strata are approximately holomorphic (for the modified connection).

First we study the sequence $Z_k$. Though for this sequence the approximate holomorphicity is obvious, we will
 give a proof that works for other sequences of strata. Indeed, by lemma \ref{lem:trivializationweak} the
  sequence of zero sections $Z_k\subset E_k$ is as required. If we prove that the natural projections
\[ \pi^{r}\colon \mathcal{J}^r_DE_k\rightarrow E_k\]
are an A.H. sequence of maps which is also $\epsilon$-transverse for some $\epsilon>0$, then the composition
of the local maps defining $Z_k\subset E_k$ with the projection $\pi^r$ are local functions for
${(\pi^{r})}^{-1}(Z_k)=Z_k\subset  \mathcal{J}^r_DE_k$ meeting the conditions of definition \ref{def:stratification}.

 More generally we prove that the natural projection  $\pi^{r}_{r-h}\colon \mathcal{J}^r_DE_k\rightarrow \mathcal{J}^{r-h}_DE_k$
  is approximately holomorphic: we fix A.H. coordinates and A.H. reference frames $j^r_D\tau_{k,x,I}^{\mathrm{ref}}$
    of  $\mathcal{J}^r_DE_k$ (resp. $j^{r-h}_D\tau_{k,x,I'}^{\mathrm{ref}}$  of  $\mathcal{J}^{r-h}_DE_k$) as in
     equation (\ref{eqn:holonframe}). Recall that proposition \ref{pro:perthol} implies that the sequences are indeed
      A.H. Using these frames we obtain A.H. coordinates $z_k^1,\dots,z_k^n,u_k^I,s_k$ (resp. $z_k^1,\dots,z_k^n,v_k^{I'},s_k$)
       for the total space of $\mathcal{J}^r_DE_k$ (resp. $\mathcal{J}^{r-h}_DE_k$). From the  equality
 \begin{equation}\label{eq:projrjet}
  \pi^{r}_{r-h}(j^r_D\tau_{k,x,I}^{\mathrm{ref}})=j^{r-h}_D\tau_{k,x,I}^{\mathrm{ref}}
  \end{equation}
 we deduce $\pi^{r}_{r-h}(j^r_D\tau_{k,x,I}^{\mathrm{ref}})=W_{I}(z_k,v_k^{I'})$, where $W_I(z_k,v_k^{I'})$ is
 A.H. with respect to the canonical CR structures associated to the coordinates. This, together with the fiberwise linearity
  of $\pi^{r}_{r-h}$ imply that in these coordinates $\pi^{r}_{r-h}$ is A.H., and hence it is A.H. with respect to the almost
  CR structures of the total spaces. It is also straightforward from equation (\ref{eq:projrjet}) that the projections
   are $\epsilon$-transverse (another way is to use rather than holonomic frames, the frames $\mu_{k,x,I}$ of equation
    (\ref{eqn:canframe}). They are also frames for the modified metric because of for example remark \ref{rem:metriccomparison},
     therefore one can check estimated transvesality using them, something which is straightforward).

We would like to do something similar with the strata $\Sigma_{k,I}$ and the projection
$j^r\pi\colon \mathcal{J}^r_DE_k^*\rightarrow \mathcal{J}^r_D(M,\mathbb{C}\mathbb{P}^m)$
(away from a uniform tubular neighborhood of the zero section, where the differential goes to infinity).
 The image of a trivialization  $j^r_D\tau_{k,x,I}^{\mathrm{ref}}$ is $j^r_D(\pi\circ\tau_{k,x,I}^{\mathrm{ref}})$,
  also  approximately holomorphic. The map is equally fiberwise holomorphic, but the difference is the
   non-linearity of the restriction to the fibers.

We adopt a different strategy that amounts to perturbing the almost CR structures into integrable ones and
then checking that $j^r\pi$ is CR with respect to them:  we take Darboux charts and trivialize $L^{\otimes k}$ with
 a unitary section $\xi_k$ whose associated connection form in the domain of Darboux charts is $A$. Next we
  trivialize  $\mathcal{J}^r_DE_k$ with the frames $\mu_{k,x,I}$ of equation (\ref{eqn:canframe}), but using
  $\xi_k$ tensored with a basis of $\mathbb{C}^{m+1}$ to trivialize $\underline{\mathbb{C}}^{m+1}\otimes L^{\otimes k}$.
   In this way $\mathcal{J}^r_DE_k$ becomes the trivial bundle $\mathcal{J}^r_{D_h,n,m+1}$ (with is canonical
   trivialization constructed out of $dz_k^1,\dots,dz_k^n$). Let us use in the base the canonical CR structure
   $(D_h,J_0)$. Proposition   \ref{pro:perthol} in the integrable case (and for curvature of type (1,1) and with
    trivial derivative, as it is the case in Darboux coordinates) implies that the modified connection defines
     a new CR structure in the total space of $\mathcal{J}^r_{D_h,n,m}$; let $(\hat{D}_h,\bar{J}_0)$ be the
     corresponding distribution and almost complex structure, and   let $(\hat{D},\hat{J})$ be the distribution
      and almost complex structure induced by the almost CR structure of $\mathcal{J}^r_DE_k$. If in the fiber
       of $\mathcal{J}^r_{D_h,n,m+1}$ we fix a ball $B(\sigma, R)$, then in $B(0,\rho)\times B(\sigma, R)$ the
        Euclidean metric is comparable with the metric carried by $\mathcal{J}^r_DE_k$. More important
\begin{equation}\label{eq:compcr}
|\mathrm{d}^j(\hat{D}-\hat{D}_h)|_{g_0}\leq O(k^{-1/2}),\; j\geq 0.
\end{equation}
If we use the orthogonal projection to push  $\hat{J}$ into $\hat{J}_h\colon \hat{D}_h\rightarrow\hat{D}_h$
 we also have
\begin{equation}\label{eq:compcr2}
|\mathrm{d}^j(\hat{J}_h-\bar{J}_0)|_{g_0}\leq O(k^{-1/2}),\; j\geq 0.
\end{equation}
 We use the same Darboux charts for   $\mathcal{J}^r_{D_h}(\mathbb{C}^n\times \mathbb{R},\mathbb{C}\mathbb{P}^m)$,
  so locally and using canonical affine charts we have identifications with $\mathcal{J}_{D_h,n,m}^r$. This is
  a trivial vector bundle (again using the basis induced by $dz_k^1,\dots,dz_k^n$ and the basis of $\mathbb{C}^m$).
   We fix the product CR structure and denote by $(\tilde{D}_h,\tilde{J}_0)$ the distribution and almost complex
    structure. Let $(\tilde{D},\tilde{J})$ be the distribution and almost complex structure induced by the almost
     CR structure of $\mathcal{J}^r_{D}(M,\mathbb{C}\mathbb{P}^m)$. By construction

\begin{equation}\label{eq:compcr3}
|\mathrm{d}^j(\tilde{D}-\tilde{D}_h)|_{g_0}, |\mathrm{d}^j(\tilde{J}_h-\tilde{J}_0)|_{g_0}\leq O(k^{-1/2}),\; j\geq 0,
\end{equation}
where $\tilde{J}_h$ is the almost complex structure on $\tilde{D}_h$ defined out of $\tilde{J}$ and the orthogonal projection.

Equations (\ref{eq:compcr}),  (\ref{eq:compcr2}), (\ref{eq:compcr3}) imply that if
$j^r(\varphi^{-1}_i\circ\pi )\colon \mathcal{J}^r_{D_h,n,m+1}\rightarrow \mathcal{J}_{D_h,n,m}^r$ is
 CR with respect to $(\hat{D}_h,\bar{J}_0)$ and $(\tilde{D}_h,\tilde{J}_0)$, then it is almost CR with respect to  the global almost CR structures.

The map $j^r(\varphi^{-1}_i\circ\pi)\colon \mathcal{J}^r_{D_h,n,m+1}\rightarrow \mathcal{J}_{D_h,n,m}^r$
 is exactly the same as in the holomorphic (or rather CR) models. It is CR with respect to the aforementioned CR
  structures because it preserves the foliations, it is fiberwise holomorphic and sends ``enough'' CR
  sections of $\mathcal{J}^r_{D_h,n,m+1}$ to CR sections of $\mathcal{J}^r_{D_h,n,m}$. To be more precise,
   for any point $\sigma\in \mathcal{J}^r_{D_h,n,m+1}$ and any vector $v$ in its tangent space along
   the leaf and not tangent to the fiber, we can find a CR section $F$ whose CR $r$-jet in $x$ is $\sigma$
    and such that the tangent space to its graph contains  $v$. Since
    $j^r(\varphi^{-1}_i\circ\pi)(j^r_{D_h}F)=j^r_{D_h}(\varphi^{-1}_i\circ\pi\circ F)$  is also a CR
    section, we deduce that  $j^r(\varphi^{-1}_i\circ\pi)_*(\bar{J}v)=\tilde{J}_0(j^r(\varphi^{-1}_i\circ\pi)_*(v))$.

The strata $\mathbb{P}\Sigma_{k}$ (or rather of its refinement)
-once we choose  A.H. coordinates and affine charts of projective
space- are identified with the strata of  (the refinement of) the CR
Thom-Boardman stratification of $\mathcal{J}_{D_h,n,m}^r$, which are
CR. The comparison between the $(\hat{D}_h,\bar{J}_0,g_0)$ and the
original almost CR structure implies that the strata of
$\mathbb{P}\Sigma_{k}$ are A.H., and hence
$\Sigma_k=j^r\pi^{-1}(\mathbb{P}\Sigma_{k})$ is A.H. That the
projections are $\epsilon$-transverse is also clear, therefore the
desired result follows.

In the almost complex setting
 $j^r\pi\colon \mathcal{J}^rE_k^*\rightarrow \mathcal{J}^r(P,\mathbb{C}\mathbb{P}^m)$ is equally
 shown  to be approximately holomorphic away from a uniform neighborhood of the zero section. In the relative case,
   and for a sequence of A.H. strata $\mathbb{P}S_k$ fulfilling the conditions of definition  \ref{def:holon},
    the approximate holomorphicity of  ${p_G}^{-1}{j^r\pi}^{-1}S_k$ follows from the commutativity
    of the diagram  \ref{cd:submersions}, and from the approximate holomorphicity of
    $j^r\pi\colon \mathcal{J}^rE_k^*\rightarrow \mathcal{J}^r(P,\mathbb{C}\mathbb{P}^m)$ and of
      $p_G\colon \mathcal{J}^r(P,\mathbb{C}\mathbb{P}^m)\rightarrow \mathcal{J}^r_G(P,\mathbb{C}\mathbb{P}^m)$.
       Recall that the strata $\mathbb{P}\Sigma_k$ come from holomorphic models (the refinement of the strata
        of the leafwise Thom-Boardman stratification), so they are A.H. But $\Sigma_k^G$  is not truly a
        quasi-stratification of $\mathcal{J}^rE_k$. To be more precise it is not true that  the strata
        only accumulate in points of  $Z_k\backslash\Theta_{Z_k}\subset Z_k$, but it is still true that the
         points of $Z_k$ in which the other strata accumulate are never hit by a section transverse to $Z_k$
         along $G$. Thus, the Whitney type reasoning can be applied as long as we work with $r$-jets along $G$
          (see the proof of theorem \ref{thm:main2}).
\end{proof}

\begin{remark}\label{rm:zerosect}
Notice that we only conclude that the strata different form the zero section are approximately holomorphic
uniformly far from $Z_k$. This is enough for our purposes, for once we obtain transversality to $Z_k$ our
 $r$-jet will be  uniformly far from $Z_k\backslash \Theta_{Z_k}$. All the remaining strata approach $Z_k$
  accumulating only on points of $Z_k\backslash \Theta_{Z_k}$. Therefore, the $r$-jet will only hit them
  outside of a uniform tubular neighborhood of $Z_k$, where the approximate holomorphicity holds.
\end{remark}

\begin{definition}\label{def:gensec}\quad
\begin{enumerate}
\item An A.H. sequence of sections   of $E_k\rightarrow (M,D,J,g_k)$ is said to be $r$-generic if  its pseudo-holomorphic
$r$-jet  is uniformly transverse along $D$ to the  Thom-Boardman-Auroux quasi-stratification of $\mathcal{J}^r_DE_k$.

\item An A.H. sequence of sections   of $E_k\rightarrow (P,J,G,g_k)$ is said to be $r$-$G$-generic over $M$
 if its pseudo-holomorphic $r$-jet  is uniformly transverse over $M$ to $\Sigma_k^G\subset\mathcal{J}^rE_k$.

\item Let $\phi_k\colon M\backslash B_k\rightarrow \mathbb{C}\mathbb{P}^m$   be  sequence of  functions which
is A.H. outside of a uniform tubular neighborhood of $g_k$-radius
$\eta>0$ of $B_k$.  It is said to be
 $r$-generic if for $k$ large enough $B_k$ is a codimension $2(m+1)$ calibrated submanifold and
 $j^r_D\phi_k\colon M\backslash B_k\rightarrow \mathcal{J}^r_D(M\backslash B_k, \mathbb{C}\mathbb{P}^m)$
 is uniformly transverse along $D$ to the Thom-Boardman-Auroux stratification. Moreover, it is required
  to intersect the strata of strictly positive codimension out of a tubular neighborhood of
   $B_k$ of $g_k$-radius $\eta$.
\end{enumerate}

\end{definition}

\begin{lemma}\label{lem:rgen} Let $\tau_k$  be an A.H. sequence of sections of  $E_k\rightarrow (M,D,J,g_k)$.
 Then if $\tau_k$ is $r$-generic its projectivization $\phi_k\colon M\backslash \tau_k^{-1}(Z_k)\rightarrow \mathbb{C}\mathbb{P}^m$
  is also $r$-generic.
\end{lemma}
\begin{proof}It is elementary from the construction of the Thom-Boardman-Auroux (quasi)-stratifications
of $\mathcal{J}^r_DE_k$ and $\mathcal{J}^r_D(M, \mathbb{C}\mathbb{P}^m)$, proposition \ref{pro:TBAsubmer}
 relating $j^r_D\tau_k$, and $j^r_D\phi_k$ and lemma \ref{lem:TBA}.

Uniform transversality of $\tau_k$ to $Z_k$ implies by remark \ref{rm:zerosect} that $\phi_k$ intersects
 the remaining strata uniformly away from the zero set.  Estimated transversality along $D$ is also preserved
  when composed with $j^r\pi$ uniformly away from $Z$; the key point is selecting appropriate local A.H. defining
  functions for the strata:  in A.H. coordinates and affine charts  $\mathbb{P}\Sigma_{k,I}$ corresponds
   to a CR stratum $\mathbb{P}\Sigma_I$. Let $f$ be a local CR function defining it. Then
   $f\circ \Pi_{k,x,i}\circ j^r(\varphi^{-1}_i\circ\pi)$ are local defining functions for $\Sigma_{k,I}$.
    Now lemma \ref{lem:localchartrans} implies that local uniform estimated transversality along $D$ of
    $j^r_D\tau_k$ to $\Sigma_{k,I}$  is equivalent to uniform transversality along $D$ to {\bf 0}
     of $f\circ j^r(\varphi^{-1}_i\circ\pi)\circ j^r_D\tau_k=f\circ j^r_D(\varphi^{-1}_i\circ\phi_k)$.
     Again by the same lemma this is equivalent to uniform transversality  along $D$ of $j^r_D\phi_k$  to
      $\mathbb{P}\Sigma_{k,I}$. The case of the points close to the boundary of the strata is just a problem
       in a vector space; it follows from
       $j^r(\varphi^{-1}_i\circ\pi)\colon \mathcal{J}^{r}_{D_h,n,m+1}\backslash Z\rightarrow \mathcal{J}^r_{D_h,n,m}$
        being a submersion which amounts to suppressing coordinates of the fiber of $\mathcal{J}^{r}_{D_h,n,m+1}$
        (and because the metrics in these coordinates are comparable with the ambient metric, so the projection
         is $\epsilon$-transverse).
\end{proof}

 Let  $(P,\Omega)$ be a symplectic manifold with $(M,D,\omega:=\Omega_{\mid M})$ 2-calibrated and  $G$
  a local $J$-complex distribution extending  $D$. Let  $\tau_k$ be an A.H.  sequence of sections of  $E_k$ and denote by ${\phi_k}$ its  projectivization  away from its zero set.

\begin{proposition}\label{pro:relgen} Using the above notation, if $j^r\tau_k\colon P\rightarrow \mathcal{J}^rE_k$ is
  uniformly transverse over $M$ to $\Sigma_k^G\subset\mathcal{J}^r_GE_k $ then  ${\phi_k}_{\mid M}$ is $r$-generic.
\end{proposition}
\begin{proof}
We will make extensive use of diagram (\ref{cd:submersions})
\[\begin{CD}\mathcal{J}^rE_k^* @>{p_G}>>  \mathcal{J}^r_GE_k^*\\
@VV{j^r\pi}V         @VV{j^r\pi}V     \\
\mathcal{J}^r(P,\mathbb{C}\mathbb{P}^m) @>{p_G}>>  \mathcal{J}^r_G(P,\mathbb{C}\mathbb{P}^m)
\end{CD}
\]

\emph{Step 1:} Study the compatibility of the Thom-Boardman-Auroux stratifications with the identification
 of $\mathcal{J}^r_D(M,\mathbb{C}\mathbb{P}^m)$ with ${\mathcal{J}^r_G(P,\mathbb{C}\mathbb{P}^m)}_{\mid M}$.

At the points of $M$ there is a canonical $J$-complex identification between  $D$ and $G$,  inducing isometries

\[\Lambda_{k,i}\colon \mathcal{J}^r_D(M,\mathbb{C}\mathbb{P}^m)\rightarrow  {\mathcal{J}^r_G(P,\mathbb{C}\mathbb{P}^m)}_{\mid M}.\]

Let $z_k^1,\dots,z_k^p$ be A.H coordinates adapted to $(M,G)$. We can rewrite them as
$z_k^1,\dots,z_k^n,x_k^{2n+1},x_k^{2n+2},z_k^{n+2},\dots,z_k^p$, where $z_k^1,\dots,z_k^n,x_k^{2n+1}$ are
 by lemma \ref{lem:ahm} A.H. coordinates for $M$. Using also  the canonical affine charts of projective space we have

\begin{align*}
 \Pi_{k,x,i}^D\colon \mathcal{J}^r_D(M,\mathbb{C}^m)_i&\rightarrow \mathcal{J}^r_{D_h,n,m}=\mathcal{J}^r_{n,m}\times \mathbb{R},\\ \Pi_{k,x,i}^G\colon \mathcal{J}^r_G(P,\mathbb{C}^m)_i&\rightarrow  \mathcal{J}^r_{\mathbb{C}^n,p,m}=\mathcal{J}^r_{n,m}\times \mathbb{C}^{p-n},\\
\end{align*}
and a canonical identification in $\mathbb{C}^n\times\mathbb{R}\subset\mathbb{C}^p$

\[\Lambda\colon \mathcal{J}^r_{D_h,n,m}\rightarrow {\mathcal{J}^r_{\mathbb{C}^n,p,m}}_{\mid \mathbb{C}^n\times\mathbb{R}}. \]
The construction of  $\Pi_{k,x,i}^D,\Pi_{k,x,i}^G$ (see  equation (\ref{eq:linidentf}) and the last paragraph in the
 proof of lemma \ref{lem:locrep}) implies the commutativity of

\begin{equation}\label{cd:TBA}\begin{CD}\mathcal{J}^r_D(M,\mathbb{C}\mathbb{P}^m) @>{\Lambda_k}>>{\mathcal{J}^r_G(P,\mathbb{C}\mathbb{P}^m)}_{\mid M} \\
@VV{\Pi_{k,x,i}^D}V         @VV{\Pi_{k,x,i}^G}V     \\
 \mathcal{J}^r_{D_h,n,m}@>{\Lambda}>> {\mathcal{J}^r_{\mathbb{C}^n,p,m}}_{\mid \mathbb{C}^n\times\mathbb{R}}
\end{CD}
\end{equation}
The restriction of $\mathcal{J}^r_{\mathbb{C}^n,p,m}$ to $\mathbb{C}^n\times\mathbb{R}\approx M$ coincides with
$\mathcal{J}^r_{n,m}\times \mathbb{R}=\mathcal{J}^r_{D_h,n,m}$.

The identification $\Lambda$ obviously preserves the  Thom-Boardman-Auroux stratifications (and even the refinements),
 and hence so $\Lambda_k$ does.

\emph{Step 2:} Check that $\Lambda_k^{-1}\circ {(j^r_G\phi_k)}_{\mid M}\approxeq j^r_D({\phi_k}_{\mid M})$.

Since $\Lambda_k$ are $J$-complex isometries preserving the Thom-Boardman-Auroux stratifications we omit them from now on.

 By using the charts $\Pi_{k,x,i}^D,\Pi_{k,x,i}^G$ it is easy to see  that for any $j\in\{1,\dots,r\}$, the degree $j$
  homogeneous component of  $j^r_D({\phi_k}_{\mid M})$ approximately coincides with $\nabla^j_D({\phi_k}_{\mid M})$.
   Similarly,  the degree $j$ homogeneous component of  $j^r_G{\phi_k}$  approximately coincides with $\nabla^j_G{\phi_k}$.
    The result follows because we also have
\[(\nabla^j_G{\phi_k})_{\mid M}\approxeq \nabla^j_D({\phi_k}_{\mid M}).\]

\emph{Step 3:} Analyze the behavior of $j^r_D({\phi_k}_{\mid M})$ near  the set of base points $B_k$.

Since $Z_k\subset\mathcal{J}^rE_k$ is an A.H. sequence of submanifolds and $j^r\tau_k$ an A.H. sequence of sections,
 by corollary \ref{cor:fulltrans} uniform transversality over $M$ is equivalent to uniform transversality along
  $G$ at the points of $M$. In A.H. coordinates adapted to $G$, we are saying that the matrix of partial derivatives
   of $\tau_k$ with respect to $z_k^1,\dots,z_k^g$ has maximum rank and norm greater than some $\eta>0$. But this
    is equivalent to saying that  is uniformly transverse to $Z_k^G$, the pullback of the zero section of $\mathcal{J}^r_GE_k$.

By construction $\Sigma_k^G\backslash Z_k=p_G^{-1}j^r\pi^{-1}(\mathbb{P}\Sigma_k)=p_G^{-1}(\Sigma_k\backslash Z_k)$,
 and the strata of $\Sigma_k^G\backslash Z_k$  when approaching the zero section   accumulate into $p_G^{-1}(\Theta_{Z_k})$,
  where here $\Theta_{Z_k}\subset\mathcal{J}^r_GE_k$. Therefore  $j^r\tau_k$ intersects the strata of $\Sigma_k^G\backslash Z_k$
   away from a tubular neighborhood in $P$ (and hence in $M$) of radius $\eta'$ of $B_k$, the zero set of $j^r\tau_k$.
   Thus  ${(j^r\phi_k)}_{\mid M}={(j^r\pi(j^r\tau_k))}_{\mid M}$ intersects the  strata of $\mathbb{P}\Sigma_k^G$
   away from a tubular neighborhood in $M$ of radius $\eta'$ of $B_k$.

In general  $p_G(j^r\phi_k)\neq j^r_G\phi_k$ but using A.H. coordinates it is easy to check that
 $p_G(j^r\phi_k)\approxeq j^r_G\phi_k$. Hence, $j^r_G\phi_k$ intersects the strata of
 $\mathbb{P}\Sigma_k\subset \mathcal{J}^r_G(P,\mathbb{C}\mathbb{P}^m)$ away from a tubular neighborhood
  in $M$ of radius $\eta'$ of $B_k$, for all $k \gg 1$.

By steps 1 and 2 we deduce that  $j^r_D({\phi_k}_{\mid M})$ intersects the strata of
$\mathbb{P}\Sigma_k\subset \mathcal{J}^r_D(P,\mathbb{C}\mathbb{P}^m)$ away from a tubular neighborhood in $M$
of radius $\eta'$ of $B_k$, for all $k \gg 1$.

\emph{Step 4:} Relate uniform transversality over $M$ of $j^r\tau_k$ to $\Sigma_k^G\backslash Z_k$ with
uniform transversality along $D$  of $j^r_D({\phi_k}_{\mid M})$  to
 $\mathbb{P}\Sigma_k\subset   \mathcal{J}^r_D(M,\mathbb{C}\mathbb{P}^m)$.

The same ideas used in the proof of lemma \ref{lem:rgen} combined with $p_G(j^r\phi_k)\approxeq j^r_G\phi_k$,
  show that  uniform transversality over $M$ of $j^r\tau_k$  to $\Sigma_k^G\backslash Z_k$ is equivalent
  to uniform transversality over $M$ of ${j^r_G\phi_k}$  to
   $\mathbb{P}\Sigma_k\subset   \mathcal{J}^r_G(P,\mathbb{C}\mathbb{P}^m)$.

 Uniform transversality over $M$ of ${j^r_G{\phi_k}}$ to $\mathbb{P}\Sigma_k\subset   \mathcal{J}^r_G(P,\mathbb{C}\mathbb{P}^m)$
  is comparable to uniform transversality of ${(j^r_G{\phi_k})}_{\mid M}$ to
  ${\mathbb{P}\Sigma_k}_{\mid M}\subset   {\mathcal{J}^r_G(P,\mathbb{C}\mathbb{P}^m)}_{\mid M}$
  (it can be easily proven in the charts $\Pi_{k,x,i}^D,\Pi_{k,x,i}^G$).

By steps 1 and 2    $j^r_D({\phi_k}_{\mid M})$ is uniformly transverse to
 $\mathbb{P}\Sigma_k\subset   \mathcal{J}^r_D(M,\mathbb{C}\mathbb{P}^m)$.

If the hypothesis on the amount of transversality over $M$ of  corollary \ref{cor:fulltrans} are met,
 then $j^r_D({\phi_k}_{\mid M})$ is uniformly transverse along $D$ to
 $\mathbb{P}\Sigma_k\subset   \mathcal{J}^r_D(M,\mathbb{C}\mathbb{P}^m)$. Observe that this requirement
  is not a problem, since the induction construction to obtain uniform transversality  over $M$ for
   $j^r\tau_k$ to $\Sigma_k^G\backslash Z_k$ can guarantee that.
\end{proof}

The vector bundles $\mathcal{J}^r_GE_k$ are endowed with hermitian metrics $\hat{g_k}$ and connections $\nabla_{k,H}$ (or just $\nabla_{H}$), which are induced by the metrics and connections on $\mathcal{J}^rE_k$ via the projection $p_G$. We do not know
 whether $\mathcal{J}^r_GE_k$ is an almost CR submanifold of $\mathcal{J}^rE_k$, but in any case we are
  not interested in doing almost complex geometry on $\mathcal{J}^r_GE_k$.

Let $\sigma_k$ be a sequence of sections of $\mathcal{J}^r_GE_k$ with  $|\nabla^j\sigma_k|_{g_k}\leq O(1),\;\forall j\geq 0$.
 Using the metric $\hat{g_k}$ we have a well defined notion of uniform transversality of $\sigma_k$ to the
  Thom-Boardman-Auroux stratification $\Sigma_k\subset \mathcal{J}^r_GE_k$ (definition \ref{def:esttrans2});
   notice that we have no notion of approximate holomorphicity neither for the sequence of sections nor for the strata.

\begin{remark}\label{rem:noah}
If $\tau_k\colon P\rightarrow E_k$ is A.H. then $|\nabla^jj^r_G\tau_k|_{g_k}\leq O(1),\;\forall j\geq 0$.
Having into account remark \ref{rem:weakstrat}, it can also be shown that  if $j^r\tau_k\colon P\rightarrow \mathcal{J}^rE_k$
 is  uniformly transverse over $M$ to $\Sigma_k^G$, then  $j^r_G\tau_k\colon P\rightarrow \mathcal{J}^r_GE_k$
  is uniformly transverse over $M$ to $\Sigma_k$.
\end{remark}

We finish this section by proving the following

  \begin{lemma}\label{lem:2calstrt}\quad
  \begin{enumerate}
  \item  Let $\mathcal{S}=(S_k^a)_{a\in A_k}$ be an approximately holomorphic  finite invariant stratification
   of $E_k$ such that in approximately holomorphic coordinates and A.H. frames each sequence of strata has a
   CR model transverse to the fibers. Let $\tau_k\colon M\rightarrow E_k$ be an A.H. sequence uniformly
   transverse along $D$ to $\mathcal{S}$. Then   $\tau_k^{-1}(\mathcal{S})$ is a stratification  of $(M,D,\omega)$
    by  2-calibrated submanifolds for all $k\gg 1$.
\item  Let $\tau_k\colon M\rightarrow E_k$ be an   A.H. uniformly transverse to $Z_k$ and whose projectivization
 $\phi_k$ is $r$-generic. Then  $B_k\cup\phi_k^{-1}(\mathbb{P}\Sigma_k)$ is a stratification by 2-calibrated submanifolds
  of $(M,D,\omega)$ for all $k\gg 1$.
\end{enumerate}

\end{lemma}
\begin{proof}
Let $S_k^a\subset E_k$. Corollary \ref{cor:transD} implies that
$\tau_k^{-1}(S_k^a)$  is uniformly transverse to $D$. Hence, if we
check that for each $x\in \tau_k^{-1}(S_k^a)$ the sequence of linear
subspaces $T_D\tau_k^{-1}(S_k^a)\subset D$ is A.H., i.e.
\[\angle_M(T_D\tau_k^{-1}(S_k^a), JT_D\tau_k^{-1}(S_k^a))\leq O(k^{-1/2})\]
 (uniformly on the point), we are done.

Let $\hat{J}$ denote the induced the almost complex structure on $E_k$.   In approximately holomorphic coordinates
 and A.H. frames, the strata $S_k\subset E_k$ have a CR model $S\subset \underline{\mathbb{C}}^m$ with respect to the canonical
 product CR structure. Recall that any almost CR structure defined out of $J_0$ in the base and the fiber, and a
 connection form with vanishing (0,1)-component, coincides with the product CR structure (this appears also in
  the proof of lemma \ref{lem:trivializationweak}). Hence the linear subspaces $T_DS=T_DS_k$ verify
   $\angle_M(T_DS,\hat{J}T_DS)\leq O(k^{-1/2})$, the bounds being uniform on the points of $\underline{\mathbb{C}}^m$,
    and hence uniform on the points of $E_k$.

The approximate holomorphicity of $\tau_k$ implies $\angle_M(T_D\tau_k,\hat{J}T_D\tau_k)\leq O(k^{-1/2})$.
Since $\angle_m(T_D\tau_k,T_DS_k)\geq \eta$, by proposition 3.7 in \cite{MPS02} for all $k \gg 1$ the intersection
 $T_D\tau_k\cap T_DS_k$ is  an A.H. sequence and  thus  also its projection to $M$, and this proves item (1).

Regarding item (2),     $B_k:={\tau_k}^{-1}(Z_k)$. Therefore item (1) applies.

The strata $\Sigma_{k,I}$  are intersected uniformly away from $B_k$. Therefore it is equivalent to work with
the projectivizations $\phi_k$ and the Thom-Boardman-Auroux stratification of $\mathcal{J}^r_D(M,\mathbb{C}\mathbb{P}^m)$,
 because $j^r_D\tau_k^{-1}(\Sigma_{k,I})=j^r_D\phi^{-1}_k(\mathbb{P}\Sigma_{k,I})$. Since for each canonical
  chart of projective space  the strata have CR models in $\mathcal{J}^r_D(M,\mathbb{C}^m)_i$,  everything reduces to item (1).
 \end{proof}

 We would like the pullback of any regular value of $\phi_k$ to be a 2-calibrated submanifold, which forces us
 to study  the behavior of an $r$-generic function near its base locus and near the pullback of the Thom-Boardman-Auroux
  strata. In our applications we would only need this analysis for the Lefschetz pencils
   $\phi_k\colon M\backslash B_k\rightarrow \mathbb{C}\mathbb{P}^1$: the same ideas used in \cite{Pr02}
   show that indeed near the base locus $|\partial \phi_k|>|\bar{\partial} \phi_k|$ and thus the regular
   ``fibers'' are 2-calibrated submanifolds. On the other hand, near the strata of the Thom-Boardman-Auroux
    stratification there is no such inequality between the holomorphic and antiholomorphic component of the
    derivative, and ad hoc modifications are needed to obtain 2-calibrated regular fibers.

In \cite{Ma05b}  the approximately holomorphic theory is appropriately modified to construct generic CR sections
 for a Levi-flat CR manifold. The complication near the base locus and degeneration loci of the
 leafwise differential does not occur (over each complex leaf the CR-Thom-Boardmann stratification is holomorphic
  and the restriction of the CR-$r$-jet holomorphic as well, therefore the former is pulled back to the leaf to
   a stratification by holomorphic strata).

\section{The main theorem}\label{sec:mainthm}

It is possible to perturb A.H. sections of  $E_k=E\otimes L^{\otimes k}\rightarrow (M,D,\omega)$ so
that their $r$-jets are transverse to an A.H. quasi-stratification of  $\mathcal{J}^r_{D}E_k$.

 \begin{theorem}\label{thm:main1}
Let  $E_k\rightarrow (M,D,\omega)$, $E_k=E\otimes L^{\otimes k}$,  and  $\mathcal{S}=(S_k^a)_{a\in A_k}$ an A.H.
 sequence of finite, Whitney (A) quasi-stratifications of $\mathcal{J}_{D}^rE_k$ transverse to the fibers. Let us
 fix $h\in \mathbb{N}$.  Let  $\delta$ be a strictly  positive constant. Then a constant $\eta>0$ exists such
  that for any A.H.  sequence    $\tau_k$ of $E_k$, it is possible to find an A.H. sequence   $\sigma_k$ of $E_k$
   so that for every  $k$ bigger than some  $k_0$,
 \begin{enumerate}
 \item
$|\nabla_{D}^j(\tau_k-\sigma_k)|_{g_k}<\delta,j=0,\dots,r+h$.
 \item $j^r_{D}\sigma_k$ is $\eta$-transverse along $D$ to $\mathcal{S}$.
 \end{enumerate}
 \end{theorem}

 Theorem \ref{thm:main2} -to be introduced- suffices for our applications; the proof of theorem \ref{thm:main1}
  -which is left to the interested reader-  is a suitable modification of the proof of theorem 1.1. in \cite{Au01},
   being the main difference the use of a result on local estimated transversality along $D_h$ to $\bf 0$  for A.H. functions
    $f_k\colon \mathbb{C}^n\times\mathbb{R}\rightarrow \mathbb{C}^m$.

Observe in theorem \ref{thm:main1} that while for any $h\in \mathbb{N}$ we can bound $|\nabla_{D}^j(\tau_k-\sigma_k)|_{g_k},\;j=0,\dots,r+h$,
by any arbitrarily small $\delta$,  we cannot do the same for the full derivative.  For the latter it can be proven that
$|\nabla^j(\tau_k-\sigma_k)|_{g_k}\leq C_j,\; \forall j\in \mathbb{N}$, where $C_j$ are constants independent of $k$ whose value we cannot control.
  Moreover the non-integrability of $D$ also forces us to work with sequences of A.H. functions all whose
  derivatives are controlled (even if we want to control the size of the perturbation along $D$ up to a finite
   order $h$); basically the derivatives along the directions of $D$ (up to some finite order $h$) will be
    arbitrarily small only if we have control for the full derivative of all the orders, and $k$ is chosen to be very large.

We can prove a strong transversality result for symplectic  manifolds with distribution $G$ along compact
 2-calibrated subvarieties.

\begin{theorem}\label{thm:main2} Let  $E_k\rightarrow (P,\Omega)$ and let $(M,D)$ be a compact $2$-calibrated
submanifold of the symplectic manifold $(P,\Omega)$ and  $G$ a $J$-complex distribution extending $D$. Let us
consider   $\mathcal{S}^G$ a $C^{h}$-A.H. sequence of finite, Whitney (A) quasi-stratifications of  $\mathcal{J}^rE_k$
 ($h\geq 2$). Let   $\delta$ be a positive constant. Then a constant  $\eta>0$ and a natural number $k_0$ exist
  such that for any   $C^{r+h}$-A.H.($C$) sequence $\tau_k$  of $E_k$, it is possible to find a  $C^{r+h}$-A.H.
   sequence   $\sigma_k$ of $E_k$ so that for any $k$ bigger than   $k_0$,
 \begin{enumerate}
 \item
$|\nabla^j(\tau_k-\sigma_k)|_{g_k}<\delta,j=0,\dots,r+h$ ($\tau_k-\sigma_k$ is $C^{r+h}$-A.H.($\delta $)).
 \item $j^r\sigma_k$ is $\eta$-transverse  over  $M$ to $\mathcal{S}^G$.
  \end{enumerate}
 \end{theorem}

 \begin{proof} We will closely follow the pattern of the proof of theorem 1.1 in \cite{Au01}, but introducing
 the appropriate modifications.

 The very basic strategy of the proof is to add a perturbation for each sequence of strata ${S^G}_k^b$, so that
  a sequence of strata is dealt with only if all the preceding ones have been already dealt with. The solution
   $\sigma_k$ will be the result of adding all the perturbations. To achieve our goal in this way we must make
   sure that at a stage corresponding to the strata ${S^G}_k^b$, the perturbation added  is such that:
 \begin{enumerate}
 \item[(i)] Uniform transversality to preceding strata is not destroyed.
 \item[(ii)] Uniform transversality to ${S^G}_k^b$ is attained.
 \end{enumerate}

 To make sure that item (i) above holds, we start by adapting the definition of local open condition
  of \cite{Au00}  to out setting:

 \begin{definition}\label{def:oppropo} Let $\eta,\bar{\eta}>0$. A family of properties $\mathcal{P}(\eta,\bar{\eta},x)_{x\in M}$
  of sections of bundles over $P$ is local and $C^q$-open, if given a section $\tau$ that verifies
   $\mathcal{P}(\eta,\bar{\eta},x)$ and a section $\sigma$ so that
 $|\tau-\sigma|_{C^q(P,g)}\leq \epsilon$, then there exist $L>0$ only depending on the $C^q$-norm of $\tau$
 so that $\tau-\sigma$ verifies $\mathcal{P}(\eta-L\epsilon,\bar{\eta}-L\epsilon,x)$.
 \end{definition}

 The advantage of a local open property is that we have an estimate on how much it varies according to the
  size of the perturbation.

 In our specific problem we say  that a $C^{r+2}$-A.H. sequence of sections  $\tau_k$ of $E_k$   verifies
  $\mathcal{P}_k(\eta,\bar{\eta},x)$, $x\in M$, if  $j^r\tau_k$ is $(\eta,\bar{\eta})$-transverse over $M$
   to ${S^G}_k^b$ at $x$. We want to show that this is a local $C^{r+2}$-open condition, because if that is
   the case we know that if at a given stage we add a perturbation with small enough  $C^{r+2}$-norm, we will still
   have a sequence of sections uniformly transverse over $M$ to ${S^G}_k^b$.

This is proven in theorem 1.1  \cite{Au01} for full transversality. For estimated transversality over $M$
the theorem is equally true because  a perturbation $\chi_k$ with  $C^{r+2}$-size bounded by $C$ gives rise to an $r$-jet
 such that (i) $|j^r\chi_k|_{g_k}\leq L'C$, (ii) $|\nabla_{TM}j^r\chi_k|_{g_k}\leq L'C$, and
  (iii)$|\nabla\nabla_{TM}j^r\chi_k|_{g_k}\leq L'C$, for some $L'>0$. Therefore small perturbations of a
   given section give rise to an $r$-jet  that remains within controlled distance of the one for the initial
    section and whose derivative along $TM$ varies in a controlled way. Similarly for a given $r$-jet  we
     can control in a ball of uniform radius its variation  up to order 2, and hence the variation of its
      derivative along $TM$ in the ball.

Next we have to make sure that the perturbation added at each stage fulfills condition (ii). We will split
the problem of achieving transversality over $M$ to ${S^G}_k^b$ into doing it for points close to the
boundary and far from the boundary. Actually, the former problem turns out to be already solved. To show
it we must check  that $(\eta_a,\bar{\eta}_a)$-transversality over $M$ of $j^r\tau_k$ to ${S^G}_k^a$, for
all $a<b$, implies the existence of $\bar{\eta}_b>0$ such that $j^r\tau_k$ is $\bar{\eta}_b$-transverse
over $M$ to ${S^G}_k^b$ at the points $\bar{\eta}_b$-close to its boundary.

In theorem 1.1  \cite{Au01} it is shown that the quasi-stratification condition together with full uniform
 transversality  can be used to show that $j^r\tau_k$ stays uniformly away from ${S^G}_k^a\backslash \Theta_{{S^G}_k^a}$,
  say at distance greater than some  $\eta'>0$; since uniform transversality over $M$ is stronger than
  uniform transversality we deduce the same result.

We now make use of  the estimated Whitney's condition (A) as in corollary \ref{cor:fulltrans}. We have the inequality
\begin{equation}\label{eq:ineq4}
\angle_\mathrm{m}(T_Mj^r\tau_k,T^{||}_M{S^G}^a_k)\leq  \angle_\mathrm{M}(T^{||}_M{S^G}_k^a,T_M{S^G}_k^b)+\angle_\mathrm{m}(T_Mj^r\tau_k,T_M{S^G}_k^b).
\end{equation}
For  $\eta''>0$ small enough  the induction hypothesis implies that for points $\eta''$-close to $\bar{\partial}{S^G}_k^b$
 there is some index $a\in A_k$ such that   \[\angle_\mathrm{m}(T_Mj^r\tau_k,T^{||}_M{S^G}^a_k)\geq \eta_a.\]
Let $\hat{M}$ denote the pullback of $TM$ to $\mathcal{J}^rE_k$. In
order to make
\[\angle_\mathrm{M}(T^{||}_M{S^G}_k^a,T_M{S^G}_k^b)<\eta_a/2\] we use
the estimated Whitney's condition (A) that gives
$\angle_\mathrm{m}(\hat{M},T{S^G}^b_k)>\gamma>0$  and
$\angle_\mathrm{M}(T^{||}{S^G}_k^a,T{S^G}_k^b)<C(\gamma)^{-1}\eta_a/2$
 (see the proof of corollary \ref{cor:fulltrans}), for $\eta''$ small enough. Then the desired result
  holds for
  \[\bar{\eta}_b:=\mathrm{min}(\eta',\eta'',\mathrm{min}_{a<b}(\eta_a/2)).\]

Therefore our task is reduced to constructing arbitrarily small
perturbations which solve the uniform transversality problem in
points $\bar{\eta}_b$-far from the boundary. We will construct such
a perturbation following Donaldson's globalization method. The key
point is the following

\begin{proposition}\label{pro:donaldsonglob} Let $\mathcal{P}_k(\eta,\bar{\eta},x)_{x\in M,\eta,\bar{\eta}>0}$
 be a family of $C^q$-open properties of sections of $E_k\rightarrow (P,g_k)$. Assume that there exist
  (uniform) constants $\rho,c',c'',p$ such that given  any $\delta>0$ small enough, any $x\in M$, and
  any sequence $\tau_k$ with uniform $C^q$-bound $C$, there exist -for all $k\gg 1$- $C^q$-bounded sections
   $\chi_{k,x}$  with the
following properties:
\begin{enumerate}
\item $|\nabla^j\chi_{k,x}|_{g_k}< c''\delta$, $j=0,\dots,q$.
\item The sections $\frac{1}{\delta}\chi_{k,x}$ have Gaussian decay away from
$x$ in $C^q$-norm.
\item $\tau_k+\chi_{k,x}$ satisfy the property $\mathcal{P}_k(\eta,\bar{\eta}-c'\delta,y)$ for all $y\in B_{g_k}(x,\rho)\cap M$,
with $\eta=c'\delta(\log (\delta^{-1}))^{-p}$.
\end{enumerate}

Then given any $\alpha>0$ and  $C^q$-bounded sections $\tau_k$ of $E_k$,
there exist -for $k\gg 1$-  $C^q$-bounded sections $\sigma_k$
of $E_k$ such that
\begin{enumerate}
\item[(i)] $|\nabla^j(\tau_k-\sigma_k)|_{g_k}< \alpha$, $j=0,\dots,q$.
\item[(ii)] The sections $\sigma_k$ satisfy $\mathcal{P}_k(\epsilon,\bar{\eta}-L\delta,x)$, for some uniform
$\epsilon,L>0$, at any $x\in M$.
\end{enumerate}
\end{proposition}

We do not give the proof of this proposition, since it is a
repetition  step by step of   Donaldson's globalization procedure
\cite{Do96}.

Hence we must check that the hypothesis of proposition
\ref{pro:donaldsonglob} hold. We will use  the following local
transversality result, which is a reformulation of lemma 5.2 and
theorem 5.4 in \cite{Moh01}.

\begin{proposition}\label{pro:mohsen} Let $F$ be a function with values in $\mathbb{C}^l$ defined over the
 ball of radius 11/10 in $\mathbb{C}^l$. Let $V$ be a vector subspace of $\mathbb{C}^l$. Let $\delta$ be
 a constant $0<\delta<1/2$. Let $\eta=\delta(P(\mathrm{log}(\delta^{-1}))^{-1}$, where $P$ is a real
 monomial depending on $n,l,V$. If in the ball of radius 11/10 we have
\[|F|_{g_0}\leq 1,\;|\bar{\partial}F|_{g_0}\leq \eta,\;|\mathrm{d}\bar{\partial}F|_{g_0}\leq \eta,\]
then there exists $u\in\mathbb{C}^p$ such that $F-u$ is $\eta$-transverse over $V$ to {\bf 0}  in the
 interior of $B(0,1)\cap V$.
\end{proposition}

We assume  that $\tau_k$ is already $\bar{\eta}_b$-transverse over $M$ at the points $\bar{\eta}_b$-close
 to the boundary. Let $0<\epsilon<\bar{\eta}_b/4$ small enough. If $x\in M$ such that
 $j^r\tau_k(x)\notin \mathcal{N}_{{S^G}_k^b}(\epsilon/2,\bar{\eta}_b)$ then $\chi_{k,x}$ is chosen to be the zero perturbation.
If  $j^r\tau_k(x)=p\in
\mathcal{N}_{{S^G}_k^b}(\epsilon/2,\bar{\eta}_b)$ then there exists
$\rho_1$ such that $j^r\tau_k(B_{g_k}(x,\rho_1))\subset
B_{\hat{g}_k}(p,\rho_\epsilon)\subset
\mathcal{N}_{{S^G}_k^b}(\epsilon,3\bar{\eta}_b/4)$.
 We consider the composition $f\circ j^r\tau_k$ pulled back to the domain of an A.H. chart adapted to $(M,G)$
  and centred at $x$. In this way we obtain a function $H_k\colon B(0,\rho_2)\subset \mathbb{C}^p\rightarrow \mathbb{C}^l$.
   If we apply  proposition \ref{pro:mohsen} directly to $H_k$, with $V=TM$, and for $\delta\ll \bar{\eta}_b/6$,
    we will obtain $\delta(P(\mathrm{log}(\delta^{-1})))^{-1}$-transversality over $M$ to {\bf 0}  for
    $H_k-u_k$ in $B_{g_k}(x,\rho_3)$. The problem is how to associate $u_k$ to a perturbation of $\tau_k$
     (the difficulty coming from the non-linearity of the strata). Instead, we consider for each index $I$
      the $\mathbb{C}^l$-valued function such that for each $y\in B_{g_k}(x,\rho_4)$
\[\Theta_I(y)=(df_1(j^r\tau_k(y))j^r\tau_{k,x,I}^{\mathrm{ref}},\dots, df_l(j^r\tau_k(y))j^r\tau_{k,x,I}^{\mathrm{ref}}),\]
with $\tau_{k,x,I}^{\mathrm{ref}}$ as defined in equation (\ref{eqn:holonframe}). There is a choice of $l$
indices $I_1,\dots,I_l$  such that the corresponding A.H. sections $j^r\tau_{k,x,I_j}^{\mathrm{ref}}$ are
 a frame for a distribution complementary to $\mathrm{Ker}df$ (and with minimal angle  bounded from below).
  Then $\Theta_{I_1},\dots,\Theta_{I_l}$ is  a frame (depending on $y$) of $\mathbb{C}^l$ comparable to the
  canonical one. We can write
\[H_k=h_k^1\Theta_{I_1}+\cdots+h_k^l\Theta_{I_l}.\]
We apply proposition \ref{pro:mohsen} (after suitable rescalings) to the $\mathbb{C}^l$-valued function
 $h_k=(h_k^1,\dots,h_k^l)$, with $V=TM$, for some $\delta$ small enough, so we get $u_k\in \mathbb{C}^l$
 such that $h_k-u_k$ is $c_1\delta(P(\mathrm{log}(\delta^{-1})))^{-1}$-transverse over $M$ to {\bf 0}
 in $B_{g_k}(x,\rho_5)$. If we multiply by the functions $\Theta_{I_1},\dots,\Theta_{I_l}$ we obtain
 $c_2\delta(P(\mathrm{log}(\delta^{-1})))^{-1}$-transversality over $M$ to {\bf 0} for
 $H_k-u_k^1\Theta_{I_1}-\cdots -u_k^l\Theta_{I_l}$.
Our perturbation is the section
\[\chi_{k,x}:=-u_k^1\tau_{k,x,I_1}^{\mathrm{ref}}-\cdots
-u_k^l\tau_{k,x,I_l}^{\mathrm{ref}}.\]

The key point is that having into account the norm of $u_k$ and the bounds on the second derivatives of $f$, the $C^1$-norm of
\[H_k-u_k^1\Theta_{I_1}-\cdots -u_k^l\Theta_{I_l}-f\circ j^r(\tau_k+s_{k,x})\]
 is bounded by $O(\delta^2)$. Since the  $C^1$-norm majorates the $C^1$-norm along $TM$
 we conclude that for $\delta$ small enough $f\circ j^r(\tau_k+s_{k,x})$ is
 $c_3\delta(P(\mathrm{log}(\delta^{-1})))^{-1}$-transverse over $M$ to {\bf 0}.
 By lemma \ref{lem:localchartrans}  we get $\mathcal{P}_k(c_4\delta(P(\mathrm{log}(\delta^{-1})))^{-1},\bar{\eta}_b-L\delta,y)$,
  for all $y\in B_{g_k}(x,\rho_5)$. Since $\mathcal{P}_k(\eta,\bar{\eta},x)$ is $C^{r+2}$-open,
  if $\delta$ is small enough compared to $\bar{\eta}_b$ and $\eta_a,\bar{\eta}_a$, we still get
  uniform transversality to the previous strata and $5\bar{\eta}_b/6$-transversality over $M$ at the points
  $3\bar{\eta}_b/4$-close to the boundary of ${S_k^G}^b$.

So we can apply proposition \ref{pro:donaldsonglob} to obtain $\mathcal{P}_k(\eta_b,3\bar{\eta}_b/4,x)$
(with respect to ${S_k^G}^b$) in all the points of $M$.

Hence we deduce the existence of  a $C^{r+2}$-A.H. sequence   $\sigma_k$ such that:
\begin{enumerate}
 \item
$|\nabla^j(\tau_k-\sigma_k)|_{g_k}<\delta,j=0,\dots,r+h$ ($\sigma_k$ is $C^{r+2}$-A.H.($\delta $)).
 \item $j^r\sigma_k$ is $\eta$-transverse  over  $M$ to $\mathcal{S}^G$.
  \end{enumerate}\end{proof}

\section{Applications}\label{sec:applications}

We begin by proving proposition \ref{pro:proposition1}, which can be also obtained as a simple corollary
of the work of  J.-P. Mohsen \cite{Moh01} together with some extra local work borrowed from  \cite{MMP03}.

\begin{proof}[Proof of proposition \ref{pro:proposition1}.]
We consider a more general situation than that of the statement of proposition  \ref{pro:proposition1}.
 Let $E$ be any rank m Hermitian vector bundle over  $(M^{2n+1},D,\omega)$, and let
 $E_k=E\otimes L_\Omega^{\otimes k}$ ($L_\Omega$ the pre-quantum line bundle of the symplectization and $E$
 is meant to be the pullback of the initial $E$ to the symplectization). We want to apply theorem
  \ref{thm:main2} to the sequence of zero sections $Z_k$, but with some changes. Basically, we want
   to start with an A.H. sequence which vanishes at $y$ and is uniformly transverse on $B_{g_k}(\rho,y)$,
   and then add perturbations not destroying these properties. We fix A.H. coordinates adapted to $(M,G)$ and
   reference sections $\tau_{k,x,j}^{\mathrm{ref}}$ centred at the points of $M\subset M\times [-\epsilon,\epsilon]$.
       In A.H. coordinates adapted to $(M,G)$ we take the sections $z_k^j\tau_{k,y,j}^{\mathrm{ref}}$, $j=1,\dots,m\leq n+1$,
        and consider their direct sum, a section of $E_k$. This sequence of sections $\tau_{k,y}$
         vanishes at $y$ and is $\eta$-transverse over $M$ to $Z_k$ in $B_{g_k}(y,\rho)$. The key point
          is to keep on adding local perturbations -as described in the proof of theorem \ref{thm:main2}-
           which vanish at $y$ and with $C^1$-norm small enough compared to $\eta$. For that we need new
            reference sections vanishing at $y$. Notice that if $d_k(x,y)\geq O(k^{1/6})$ then
            $\tau_{k,x,j}^{\mathrm{ref}}$ is already vanishing at $y$, so we do not need to change the reference
            section. Assuming $d_k(x,y)\leq O(k^{1/6})$ once we go to A.H. coordinates adapted to $(M,G)$
             and centred at $x$, the point $y$ belongs to $B(0,\rho'k^{1/6})\subset \mathbb{C}^{n+1}$. Consider
              the polynomial $P(z_k^1,\dots,z_k^{n+1})=1-z_k^1$.
Let $L_{k,y,x}\in Gl(n+1,\mathbb{C})$ be the composition of homothety and then a rotation sending $y$ to $(1,0,\dots,0)$.
We define $P_{k,y,x}=P\circ L_{k,y,x}$ and $\xi_{k,x,j}^{\mathrm{ref}}:=P_{k,y,x}\tau_{k,x,j}^{\mathrm{ref}}$.
 For any $\gamma>0$, if we suppose $d_k(x,y)\geq \gamma$ then  $\xi_{k,x,j}^{\mathrm{ref}}$ becomes an A.H.
 sequence (with bounds independent of $x$) that vanishes at $y$ and  so that $\xi_{k,x,j}^{\mathrm{ref}}$, $j=1,\dots,m$,
  fits into a local frame of $E_k$ over $B_{g_k}(x,\rho(\gamma))$ (we chose the linear map to arrange that the vanishing (affine) hyperplane of $P_{k,y,x}$ is at distance of the origin bounded from below).
Since $\tau_{k,y}$ is $\eta$-transverse over $M$ to $Z_k$ in $B_{g_k}(y,\rho)$, we only need to add
perturbations centred at points away from $B_{g_k}(y,\rho/2)$, and thus the globalization procedure can
be applied  with reference sections vanishing at $y$.

Thus it is possible to find sequences of A.H. sections $\tau_k$  of $E_k$ uniformly transverse over $M$
to $Z_k$ and vanishing at  $y$. Hence ${\tau_k}_{\mid M}$ are uniformly transverse to $Z_k$  and vanishing
 at $y$. Let $W_k={\tau_k}_{\mid M}^{-1}(Z_k)$.  For all $k \gg 1$ by corollary \ref{cor:transD} they are
 uniformly transverse to $D$, and by  lemma \ref{lem:2calstrt} approximately almost complex and therefore 2-calibrated.

The study of its topology is done very much as in the symplectic and contact cases (see the proofs in  \cite{Do96,Au97,IMP01}).
\end{proof}

The next result we want to prove is the existence of determinantal
submanifolds (proposition \ref{pro:determinantal}), which is still a
transversality result for  $0$-jets (vector bundles  $E_k$), but not
anymore to the  {\bf 0} section but to a  sequence of non-linear
approximately holomorphic stratifications.

\begin{proof}[Proof of proposition \ref{pro:determinantal}.]
Let $E,F\rightarrow M$ be Hermitian bundles with connection and let
us define the sequence of very ample  vector bundles
$I_k:=E^*\otimes F\otimes L^{\otimes k}$. In the total space of
$I_k$ we consider the sequence of stratifications $S_k$ whose strata
are
 $S_{k,i}=\{A\in I_k |\; \mathrm{rank}(A)=i\}$, where $A\in \mathrm{Hom}(E,F\otimes L^{\otimes k})$.

Let $E,F$ still denote the pullback of $E,F$ to the symplectization $(M\times [-\epsilon,\epsilon],\Omega)$.
Let  $I_{k,\Omega}\rightarrow M\times [-\epsilon,\epsilon]$ be
$E^*\otimes F\otimes L^{\otimes k}_\Omega=\mathrm{Hom}(E,F\otimes
L^{\otimes k}_\Omega)$. Let $G$ be as usual a $J$-complex distribution
defined on  $M\times [-\epsilon,\epsilon]$ that extends $D$, and let
\[S_{k,i}^G=\{A\in I_{k,\Omega} |\; \mathrm{rank}(A)=i\},\;\;A\in
\mathrm{Hom}(E,F\otimes L^{\otimes k}_\Omega).\]
By lemma \ref{lem:trivializationweak} (applied to almost complex
manifolds)  $S_{k,i}^G$ is an  approximately holomorphic  sequence
of finite, Whitney (A) stratifications. Therefore we can apply theorem
\ref{thm:main2} to construct an A.H. sequence of sections $\tau_k$
of $I_{k,\Omega}$ uniformly transverse  over $M$ to $S_k^G$, and
thus along $D$.

Hence   for all $k$ large enough $M$ is stratified by the
submanifolds $S_{i}(\tau_k)=\{x\in M|\;
\mathrm{rank}(\tau_k(x))=i\}$, which are uniformly transverse  to
$D$ and 2-calibrated by lemma \ref{lem:2calstrt}.
 \end{proof}

Corollary \ref{cor:determinantal2} follows from the fact that in the
contact case the 2-form  is exact and hence the cohomological
computations are those of the bundle $E^*\otimes F$.

\begin{theorem}\label{thm:genericity}
Let  $(M,D,\omega)$ be a closed integral 2-calibrated manifold,
set $E_k=\underline{\mathbb{C}}^{m+1}\otimes L^{\otimes k}$, and let
$r$ be any natural number. Any A.H. sequence of sections of
$\underline{\mathbb{C}}^{m+1}\otimes L^{\otimes k}_\Omega\rightarrow
(M\times [-\epsilon,\epsilon],\Omega,G)$ admits an arbitrarily small
$C^{r+h}$-perturbation  such that ${\phi_k}_{\mid M}\colon
M\backslash B_k\rightarrow \mathbb{C}\mathbb{P}^m$ -the restriction
to $M$ of its projectivization- is an $r$-generic A.H. sequence.
\end{theorem}

\begin{proof} The proof is just theorem \ref{thm:main2} applied to the Thom-Boardman-Auroux
 quasi-stratification along $G$ of $\mathcal{J}^rE_k\rightarrow (M\times [-\epsilon,\epsilon],J,G,g_k)$,
  combined with proposition \ref{pro:relgen}.
\end{proof}

It must be pointed  out that the behavior of A.H. functions at the points close to the   degeneration loci  is more complicated than
that of the leafwise holomorphic model: firstly, and similarly to
what happens for even dimensional  almost complex manifolds, to
obtain normal forms it is necessary to add perturbations so that the
function becomes holomorphic (at least in certain directions);
otherwise the approximate holomorphicity is not significative due to
the vanishing (degeneracy) of the holomorphic part. Secondly, we
have an extra non-holomorphic direction that we do not control.  At
most, we can apply the usual genericity results to that direction
(the perturbations at most of size $O(k^{-1/2})$ so as not to
destroy the other properties).

One instance of the preceding theorem is when the target space has large dimension,
 so that the generic map is an immersion along the directions of $D$.

\begin{proof}[Proof of corollary \ref{cor:embedding}.]
Set $E_k=\underline{\mathbb{C}}^{m+1}\otimes L^{\otimes k}$, where $m\geq 2n$.
Theorem \ref{thm:main2} is applied to the Thom-Boardman-Auroux
quasi-stratification  along $G$ of   $\mathcal{J}^1E_k\rightarrow
(M\times[-\epsilon,\epsilon],J,G,g_k)$), to obtain 1-generic A.H.
maps  ${\phi_k}\colon M\rightarrow \mathbb{C}\mathbb{P}^m$. From the
choice of  $m$ it follows that the set of base points and of points
where  $\partial{\phi_k}$ is not injective is empty. It is clear
that by construction that $\phi_k^*[\omega_{FS}]=[\omega_k]$.
\end{proof}

This is a non-trivial result because the property of being an
immersion along  $D$ is not generic (for smooth maps to
$\mathbb{C}\mathbb{P}^{2n}$). Notice that if for example $D$ is
integrable the property is generic for each leaf (locally), but not
for the 1-parameter family.

As mentioned in the introduction, the previous corollary can be improved in two different ways.

\begin{proof}[Proof of corollary \ref{cor:Tischler}.]
Let us assume that any $2$-form in the path
$\rho_{k,t}=(1-t)\omega_k+t\phi_k^*\omega_{FS}$  is non-degenerate
over $\mathcal{D}$, where $\omega_{FS}$ is be the Fubini-Study
2-form. Then  Moser's trick can be applied leafwise: if $\alpha$
is a 1-form such that $d\alpha=-(\phi_k^*\omega_{FS}-\omega_k)$, the
vector fields tangent to $\mathcal{D}$ defined by the condition
$-i_{X_t}\rho_{k,t}=-\alpha$ generate a 1-parameter family of
diffeomorphisms preserving each leaf and sending $\rho_{k,t}$ to
$\omega_k$.

The non-degeneracy over $\mathcal{D}$ of $\rho_t$ follows from the
estimated  transversality of $\phi_k$ together with the approximate
holomorphicity. For any $v\in D_x$ of $g_k$-norm 1,

\[\rho_{k,t}(v,Jv)=(1-t)\omega_k(v,Jv)+t\omega_{FS}(\phi_{k*}v,\phi_{k*}Jv)\geq (1-t)+t\eta>0.  \]
\end{proof}

In general a closed Poisson manifold with codimension one leaves
does not admit a  lift to a 2-calibrated structure (for example any
non-taut smooth foliation in $M^3$). The previous corollary can be
used to state the following result:

\begin{corollary} Let $(M^{2n+1},\mathcal{D},\omega_{\mathcal{D}})$ be a closed
 Poisson manifold with co-oriented codimension one leaves. Then the Poisson
  structure admits a lift to a (rational) 2-calibrated structure, if and only
  if a multiple of $\omega_{\mathcal{D}}$ is induced by a leafwise immersion
  in $\mathbb{C}\mathbb{P}^{2n}$ (by pulling back $\omega_{FS}$).
\end{corollary}

It is worth mentioning that it is possible to obtain uniform
transversality to a finite number of quasi-stratifications of the
same sequences of bundles. For example, and this leads to the second
improvement of corollary \ref{cor:embedding}, we can obtain the
1-genericity result that gives rise to embeddings in
$\mathbb{C}\mathbb{P}^m$ transverse to a finite number of complex
submanifolds of $\mathbb{C}\mathbb{P}^m$.

 We just need to consider for each submanifold the sequence of stratifications
 $\mathbb{P}\mathcal{S}$ of  $\mathcal{J}^1_G(M,\mathbb{C}\mathbb{P}^m)$, whose
 unique stratum (for each $k$) is defined to be the  1-jets along $G$ whose degree
  $0$ component is a point of the submanifold; next we pull it back to a stratification
   $\mathcal{S}$ of  $\mathcal{J}^1_GE_k^*$ and finally to a stratification
    $\mathcal{S}^{G}$ of $\mathcal{J}^1E_k^*$ (the structure near  $Z_k$ is not
    relevant because transversality to the Thom-Boardman-Auroux quasi-stratification
    along $G$ implies that the sections stay away from $Z_k$). Therefore, we have defined
     a stratification  of $\mathcal{J}^1E_k$ which is trivially approximately holomorphic
      because it is the pullback by A.H. maps of an initial approximately holomorphic stratification
       of $\mathcal{J}^0_G(M,\mathbb{C}\mathbb{P}^m)$.
Any 1-generic sequence of A.H. sections of  $E_k$ uniformly
transverse to  $\mathcal{S}^{G}$, when  restricted to $M$ gives rise
to maps $\phi_k\colon M\hookrightarrow  \mathbb{C}\mathbb{P}^m$
uniformly transverse along $D$ to the submanifold.

\begin{proof}[Proof of Theorem \ref{thm:pencil}] We first apply theorem \ref{thm:genericity}
 to obtain  ${\phi_k}_{\mid M}\colon M\backslash B_k\rightarrow \mathbb{C}\mathbb{P}^1$ 1-generic.

Near the base points and the points where $\nabla_D{\phi_k}_{\mid M}$ vanishes, we apply the
 perturbations defined in \cite{Pr02} to obtain the required local models.
\end{proof}

Another possible application is, as proposed by  D. Auroux for
symplectic manifolds  \cite{Au00,Au01}, to obtain $r$-generic
applications to $\mathbb{C}\mathbb{P}^m$ whose composition with
certain projections  $\mathbb{C}\mathbb{P}^m\rightarrow \mathbb{C}
\mathbb{P}^{m-h}$ are still $r$-generic (the corresponding
stratifications are approximately holomorphic because they are
pullback of approximately holomorphic stratifications by A.H. maps;
the structure near $Z_k$ is also seen to be appropriate).

It is also possible to develop an analogous construction but for
A.H. maps to Grassmannians  $Gr(r,m)$, starting from sections of
$\underline{\mathbb{C}}^{r}\otimes E_k$, $E_k$ of rank $m$ (see
\cite{MPS02,Au02}).

Our techniques can be applied to any closed 2-calibrated manifold to
give a  finer topological description of the 2-calibrated structure.
It is possible to apply the same idea to manifolds for which the
2-calibrated structure enters as an auxiliary tool. This point of
view has already been adopted in \cite{MMP03}.

We recall the following result.

\begin{theorem} (Gromov) Let $M^{2n+1}$ be a closed manifold whose structural group reduces
 to $U(n)$, and let $a\in H^2(M;\mathbb{Z})$. Then there exists $\omega$ a closed maximally non-degenerate 2-form such that $[\omega]=a$.
\end{theorem}
\begin{proof} The structural group of the open manifold $M\times \mathbb{R}$ reduces to $U(n+1)$. Then by \cite{Gr86} it carries a symplectic form representing any given cohomology class, in particular the pullback of $a$ to $M\times \mathbb{R}$. Its restriction to $M\times \{0\}$ is $\omega$.
\end{proof}

So by selecting any codimension one distribution transverse to the kernel of $\omega$, we have:
\begin{corollary}\label{cor:upgrade} Let $M^{2n+1}$ be a closed manifold whose structural group
 reduces to $U(n)$, and let $a\in H^2(M;\mathbb{Z})$. Then $M$ admits 2-calibrated
  structures $(D,\omega)$ for which $[\omega]=a$.
\end{corollary}

Notice that if we apply any of the previous constructions to $(M,D,w)$, we obtain
submanifolds and more generally stratifications of $M$ by 2-calibrated submanifolds.
 Regarding the initial structure, which was just a reduction of the structural group
 to $U(n)$, we can conclude that the corresponding strata also admit such a reduction.

\renewcommand{\thesection}{Appendix \Alph{section}}
\addtocounter{section}{-8}

\section[\hspace{1.5 cm} Proof of proposition \ref{pro:perthol}]{Proof of proposition \ref{pro:perthol}}\label{sec:newconn}

\renewcommand{\thesection}{\Alph{section}}

We write down the proof for the bundle $\mathcal{J}^rE_k$ because it is
a necessary  ingredient in the proof of theorem \ref{thm:main2}. The
case of $\mathcal{J}^r_DE_k$ bears no further complications and it
is left to the interested reader.

We omit the subindices  $k$ and $r$ for the connections whenever there
is no risk of confusion.

Recall that in coordinates the curvature can be computed as follows:
in a chart where $T^*P$ is trivialized using the derivatives of the
coordinates, we have the corresponding flat connection $\textrm{d}$
on $T^*P$. We have the operator

\begin{eqnarray*}
\nabla^1\colon T^*P\otimes E_k & \longrightarrow & T^*P\otimes T^*P\otimes E_k\\
\nabla^1&:=&\mathrm{d}\otimes \mathrm{I}-\mathrm{I}\otimes \nabla
\end{eqnarray*}
and the antisymmetrization map

\begin{eqnarray*}
\mathrm{asym}_2\colon T^*P\otimes T^*P & \longrightarrow & \wedge^2T^*P\\
 \alpha\otimes\beta &\longmapsto & \alpha\wedge\beta,\\
 \alpha\wedge\beta(u,v):=\alpha(u)\beta(v)-\alpha(v)\beta(u).&&
 \end{eqnarray*}
The curvature is the composition $\mathrm{asym}_2(\nabla^1\circ \nabla)$.

Let $\sigma_k=(\sigma_{k,0},\sigma_{k,1})$ be a section (maybe
local) of $\mathcal{J}^1 E_k$.  The modified connection is
$\nabla_{H_1}(\sigma_{k,0},\sigma_{k,1})= (\nabla\sigma_{k,0},\nabla
\sigma_{k,1})+ (0,-F^{1,1} \sigma_{k,0})$, where $-F^{1,1}
\sigma_{k,0}\in T^{*0,1}P\otimes T^{*1,0}P\otimes E_k$ (see
\cite{Au02}). For jets along $D$ we add $-F^{1,1}_D$.

The previous formula defines a connection.

\begin{lemma}\label{lem:holcurv} Let $\underline{\mathbb{C}}^m\rightarrow \mathbb{C}^p$
be the trivial bundle endowed with a connection $\nabla$ whose
curvature is of type (1,1) with respect to the canonical complex
structure $J_0$; the connection splits  into
$\partial_{\nabla}+\bar{\partial}_{\nabla}$.  Let $\tau$ be a
holomorphic section of $\underline{\mathbb{C}}^m$ (with respect to to the
holomorphic structure induced by $\nabla$). Then

\begin{equation}\label{eq:holmod}
\nabla_H(\tau,\partial_{\nabla}\tau)=\nabla(\tau,\partial_{\nabla}\tau)-(0,\bar{\partial}_{\nabla}\partial_{\nabla}\tau)
\end{equation}
and ${\bar{\partial}}_{\nabla_{H}}(\tau,\partial_{\nabla}\tau)=0$.
\end{lemma}
\begin{proof}
By definition

\begin{equation*}\label{eq:comp1}
F\tau=\mathrm{asym}_2(\nabla^1\nabla \tau).
\end{equation*}
Let us denote the trivialization of the bundle that identifies it with $\underline{\mathbb{C}}^m$  by $\xi_1,\dots,\xi_m$. Since $\tau$ is holomorphic

\begin{equation*}\label{eq:comp2}
F\tau=\mathrm{asym}_2((\mathrm{d}\otimes \mathrm{I}-\mathrm{I}\otimes \nabla) \partial_\nabla \tau).
\end{equation*}
If we write  $\partial_\nabla\tau=dz^i h_i^j\xi_j$ then

\begin{equation*}\label{eq:comp3}
F\tau=\mathrm{asym}_2(-(\mathrm{I}\otimes \nabla)dz^i h_i^j\xi_j).
\end{equation*}
But being the curvature of type (1,1) we can write

\begin{equation}\label{eq:comp3a}
F\tau=\mathrm{asym}_2(-(\mathrm{I}\otimes \bar{\partial}_\nabla)dz^i h_i^j\xi_j).
\end{equation}

Recall that $F\tau$ has to be understood as an element of
$T^{*0,1}\mathbb{C}^p\otimes T^{*1,0}\mathbb{C}^p\otimes
\underline{\mathbb{C}}^m$. That amounts to switch the $d\bar{z}^l$'s
with the $dz^l$'s, which cancels the negative sign on the right hand side of
equation (\ref{eq:comp3a}). Thus what we obtain is:

\begin{equation}\label{eq:comp3b}
F\tau:=(\mathrm{I}\otimes \bar{\partial}_\nabla)dz^i h_i^j\xi_j\in\Gamma(T^{*0,1}\mathbb{C}^p\otimes T^{*1,0}\mathbb{C}^p\otimes \underline{\mathbb{C}}^m).
\end{equation}
But equation (\ref{eq:comp3b}) equals
\[(\bar{\partial}_0\otimes \mathrm{I}+\mathrm{I}\otimes \bar{\partial}_\nabla)dz^i h_i^j\xi_j\]
which by definition is
\begin{equation}\label{eq:comp4}
\bar{\partial}_\nabla\partial_\nabla\tau.
\end{equation}

By equation (\ref{eq:comp4})

\[{\bar{\partial}}_{\nabla_{H}}(\tau,\partial_{\nabla}\tau)=
(\bar{\partial}_\nabla\tau,\bar{\partial}_\nabla\partial_{\nabla}\tau-\bar{\partial}_\nabla\partial_\nabla\tau)=0.\]
\end{proof}

It is also clear that  $\partial_\nabla=\partial_{\nabla_H}$ and therefore they define the same coupled holomorphic jets.

Lemma \ref{lem:holcurv} has an obvious approximately holomorphic
version:  if we have a very ample sequence of rank $m$ vector bundles
by definition the sequences of  curvatures is approximately of type
(1,1). Then we can fix approximately holomorphic coordinates and the
first part of lemma  \ref{lem:holcurv} implies that for $\tau_k$ a
sequence of A.H. sections of $E_k$,  one has

\[F\tau_k\approxeq \bar{\partial}\partial\tau_k,\]
and by the second part

\[\bar{\partial}_Hj^1\tau_k\approxeq 0.\]

We now move into computing the curvature of the modified connection in the
integrable case.  We will denote the coupled holomorphic $r$-jet in
the integrable model by $j^r_{\mathrm{hol}}\tau$.

\begin{lemma} \label{lem:holcurv2} Let $\underline{\mathbb{C}}^m\rightarrow \mathbb{C}^p$  be
 the trivial bundle as in lemma  \ref{lem:holcurv}. Assume also that for the fixed
   trivialization $\xi_1,\dots,\xi_m$ the curvature is a matrix with constant
    coefficients and that we have a frame given by holomorphic sections $\tau_1,\dots,\tau_m$.
Then $F_\nabla=F_{\nabla_H}$.
\end{lemma}
\begin{proof}
If the holomorphic sections $\tau_1,\dots,\tau_m$ generate the
bundle,  then the holomorphic 1-jets of $z^l\tau_j$, $\tau_j$,
$1\leq l\leq p$, $1\leq j\leq m$ are a basis of
$\mathcal{J}^1_{p,m}$ (at least on $B(0,\rho)$). By lemma
\ref{lem:holcurv}, they are a holomorphic basis.

\begin{equation}\label{eq:holcurv4}
\nabla_Hj^1_{\mathrm{hol}}z^l\tau_j=(\partial_\nabla(z^l\tau_j), \nabla\partial_\nabla(z^l\tau_j))-(0,Fz^l\tau_j)=\nabla j^1_{\mathrm{hol}}z^l\tau_j-(0,Fz^l\tau_j).
\end{equation}

Let us write again $\partial_\nabla\tau_j=dz^i h_{i,j}^s\xi_s$,  and
$F=a_{ts}d\bar{z}^t dz^s\in \Gamma(T^{*0,1}\mathbb{C}^p\otimes
T^{*1,0}\mathbb{C}^p)$. If we apply to $\nabla
j^1_{\mathrm{hol}}z^l\tau_j$ the operator $\mathrm{asym}_2
\nabla_H^1$, $\nabla^1_H:=\textrm{d}\otimes
\mathrm{I}-\mathrm{I}\otimes \nabla_H$,  we get:

\begin{equation}\label{eq:holcurv5}
F_\nabla j^1_{\mathrm{hol}}z^l\tau_j+(0,\mathrm{asym}_2(dz^l a_{ts}d\bar{z}^t dz^s\tau_j+z^ldz^ia_{ts}d\bar{z}^t dz^s h_{i,j}^s\xi_s)).
\end{equation}
When we apply the same operator to  $(0,F z^l\tau_j)$, if  recall
that the  $a_{ts}$ are constant  and that $z^l\tau_j$ is a
holomorphic section we get

\begin{equation}\label{eq:holcurv6}
\mathrm{asym}_2 \nabla_H^1(0,F z^l\tau_j)=(0,\mathrm{asym}_2(-a_{ts}d\bar{z}^t dz^l dz^s\tau_j-a_{ts}d\bar{z}^t z^ldz^i h_{i,j}^s\xi_s)),
\end{equation}
and the right hand side of equation (\ref{eq:holcurv6}) equals
\begin{equation}\label{eq:holcurv7}
(0,\mathrm{asym}_2(dz^l F\tau_j+z^ldz^iF h_{i,j}^s\xi_s)).
\end{equation}
If we put together equations (\ref{eq:holcurv4}), (\ref{eq:holcurv5}), and (\ref{eq:holcurv7}) we obtain
\[F_{\nabla_H}\tau_j=F_\nabla\tau_j.\]
\end{proof}

We want to use a recursive construction based on lemmas
\ref{lem:holcurv}  and \ref{lem:holcurv2} to introduce the desired
connection on $\mathcal{J}^r_{p,m}$.

Before doing that we recall that the coupled holomorphic jets are sections of $\mathcal{J}^r_{p,m}$.

We now prove how to modify the connection on $\mathcal{J}^2_{p,m}$.

\noindent\emph{Step 1:} We identify $\mathcal{J}^2_{p,m}$ with the
subbundle  of $\mathcal{J}^1\mathcal{J}^1_{p,m}$ spanned by
holonomic sections, i.e. sections of the form
$j^1_{\mathrm{hol}}j^1_{\mathrm{hol}}\tau$, where $\tau$ is any
holomorphic section of $\underline{\mathbb{C}}^m$.
Pointwise, an element $\gamma$ of the fiber of $\mathcal{J}^1\mathcal{J}^1_{p,m}$ is of the form

\[(\gamma_{0,0},\gamma_{0,1},\gamma_{1,0},\gamma_{1,1})\in (\mathbb{C}\oplus  T^{*1,0}\mathbb{C}^p \oplus T^{*1,0}\mathbb{C}^p\oplus (T^{*1,0}\mathbb{C}^p\otimes T^{*1,0}\mathbb{C}^p) )\otimes \mathbb{C}^m,\]
and belongs to $\mathcal{J}^2_{p,m}$ if and only if
$\gamma_{1,1}\in T^{*1,0}\mathbb{C}^p\odot
T^{*1,0}\mathbb{C}^p\otimes \mathbb{C}^m$, and
$\gamma_{1,0}=\gamma_{0,1}$.

Using the metric induced by the Euclidean one on the base and fiber
and the connection,  we have a orthogonal projection $r\colon
\mathcal{J}^1\mathcal{J}^1_{p,m}\rightarrow\mathcal{J}^2_{p,m}$.

\noindent\emph{Step 2:} We introduce a new connection on  $\mathcal{J}^1\mathcal{J}^1_{p,m}$.

 On  $\mathcal{J}^1_{p,m}$ we use the modified connection $\nabla_{H_1}$. This, together with
  the flat connection $\textrm{d}$ on $T^*\mathbb{C}^p$ defines a connection  $\nabla_{H_{1,1}}$
   on $\mathcal{J}^1\mathcal{J}^1_{p,m}$. Notice that on $\mathcal{J}^1\mathcal{J}^1_{p,m}$
   we also have a connection $\nabla_2$ coming from $\textrm{d}$ and $\nabla_1$.

We consider the trivialization of  $\mathcal{J}^1_{p,m}$ furnished
by the sections  $\xi_j$, $dz^i\xi_j$, $1\leq j\leq m$, $1\leq i\leq
p$, so we can identify the bundle with
$\underline{\mathbb{C}}^{mp+m}$. This is a trivial bundle with
connection $\nabla_{H_1}$. By lemma \ref{lem:holcurv2}
$F_{\nabla_{H_1}}=F_{\nabla_1}$. Recall also that in the basis
$\xi_j$, $dz^i\xi_j$ the curvature $F_{\nabla_1}$ is a matrix that
decomposes into $p+1$ blocks corresponding to $\xi_1,\dots,\xi_m$ and
to $dz^i\xi_1,\dots,dz^i\xi_m$, $1\leq i\leq p$. For each such block
the corresponding matrix is the one for $F_\nabla$ in the basis
$\xi_j$. Therefore $F_{\nabla_{H_1}}$ is still of type (1,1) and has
constant entries in the aforementioned basis.

 Let $\nabla_{H_2}$ be the result of modifying $\nabla_{H_{1,1}}$.
Since $\nabla_{H_{1}}$ is of type (1,1) by lemma \ref{lem:holcurv}
applied  to $(\underline{\mathbb{C}}^{mp+m},\nabla_{H_1})$, if
$\tau^1\in\Gamma(\mathcal{J}^1_{p,m})$ is holomorphic with respect to
$\nabla_{H_1}$, then $j^1_{\mathrm{hol}}\tau^1$ is holomorphic
with respect to $\nabla_{H_2}$. In particular
$j^1_{\mathrm{hol}}(j^1_{\mathrm{hol}}z^i\tau_j)$,
$j^1_{\mathrm{hol}}(z^lj^1_{\mathrm{hol}}z^i\tau_j)$ are a local
holomorphic frame of
$(\mathcal{J}^1\mathcal{J}^1_{p,m},\nabla_{H_2})$ (recall that
$\tau_j$ was a local holomorphic frame of
$\underline{\mathbb{C}}^m$).

Having into account that the curvature of
$(\underline{\mathbb{C}}^{mp+m}$, $\nabla_{H_1})$ is of type  (1,1)
and with constant entries, and that
$(\underline{\mathbb{C}}^{mp+m}$, $\nabla_{H_1})$ has a local
holomorphic basis, lemma \ref{lem:holcurv2} gives
$F_{\nabla_{H_2}}=F_{\nabla_{H_{1,1}}}$. From
$F_{\nabla_{H_1}}=F_{\nabla_1}$ it follows that
$F_{\nabla_{H_{1,1}}}=F_{\nabla_2}$. Therefore

\begin{equation*}\label{eq:curv2jet}
F_{\nabla_{H_2}}=F_{\nabla_2}\, \textrm{on} \, \mathcal{J}^1\mathcal{J}^1_{p,m}.
\end{equation*}

\noindent\emph{Step 3:} Check that $\nabla_{H_2}$ restricts to
$\mathcal{J}^2_{p,m}\hookrightarrow
\mathcal{J}^1\mathcal{J}^1_{p,m}$ with the desired properties.

Let $I=(i_0,i_1,\dots,i_p)$ with $1\leq i_0\leq m$, $0\leq i_j\leq
2$, $i_1+\cdots +i_p\leq 2$,  and let $\tau_I:=z_1^{i_1}\dots
z_p^{i_p}\tau_{i_0}$. We consider the sections
$j^1_{\mathrm{hol}}j^1_{\mathrm{hol}}\tau_I$, which are a local
holomorphic frame $\mathcal{J}^2_{p,m}$ (using the identification
described in step 1).  We will see that  $\nabla_{H_2}
j^1_{\mathrm{hol}}j^1_{\mathrm{hol}}\tau_I\in
\Gamma(T^{*1,0}\mathbb{C}^p\otimes \mathcal{J}^2_{p,m})$, and
therefore that the connection $\nabla_{H_2}$ preserves
$\mathcal{J}^2_{p,m}$.

We just proved in step 2 that
$j^1_{\mathrm{hol}}j^1_{\mathrm{hol}}\tau_I$ is holomorphic  with respect to
$\nabla_{H_2}$ and that
$\partial_{\nabla_{H_2}}=\partial_{\nabla_{H_{1,1}}}=\partial_{\nabla_2}$.
Let us write
$j^1_{\mathrm{hol}}j^1_{\mathrm{hol}}\tau_I=(\tau_I,\partial_\nabla\tau_I,\partial_\nabla\tau_I,\partial^2_\nabla\tau_I)$.
Then
\begin{align*}
\nabla_{H_2}j^1_{\mathrm{hol}}j^1_{\mathrm{hol}}\tau_I &=\partial_{\nabla_{H_2}}j^1_{\mathrm{hol}}j^1_{\mathrm{hol}}\tau_I=\\
\partial_{\nabla_{2}}(\tau_I,\partial_\nabla\tau_I,\partial_\nabla\tau_I,\partial^2_\nabla\tau_I)&=(\partial_\nabla\tau_I,\partial_\nabla\partial_\nabla\tau_I,\partial_\nabla\partial_\nabla\tau_I,\partial_\nabla\partial^2_\nabla\tau_I),
\end{align*}
which belongs to $\Gamma(T^{*1,0}\mathbb{C}^p\otimes\mathcal{J}^2_{p,m})$.

Therefore, the curvature of the restriction of $\nabla_{H_2}$ to
$\mathcal{J}^2_{p,m}$ is of  course of type $(1,1)$. The last
observation is its expression in a suitable basis. The curvature of
$\nabla_2$ on $\mathcal{J}^1\mathcal{J}^1_{p,m}$ splits on blocks
corresponding to the basis $\xi_1,\dots,\xi_m$,
$dz^i\xi_1,\dots,dz^i\xi_m$, $dz^l\xi_1,\dots,dz^l\xi_m$,
$dz^r\otimes dz^t\xi_1,\dots,dz^r\otimes dz^t\xi_m$, $1\leq
i,l,r,t\leq p$. Each submatrix is $F_\nabla $. If we use the basis
$\xi_1,\dots,\xi_m$, $dz^i\xi_1,\dots,dz^i\xi_m$, $dz^r\odot
dz^t\xi_1,\dots,dz^r\odot dz^t\xi_m$, $1\leq i,r,t\leq p$, the
curvature equally splits into blocks each matching $F_\nabla $.

The general case uses the following induction step: on
$\mathcal{J}^r_{p,m}$ there exists a  connection $\nabla_{H_r}$ with
the following properties:
\begin{enumerate}
\item $\partial_{H_r}=\partial_r$.
\item $F_{\nabla_{H_r}}=F_{\nabla_r}$ and therefore $F_{\nabla_{H_r}}$ is of type $(1,1)$.
\item If $\bar{\partial}_\nabla\tau=0$ then $\bar{\partial}_{H_r}j^r_{\mathrm{hol}}\tau=0$.
\item In the basis $\xi_I:={(dz_k^1)}^{\odot i_1}\cdots {(dz_k^n)}^{\odot i_n}\xi_{i_0}$ the
 curvature splits into blocks each matching $F_\nabla $.
\end{enumerate}

To define $\nabla_{H_{r+1}}$ on $\mathcal{J}^{r+1}_{p,m}$ we reproduce the previous 3 steps.

Firstly we consider the identification of  $\mathcal{J}^{r+1}_{p,m}$
with the subbundle of  $\mathcal{J}^1\mathcal{J}^{r}_{p,m}$ spanned
by sections of the form $j^1_{\mathrm{hol}}j^r_{\mathrm{hol}}\tau$,
$\tau$ a holomorphic section of $\underline{\mathbb{C}}^m$.

Secondly we consider the connection $\nabla_{H_{1,r}}$ on
$\mathcal{J}^{r+1}_{p,m}$ constructed  out of $\textrm{d}$ and
$\nabla_{H_r}$ and modify it to $\nabla_{H_{r+1}}$. By the induction hypothesis using the basis  $\xi_I$ we are in the
situation of lemma \ref{lem:holcurv2},  for $\mathcal{J}^{r}_{p,m}$
identifies with $\underline{\mathbb{C}}^{N_r}$ with a connection
whose curvature is of type (1,1) and with constant coefficients, and
with a frame of holomorphic sections. Therefore
$F_{\nabla_{H_{r+1}}}=F_{\nabla_{H_{1,r}}}=F_{\nabla_{r+1}}$.
Since we can also apply lemma  \ref{lem:holcurv},  for any
$\tau^r\in \Gamma(\mathcal{J}^{r}_{p,m})$  the 1-jet
$j^1_{\mathrm{hol}}\tau^r$ is holomorphic with respect to $\nabla_{H_{r+1}}$.

The third step is to check that the modified connection restricts to
$\mathcal{J}^{r+1}_{p,m}\hookrightarrow
\mathcal{J}^1\mathcal{J}^{r}_{p,m}$. Using that
$\partial_{\nabla_{H_{r+1}}}=\partial_{\nabla_{r+1}}$, any frame of
sections of the form $j^1_{\mathrm{hol}}j^r_{\mathrm{hol}}\tau_I$,
$\tau_I$ holomorphic, is sent by the connection to sections of
$\mathcal{J}^{r+1}_{p,m}$.

It is also routine to check that in the basis $\xi_I$ the curvature matrix is made of blocks of the form $F_\nabla$.

The almost complex counterpart of the result we just proved is done
exactly in the same way.  The only modification is that the
connection on $\mathcal{J}^1\mathcal{J}^{r}E_k$ does not descend
automatically to a connection on
$\mathcal{J}^{r+1}E_k\hookrightarrow
\mathcal{J}^1\mathcal{J}^{r}E_k$. We have to project via $r\colon
\mathcal{J}^1\mathcal{J}^rE_k\rightarrow\mathcal{J}^{r+1}E_k$, but
this is seen to introduce an error which is approximately vanishing.
It might happen that the resulting connection amounts to adding also
a pseudo-holomorphic part. If that is the case we forget about this
contribution (which again would be approximately vanishing).
Therefore, we obtain a connection with all the desired properties.

 Using similar considerations to the ones for $1$-jets, it can be deduced  that the   (r+1)-jet
  of a  $C^{r+1+h}$-A.H. sequence of sections of  $E_k$ is a $C^h$-A.H. sequence of sections of
  $(\mathcal{J}^{r+1} E_k,\nabla_{H_{r+1}})$.

\renewcommand{\thesection}{Appendix \Alph{section}}

\section[\hspace{1.5 cm} Chern classes and top Chern classes]{Chern classes and top Chern classes}\label{sec:det}

Corollary \ref{cor:determinantal2} proves the existence of contact determinantal submanifolds, which we expect to be more general than  those coming from zeroes of vector bundles constructed in \cite{IMP01}. To support this  we recall that it is known that in the algebro geometric setting that determinantal varieties are more general that zeroes
    of vector bundles (see for example \cite{Ha74,Ar06}), and a similar result
    should be  expected to hold in the smooth category. A way to prove it would
     be exhibiting  a manifold in which there exist a cohomology class $a$ which
      is the Chern class of a complex vector bundle $F$ but it is not the top
       Chern class of any complex vector bundle (i.e. showing Chern classes are more
       general than top Chern classes), the reason being that if we choose as
        $E$ the trivial complex vector bundle of the appropriate rank and the appropriate
         determinantal locus, we have
\[\Delta_{E,F,i}=a.\]

As far as the author knows such a question has not been addressed.
 A lot is known about cohomology classes which can be Chern classes,
 mainly because for a given finite CW complex of dimension $n$ there is a
  rather clear picture of complex vector bundles of rank $\geq [n/2]$
  (the so called stable rank) \cite{BP06}; much less is known about lower
   ranks and that is what makes it difficult to discard a Chern class as a top Chern class
    (besides, according to Thom \cite{Th54}, Theorem II.25, in a (compact, oriented) manifold
     any $a\in H^{2k}(M;\mathbb{Z})$ has a multiple which is a top Chern class).
     In any case, finding manifolds with certain cohomological properties
      would prove that Chern classes are more general than top Chern classes.
       For example, according to \cite{BP06} for a (compact, oriented) manifold $X$ of
       dimension $\leq$ 7, any $a\in H^4(X;\mathbb{Z})$ is the second Chern class
       of a rank 3 complex vector bundle. If it were the top Chern class of some $F$,
        then corollary 2.2 in \cite{BP06} applied to the direct sum of $F$ with the trivial line bundle would imply
\begin{equation}\label{eq:stsquare}
c_1(F)a+\mathrm{Sq}^2a\equiv 0\;\; \mathrm{in}\;\; H^6(X;\mathbb{Z}_2).
\end{equation}
Therefore, if $H^2(X;\mathbb{Z}_2)=0$ and there exist a class $a$ with non-vanishing
second Steenrod square, equation (\ref{eq:stsquare}) could not hold and hence $a$ would not be a top Chern class.

\end{document}